\documentclass[11pt]{article} 

\usepackage{latexsym, amsmath, 
amscd, amsthm, graphicx, epsfig, amssymb}
\usepackage{xypic}
\usepackage{psfrag}
\DeclareGraphicsExtensions{.ps}

\newtheorem{thm}{Theorem}[section]
\newtheorem{defn}[thm]{Definition}
\newtheorem{prop}[thm]{Proposition}
\newtheorem{cor}[thm]{Corollary}
\newtheorem{rema}[thm]{Remark}
\newtheorem{lemma}[thm]{Lemma}
\newtheorem{ass}[thm]{Assumption}

\newcommand{\halmos}{\rule{1ex}{1.4ex}}

\newcommand{\beq}{\begin{equation}}
\newcommand{\eeq}{\end{equation}}
\newcommand{\bnu}{\begin{enumerate}}
\newcommand{\enu}{\end{enumerate}}
 \newcommand{\bea}{\begin{eqnarray}}
\newcommand{\eea}{\end{eqnarray}}
 \newcommand{\nn}{\nonumber \\}
	
	\newcommand{\lbar}{\bigg\vert}

 \newcommand{\res}{\mbox{\rm Res}}
\renewcommand{\hom}{\mbox{\rm Hom}}

\newcommand{\edo}{\mbox{\rm End}\;}

 \newcommand{\pf}{{\it Proof.}\hspace{2ex}}
 \newcommand{\epf}{\hspace*{\fill}\mbox{$\halmos$}}

\newcommand{\wt}{\mbox{\rm wt}\ }

\newcommand{\lwt}{\mbox{\rm wt}^{L}\;}
\newcommand{\rwt}{\mbox{\rm wt}^{R}\;}

\newcommand{\id}{\mbox{\rm id}}
\newcommand{\one}{\mathbf{1}}

\newcommand{\C}{\mathbb{C}}
 
\newcommand{\HH}{\mathbb{H}}
\newcommand{\N}{\mathbb{N}}

\newcommand{\R}{\mathbb{R}}
\newcommand{\Z}{\mathbb{Z}}
\newcommand{\Y}{\mathcal{Y}}

\title{ {\bf Open-closed field algebras } }
\date{}
\author{Liang Kong}

\begin{document}

\setlength{\unitlength}{1cm}

\bibliographystyle{alpha}
\maketitle

\begin{abstract} 
We introduce the notions of open-closed field algebra
and open-closed field algebra over a vertex operator algebra $V$.
In the case that $V$ satisfies 
certain finiteness and reductivity conditions, we show that
an open-closed field algebra over $V$ canonically
gives an algebra over a $\C$-extension of Swiss-cheese partial operad. 
We also give a tensor-categorical formulation and constructions
of open-closed field algebras over $V$. 
\end{abstract}


\renewcommand{\theequation}{\thesection.\arabic{equation}}
\renewcommand{\thethm}{\thesection.\arabic{thm}}
\setcounter{equation}{0}
\setcounter{thm}{0}
\setcounter{section}{-1}

\section{Introduction}

Open-closed (or boundary) conformal field theory
was first developed in physics 
by Cardy \cite{C1}-\cite{C4}\cite{CL}. 
It has important applications to certain critical phenomena on 
surfaces with boundaries in condensed matter physics. 
It is also a powerful tool in the study of ``D-branes'' in 
string theory. A mathematical definition of such theory, 
which generalizes Segal's definition 
of conformal field theory \cite{S1}\cite{S2},
is given by Huang \cite{H7}, Hu and Kriz \cite{HKr}. 
However, a thorough mathematical study of this subject is still lacking.

In \cite{V}, Voronov introduced the so-called ``Swiss-cheese
operad'' to describe the interactions among 
open strings and closed strings. 
In order to incorporate the full conformal symmetry which is not
included in Voronov's notion of Swiss-cheese operad, 
Huang and the author introduced in \cite{HKo1} a notion called
Swiss-cheese partial operad. 
One of the main goals of this work is to construct examples of 
algebras over certain $\C$-extension of Swiss-cheese partial operad.

Swiss-cheese partial operad consists of disks 
with strips and tubes, which are conformal equivalence classes of 
disks with oriented punctures 
in the interior of the disks (called tubes) 
and those on the boundaries of the disks (called strips), 
together with a local coordinate map around each puncture.
On the moduli space of disks with strips and tubes, 
which is denoted as $\mathfrak{S}$,  
there are two kinds of sewing operations.  
One is called boundary sewing operation, which sews 
two elements in $\mathfrak{S}$ 
along two oppositely oriented strips. 
The other one is called interior sewing operation, which 
sews a disk with strips and tubes with a sphere with tubes 
\cite{H4} along two oppositely oriented tubes.

Disks without interior punctures are
closed under boundary sewing operations. This closed
structure is nothing
but the operad of disks with strips, 
which was introduced and studied in \cite{HKo1} and
denoted as $\Upsilon$. In \cite{HKo1}, Huang and the author
introduced the notion of open-string vertex operator algebra and 
showed that open-string vertex operator algebras 
of central charge $c$ are algebras over $\tilde{\Upsilon}^c$, which is
a $\C$-extension of $\Upsilon$.

The moduli space of spheres with tubes 
equipped with sewing operations
has a structure of partial operad, 
called sphere partial operad and denoted as $K$. 
It is well understood by the works of Huang 
\cite{H1}\cite{H2}\cite{H4}. 
Similar to the study of $\Upsilon$, 
we can study $\mathfrak{S}$
in the framework of sphere partial operad $K$
by embedding $\mathfrak{S}$ into $K$ 
via a doubling map $\delta$ \cite{A}\cite{C1}\cite{HKo1}.
The interior sewing operations on $\mathfrak{S}$
correspond to a double-sewing ``action'' of 
$K$ on the image of $\delta$ (see (\ref{double-equ})). 
As a consequence, the closed theory in an open-closed
conformal field theory must contain
both chiral part and anti-chiral part. 
Such closed theories were studied in \cite{HKo2}\cite{Ko1}\cite{HKo3}
in terms of the so-called (conformal) full field algebra and
variants of it. More precisely, we showed in \cite{Ko1} that
conformal full field algebras over $V^L\otimes V^R$, 
where $V^L$ and $V^R$ are vertex operator algebras
of central charge $c^L$ and $c^R$ respectively and 
satisfy certain finiteness and reductivity conditions, 
are algebras over the partial operad $\tilde{K}^{c^L}\otimes 
\overline{\tilde{K}^{\overline{c^R}}}$, where 
$\tilde{K}^{c^L}$ and $\overline{\tilde{K}^{\overline{c^R}}}$ 
are line bundles over the base space $K$ \cite{HKo1}\cite{Ko1}.

We are interested in algebras over $\tilde{\mathfrak{S}}^c$, 
which is a $\C$-extension of $\mathfrak{S}$. 
There is a natural action of $\tilde{K}^{c}\otimes 
\overline{\tilde{K}^{\overline{c}}}$ on $\tilde{\mathfrak{S}}^c$
induced by interior sewing operations \cite{HKo1}.
Therefore, it is natural
to expect that an algebra over $\tilde{\mathfrak{S}}^c$
with certain natural properties
should contain a conformal full field algebra $V_{cl}$
and an open-string vertex operator algebra $V_{op}$.
Moreover, $V_{cl}$ should ``act'' on $V_{op}$
according to the action of $\tilde{K}^{c}\otimes 
\overline{\tilde{K}^{\overline{c}}}$ on $\tilde{\mathfrak{S}}^c$.
This requires certain compatibility between $V_{cl}$ and $V_{op}$.
One of the compatibility conditions is the so-called
conformal invariant boundary condition. It roughly means that
the two vertex operator algebras generated by the 
left and right Virasoro elements in $V_{cl}$ 
must match in some way with that 
generated from the Virasoro element in $V_{op}$. In this case,
we also call it a boundary condition preserving conformal symmetry.  
In this work, we only study those boundary conditions
which preserve an enlarged symmetry given 
by a vertex operator algebra $U$ 
(see \cite{FS1}\cite{FS2} for general situations). 
Such structure is formalized by a notion called 
open-closed field algebra over $U$ (see Definition \ref{def-ocfa-U}). 
When $U$ satisfies conditions in Theorem \ref{ioa}, we 
show that an open-closed field algebra over $U$ canonically gives
an algebra over $\tilde{\mathfrak{S}}^c$.

Open-closed field algebras over $U$
are difficult to study for general $U$. 
The following Theorem proved by Huang in \cite{H8} 
is important for us to simplify the situation.  
\begin{thm} \label{ioa}
Let $V$ be a vertex operator algebra with central charge $c$ satisfying 
the following conditions: 
\bnu
\item Every $\C$-graded generalized $V$-module is a direct sum of 
$\C$-graded irreducible $V$-modules,
\item There are only finitely many inequivalent $\C$-graded
irreducible $V$-modules, 
\item Every $\R$-graded irreducible $V$-module satisfies the 
$C_1$-cofiniteness condition. 
\enu
Then the direct sum of all (in-equivalent) 
irreducible $V$-modules has
a natural structure of intertwining operator algebra and 
the category of $V$-modules, denoted as 
$\mathcal{C}_V$, has a natural structure of 
vertex tensor category. 
In particular, $\mathcal{C}_V$ 
has natural structure of braided tensor category. 
\end{thm}

\begin{ass}   {\rm
In this work, we fix a vertex operator algebra $V$, 
which is assumed to satisfy the conditions in
Theorem \ref{ioa} without further announcement. 
}
\end{ass}

The notion of open-closed field algebra
over $V$ can be described by very few data and
axioms in the framework of intertwining operator algebra
(see Theorem \ref{prop-op-cl-V} for a precise
statement). Moreover, such algebra has a very 
simple categorical formulation. We will introduce a tensor-categorical notion
called open-closed $\mathcal{C}_V|\mathcal{C}_{V\otimes V}$-algebra 
(see Definition \ref{w-cl-op-alg}) and
show that the category of open-closed field algebras over
$V$ is isomorphic to the category of
open-closed $\mathcal{C}_V|\mathcal{C}_{V\otimes V}$-algebras.
The categorical formulation of open-closed field algebra over $V$
makes it very easy to construct examples.

One of important developments on open-closed conformal field theories recently
is a series of works \cite{FFFS}\cite{FS3}\cite{FRS1}-\cite{FRS4}
\cite{FjFRS} by Felder, Fr\"{o}hlich, 
Fuchs, Runkel, Schweigert, Fjelstad, on
rational open-closed conformal field theories using 
the theories of modular tensor category and 
3-dimensional topological field theories. 
Assuming the existence of the structure of a modular tensor category 
on the category of modules for a vertex operator algebra and
the existence of conformal blocks with monodromies compatible with
the modular tensor category and all the necessary 
convergence properties, they constructed conformal blocks for
open-closed conformal field theories of all genus and 
proved their factorization properties and invariance properties
under the actions of mapping class groups. 
Our approach \cite{HKo1}-\cite{HKo3}\cite{Ko1}\cite{Ko2}
is somewhat complement to theirs.
We leave a detailed discussion of the relation between 
these two approaches to future publications.

The layout of this work is as follows.   
In Section 1, we review some old notions such as 
open-string vertex (operator) algebra, (conformal) full field algebra
and variants of them. We also introduce the notions of
open-closed field algebra and 
open-closed field algebra over $U$, and study their
basic properties.  
In Section 2, we give an operadic formulation of 
open-closed field algebra over $V$. 
In particular, we show that 
an open-closed field algebra over $V$ 
canonically gives an algebra over $\tilde{\mathfrak{S}}^c$. 
In Section 3, we give a categorical formulation 
of the notion of open-closed field algebra over $V$.

Convention of notations: $\N, \Z, \Z_+, \R, \R_+, \C$ 
denote the set of 
natural numbers, integers, positive integers, real numbers, 
positive real numbers, complex numbers, respectively. 
Let $\HH = \{ z\in \C | \text{Im} z >0\}$, 
$\overline{\HH} = \{ z\in \C | \text{Im} z < 0\}$, 
$\tilde{\HH} = \HH \cup \R$, 
$\tilde{\overline{\HH}} = \overline{\HH} \cup \R$. 
We use $\hat{\R}, \hat{\HH}, \hat{\overline{\HH}}$ to denote
the one point compactification of $\R$, 
$\tilde{\HH}$, $\tilde{\overline{\HH}}$ respectively. 
The following notations will also be used: $\forall n\in \Z_+$,
\bea
\Lambda^n &:=& \{ (r_1, \dots, r_n) \in \R^n| r_1 > \dots >r_n \geq 0\},  
\label{Lambda}  \\
M_{\HH}^{n} &:=& \{ (z_1, \dots, z_n) \in \HH^n| z_i\neq z_j, \,  \,
\mbox{for $i,j=1, \dots, n$ and $i\neq j$} \},  \label{M-HH} \\
M_{\C}^{n} &:=& \{ (z_1, \dots, z_n) \in \C^n | z_i\neq z_j, \,  \,
\mbox{for $i,j=1, \dots, n$ and $i\neq j$} \}.  \label{M-C}
\eea
The ground field is always assumed to be $\C$.

\paragraph{Acknowledgment}
This work grows from a chapter in author's thesis. 
I want to thank my advisor Yi-Zhi Huang 
for introducing this interesting field to me,  for his many valuable 
suggestions for improvement and for his constant support.
I also want to thank J. Lepowsky, J. Fuchs and C. Schweigert for 
many inspiring conversations, and I. Runkel for kindly
sending me his thesis and one of his computer programs.

\renewcommand{\theequation}{\thesection.\arabic{equation}}
\renewcommand{\thethm}{\thesection.\arabic{thm}}
\setcounter{equation}{0}
\setcounter{thm}{0}

\section{Open-closed field algebras}

Let $G$ be an abelian group. 
For any $G$-graded vector space $F=\oplus_{n\in G} F_{(n)}$ and any
$n\in G$, we shall use $P_{n}$ to denote the projection {}from 
$F$ or $\overline{F}=\prod_{n\in G}F_{(n)}$ to $F_{(n)}$.  We give 
$F$ and its
graded dual $F'=\oplus_{n\in G}F^{*}_{(n)}$ the topology induced
from the pairing between $F$ and $F'$. We also give $\hom(F,
\overline{F})$ the topology induced from the linear functionals on
$\hom(F, \overline{F})$ given by 
$f\mapsto \langle v', f(v)\rangle$ for
$f\in \hom(F, \overline{F})$, $v\in F$ and $v'\in F'$. 
For any vector space $F$ and any set $S$, 
$F^{\otimes 0}=\C$ and $S^{0}$ is an one-point set $\{ * \}$.

\subsection{Boundary field algebras and open-string vertex algebras}

\begin{defn} \label{bfa}
{\rm 
A {\it boundary field algebra} is a $\R$-graded vector space $V_{op}$ 
with grading operator being $\mathbf{d}_{op}$, 
together with a correlation-function map for each $n\in \N$:
\begin{eqnarray*}
m_{op}^{(n)}: V_{op}^{\otimes n}\times
\Lambda^{n}&\to& \overline{V}_{op}\\
(u_{1}\otimes \cdots \otimes u_{n}, (r_{1}, \dots, r_{n}))&\mapsto&
m_{op}^{(n)}(u_{1}, \dots, u_{n}; r_{1}, \dots, r_{n})
\end{eqnarray*}
and an operator $D_{op}\in \edo V_{op}$, 
satisfying the following axioms:
\begin{enumerate}

\item For each $n\in \Z_+$, 
$m_{op}^{(n)}(u_{1}, \dots, u_{n}; r_{1}, \dots, r_{n})$ is 
linear in $u_{1}, \dots, u_{n}$ and smooth in $r_{1}, \dots, r_{n}$.

\item $\forall u\in V_{op}$, $m_{op}^{(1)}(u; 0) = u$ and 
$\one_{op}:=m_{op}^{(0)}(1) \in (V_{op})_{(0)}$. 

\item {\it  Convergence property}:  
For $n\in \Z_+, k\in \N$ and $i=1, \dots, n$,  
$u_1, \dots, u_n, u_1^{(i)}, \ldots, u_{k}^{(i)} \in V_{op}$,  the following
series\footnote{That the number of
nonvanishing terms in the sum is countable is automatically assumed. }
\bea \label{bfa-conv-1}
&&\hspace{0cm} \sum_{s\in \R} 
m_{op}^{(n)}(u_1, \dots, u_{i-1}, P_s  m_{op}^k(u_1^{(1)}, 
\ldots, u_{k}^{(i)}, r_1^{(i)}, \ldots, r_{k}^{(i)}), \nn
&&\hspace{5cm} u_{i+1}, \dots, u_n; r_1, \ldots, r_{n})
\eea
converges absolutely, whenever  
$r_1^{(i)} < |r_i-r_j|$ for all $j\neq i$, to  
\bea  \label{bfa-conv-2}
&&m_{op}^{(n+k-1)}
(u_{1}, \dots, u_{i-1}, u_{1}^{(i)}, \dots, u_k^{(i)}, u_{i+1}, 
\dots, u_n; \nn
&&\hspace{2cm} r_1, \dots,  r_{i-1}, r_i+r_{1}^{(i)},
\dots, r_i+r_k^{(i)},  r_{i+1}, \dots,  r_{n}).
\eea

\item {\it $\mathbf{d}_{op}$-bracket property}: 
\beq
e^{a\mathbf{d}_{op}} m_{op}^{(n)}(u_1, \dots, u_n; r_1, \dots, r_n)
= m_{op}^{(n)}(e^{a\mathbf{d}_{op}} u_1, \dots, e^{a\mathbf{d}_{op}}u_n; 
e^ar_1, \dots, e^ar_n). \label{bfa-d-conj}
\eeq
for $n\in \Z_+, r_1>\dots >r_n\geq 0, r\in \R$ 
and $u_1\dots, u_n\in V_{op}$. 

\item {\it $D_{op}$-property}:
For $u_1, \dots, u_n \in V_{op}$, $r_1>\dots >r_n\geq 0$ 
and $r_n+a\geq 0$,
\beq  \label{bfa-D-conj}
e^{aD_{op}}m_{op}^{(n)}(u_{1}, \dots, u_{n}; r_{1}, \dots, r_{n}) 
=m_{op}^{(n)}(u_1, \dots, u_n, r_1+a, \dots, r_n+a).
\eeq

\end{enumerate}
}
\end{defn}

We denote such a boundary field algebra as
$(V_{op}, m_{op}, \mathbf{d}_{op}, D_{op})$.
{\it Homomorphisms}, {\it isomorphisms} and 
{\it subalgebras} of boundary field algebra 
are defined in the obvious way.

Let the map 
$Y_{op}: V_{op}^{\otimes (2)} \times \R_+ \rightarrow\overline{V}_{op}$ 
be defined by
\beq   \label{Y-op-bfa}
Y_{op}: \, \, (u\otimes v, r) \mapsto Y_{op}(u,r)v = m_2(u,v; r,0).
\eeq
Then by the convergence property, we have 
\bea 
Y_{op}(\one_{op}, r)&=& \id_{F} , \label{id-Y-op-bfa-1} \\
\lim_{r\to 0}Y_{op}(u, r) \one_{op} &=& u , \label{id-Y-op-bfa-2} 
\eea
for $u \in V_{op}$. (\ref{id-Y-op-bfa-2}) implies that
the map $u\mapsto Y_{op}(u, r)$ is one-to-one. 
By (\ref{bfa-d-conj}), we have 
\beq  \label{d-Y-op-bfa}
e^{a\mathbf{d}_{op}} Y_{op}(u, r) e^{-a\mathbf{d}_{op}} =
Y_{op}(e^{a\mathbf{d}_{op}} u, ar)
\eeq
for $u\in V_{op}, r\in \R_+$ and $a\in \R$.  
By (\ref{bfa-D-conj}), we also have 
\beq   \label{D-Y-op-bfa}
e^{aD_{op}} Y_{op}(u, r) e^{-aD_{op}}  = Y_{op}(u, r+a) 
\eeq
for $u\in V_{op}, r>0, r+a>0$. Moreover, using 
(\ref{id-Y-op-bfa-2}) and (\ref{bfa-D-conj}), we obtain the identity: 
\beq \label{s-creat-bfa}
Y_{op}(v, r)\one_{op} = m_{op}^{(1)}(v;r) = e^{rD_{op}} v. 
\eeq

If we assume some analytic properties on $m_{op}^{(n)}$, it is
possible to use $Y_{op}$ to generate all $m_{op}^{(n)}$. 
This motivate us to introduce the notion of 
open-string vertex algebra in \cite{HKo1}. 
The definition given here is a refinement of that in 
\cite{HKo1}.

\begin{defn} \label{op-va-defn}
{\rm An {\it open-string vertex algebra} 
is an $\R$-graded vector space $V_{op}=\oplus_{n\in \R}
(V_{op})_{(n)}$ (graded
by {\it weights}) equipped with a {\it vertex map}: 
\bea 
Y_{op}: (V_{op}\otimes V_{op}) \times \R_+ &\to& \overline{V}_{op} \nn
(u\otimes v, r)&\mapsto & Y_{op}(u, r)v,
\eea
a {\it vacuum} $\one_{op}\in V_{op}$
and an operator $D_{op}\in \edo V_{op}$ of weight $1$, satisfying 
the following conditions:

\bnu

\item {\it Vertex map weight property}: 
For $s_1, s_2\in \R$, 
there exists a finite subset $S(s_1,s_2)\subset \R$ such
that the image of 
$(\oplus_{s\in s_1+\Z} (V_{op})_{(s)}) \otimes 
(\oplus_{s\in s_2+\Z} (V_{op})_{(s)})$ under the vertex map $Y_{op}$ is in 
$\prod_{s\in S(s_1,s_2)+\Z} (V_{op})_{(s)}$. 

\item {\it Vacuum properties}: 
\bnu
\item {\it identity property}: For any $r\in \R_{+}$,
$Y_{op}(\one_{op},r)=\id_{V_{op}}$,
\item {\it creation property}: $\forall u\in V_{op}$, 
$\lim_{r\rightarrow 0} Y_{op}(u, r)\one_{op}=u$.
\enu


\item {\it Convergence properties}: 
\bnu
For $n\in \N$, $u_{1}, \dots, u_{n}, v\in V_{op}$ and $v'\in V'_{op}$, 
the series 
\bea
&&\langle v', Y_{op}(u_1, r_1)\dots Y_{op}(u_n, r_n) v\rangle  \nn
&&\hspace{0.2cm}=\sum_{m_1, \dots, m_{n-1}} \langle v', 
Y_{op}(u_1, r_1)P_{m_1}Y_{op}(u_2, r_2) \dots P_{m_{n-1}}Y_{op}(u_n, r_n)v\rangle
\nonumber
\eea
converges absolutely when $r_1>\dots >r_n>0$ and is a restriction
to the domain $\{ r_1>\dots >r_n>0 \}$ 
of an (possibly multivalued) analytic function in $(\C^{\times})^n$
with only possible singularities at $r_i=r_j$ for $1\leq i,j \leq n$
and $i\neq j$ .

\item For $u_1, u_2, v\in V_{op}, v'\in V'_{op}$, the series
\beq
\langle v', Y_{op}(Y_{op}(u_1, r_0)u_2, r_2)v\rangle 
= \sum_m \langle v', Y_{op}(P_m Y_{op}(u_1, r_0)u_2, r_2)v\rangle \nonumber
\eeq
converges absolutely when $r_2>r_0>0$. 
\enu

\item {\it Associativity}: For $u_{1}, u_{2}, v\in V_{op}$
and $v'\in V'_{op}$, 
$$\langle v', Y_{op}(u_{1}, r_1)Y_{op}(u_{2}, r_2)v\rangle =
\langle v', Y_{op}(Y_{op}(u_{1}, r_1-r_2)u_{2},r_2)v\rangle$$
for  $r_{1}, r_{2}\in \R$ satisfying 
$r_1>r_2>r_1-r_2>0$.

\item {\it $\mathbf{d}_{op}$-bracket property}: 
Let $\mathbf{d}_{op}$ be the grading 
operator on $V_{op}$.  For $u\in V_{op}$ and $r\in \R_{+}$,
\begin{equation} \label{d-bra}
[\mathbf{d}_{op}, Y_{op}(u,r)]= Y_{op}(\mathbf{d}_{op}u,r)+
r\frac{d}{dr}Y_{op}(u,r).
\end{equation}

\item {\it $D_{op}$-derivative property}: We still use $D_{op}$ to denote
the natural extension of $D_{op}$ to 
$\hom(\overline{V}_{op}, \overline{V}_{op})$.
For $u\in V_{op}$, 
$Y_{op}(u, r)$ as a map from $\R_{+}$ to 
$\hom(V_{op}, \overline{V}_{op})$
is differentiable and 
\begin{equation}\label{D-der}
\frac{d}{dr}Y_{op}(u,r)=[D_{op}, Y_{op}(u,r)]
= Y_{op}(D_{op}u,r).
\end{equation}

\enu

{\it Homomorphisms}, {\it isomorphisms}, {\it subalgebras}
of open-string vertex algebras are defined in the obvious way.  }
\end{defn}

We denote such algebra by 
$(V_{op}, Y_{op}, \one_{op}, \mathbf{d}_{op}, D_{op})$
or simply $V_{op}$.  For $u\in V_{op}$, 
it was shown in \cite{HKo1} that there is a formal vertex operator
\beq
Y_{op}^f(u, x) = \sum_{n\in \R} u_n x^{-n-1},     \label{formal-Y-osva-1}
\eeq
where $u_n\in \edo V_{op}$, so that 
\beq
Y_{op}^f(u, x)|_{x=r} = Y_{op}(u,r).    \label{formal-Y-osva-2}
\eeq
We can also replace $x$ by complex variable $z$ if we choose
a branch cut. For any $z\in \C^{\times}$ and $n\in \R$, we
define
\beq \label{branch-cut}
\log z := \log |z| + \text{arg} z, 0\leq \text{arg} z < 2\pi.
\eeq
But for power functions, we distinguish two types of complex
variables, $z$ (or $z_1, z_2, \dots$) and 
$\bar{z}_1, \zeta$ 
(or $\bar{z}_1, \zeta_1, \bar{z}_2, \zeta_2, \dots$). 
We define
\beq  \label{power-fun}
z^n: = e^{n\log z}, \quad \bar{z}^n := e^{n\overline{\log z}}, 
\quad  \zeta^n := e^{n\overline{\log \zeta}}. 
\eeq

\begin{prop}
An open-string vertex algebra canonically gives a boundary field algebra.
\end{prop}
\pf
The proof is standard. We omit it here.  
\epf

\begin{defn} 
{\rm An {\it open-string vertex operator
algebra} is an open-string vertex algebra, together
with a conformal element $\omega_{op}$,  satisfying the 
following conditions:

\bnu

\setcounter{enumi}{6}

\item {\it grading-restriction conditions}: For all $n\in \R$, 
$\dim (V_{op})_{(n)}<\infty$ and $(V_{op})_{(n)}=0$ when 
$n$ is sufficiently negative.

\item {\it Virasoro relations}: The vertex operator associated to
$\omega_{op}$ has the following expansion:
$$Y_{op}(\omega_{op}, r) =\sum_{n\in \Z} L(n) r^{-n-2}$$
where $L(n)$ satisfying the following condition: $\forall m,n\in \Z$,
$$[ L(m), L(n)]  = (m-n)L(m+n) - \frac{c}{12}(m^3-m)\delta_{m+n,0},$$
for some $c\in \C$.

\item {\it commutator formula for Virasoro operators and 
formal vertex operators (or component operators)}: For $v\in V_{op}$,
$Y_{op}^{f}(\omega_{op}, x)v$ 
involves only finitely many negative powers
of $x$ and
$$[Y_{op}^f(\omega_{op}, x_{1}), Y_{op}^f(v, x_{2})]
=\res_{x_0}x_{2}^{-1}\delta\left(\frac{x_{1}-x_{0}}{x_{2}}\right)
Y_{op}^f(Y_{op}^f(\omega_{op}, x_{0})v, x_{2}).$$

\item {\it $L(0)$-grading property} and 
{\it $L(-1)$-derivative property}:  $L(0)=\mathbf{d}_{op}$,
$L(-1)=D_{op}$.

\enu   
}
\end{defn}

We shall denote the open-string vertex operator algebra defined above
by $(V_{op}, Y_{op}, \one_{op}, \omega_{op})$ 
or simply $V_{op}$. The complex number $c$ 
in the definition is called {\it central charge}.

The {\it meromorphic center of $V_{op}$} is defined by
\begin{eqnarray*}
\lefteqn{C_{0}(V_{op})=\Bigg\{ u\in \oplus_{n\in \Z}(V_{op})_{(n)}
\;\Big|\; Y_{op}^{f}(u, x)\in (\edo V_{op})[[x, x^{-1}]],}\nn
&&\quad\quad\quad\quad\quad\quad\quad\quad\quad\quad\quad\quad
Y_{op}^{f}(v, x)u=e^{xD} Y_{op}^{f}(u, -x)v,\;\forall v\in V_{op} \Bigg\}.
\end{eqnarray*}
It was shown in \cite{HKo1} that the meromorphic center of 
a grading-restricted open-string vertex (operator) algebra is a
grading-restricted vertex (operator) algebra.

Let $U$ be a grading-restricted vertex (operator) algebra. 
If there is a monomorphism
$\iota_{op}: U\hookrightarrow C_0(V_{op})$ of 
grading-restricted vertex (operator) algebra, 
we call $V_{op}$ an open-string vertex (operator) algebra over $U$,
and denote it by $(V_{op}, Y_{op}, \iota_{op})$ or simply
by $V_{op}$. In this case, the formal vertex operator
$Y_{op}^f$ is an intertwining operator of type 
$\binom{V_{op}}{V_{op}V_{op}}$, where $V_{op}$ is a $U$-module.

\begin{rema} {\rm
Open-string vertex operator algebra can be viewed as a
noncommutative generalization of vertex operator algebra. 
Other noncommutative generalizations of vertex (operator) algebra
were also studied in the literature 
\cite{FR}\cite{B2}\cite{EK}\cite{BKa}\cite{L1}\cite{L2}\cite{L3}.
}
\end{rema}


\subsection{Full field algebras}

\begin{defn} \label{ffa}
{\rm A $\R\times \R$-graded 
{\it full field algebra} is an $\R\times \R$-graded vector space 
$V_{cl}=\oplus_{m,n\in\R}(V_{cl})_{(m,n)}$ 
(graded by {\it left weight} $\text{\rm wt}^L$ and 
{\it right weight} $\text{\rm wt}^R$ with left and right 
grading operators $\mathbf{d}^{L}$ and $\mathbf{d}^{R}$),   
equipped with  {\it correlation-function maps}
$$\begin{array}{rrcl}
m_{cl}^{(n)}: & V_{cl}^{\otimes n}\times
M_{\C}^n &\to& \overline{V}_{cl}\\
&(u_{1}\otimes \cdots \otimes u_{n}, (z_{1}, \dots, z_{n}))&\mapsto&
m_{cl}^{(n)}(u_{1}, \dots, u_{n}; z_{1}, \bar{z}_{1}, \dots, z_{n}, \bar{z}_n),
\end{array}$$
for each $n\in \N$, and operators $D^{L}$ and $D^{R}$ of weights
$(1,0)$ and $(0,1)$ respectively, 
satisfying the following axioms:

\begin{enumerate}

\item {\it Single-valuedness property}: 
$e^{2\pi i (\mathbf{d}^L-\mathbf{d}^R)} = \id_{V_{cl}}$. 

\item For $n\in \Z_+$, 
$m_{cl}^{(n)}(u_{1}, \dots, u_{n}; z_{1}, 
\bar{z}_{1}, \dots, z_{n}, \bar{z}_n)$ is 
linear in $u_{1}, \dots, u_{n}$ and smooth in the real and imaginary 
parts of $z_{1}, \dots, z_{n}$.

\item {\it Identity properties}: $\forall u\in V_{cl}$, 
$m_{1}(u; 0, 0)=u$ and 
$\one_{cl}:= m_{cl}^{(0)}(1) \in (V_{cl})_{(0,0)}$.


\item {\it  Convergence property}:   
For $k\in \Z_+$ and 
$u_1,\dots, u_n, u_{1}^{(i)}, \dots,
u_{k}^{(i)} \in V_{cl}$ and $i=1, \dots, n$,  the series 
\begin{eqnarray}   \label{ffa-conv-1}
&&\hspace{-1cm} \sum_{p,q\in \R} 
m_{cl}^{(n)}(u_1, \dots, u_{i-1}, P_{(p,q)} m_{cl}^{(k)} (v_1,\dots, v_k; 
z_1^{(i)}, \bar{z}_1^{(i)}, \dots, z_{k}^{(i)}, \bar{z}_{k}^{(i)}),  \nn
&&\hspace{3cm} 
u_{i+1}, \dots, u_n; z_1, \bar{z}_1, \dots, z_n, \bar{z}_n )
\end{eqnarray}
converges absolutely to  
\begin{eqnarray}\label{ffa-conv-2}
&&\hspace{-1cm} m_{cl}^{(n+k-1)}(u_1, \dots, u_{i-1}, v_1, \dots, v_k, 
u_{i+1}, \dots, u_n; z_1, \bar{z}_1, \dots, z_{i-1}, \bar{z}_{i-1}, \nn
&&\hspace{-0.5cm}z_i+z_{1}^{(i)}, \bar{z}_i + \bar{z}_{1}^{(i)}
\dots, z_i+z_{k}^{(i)},
\bar{z}_i+ \bar{z}_{k}^{(i)}, z_{i+1}, \bar{z}_{i+1}, \dots, 
z_n, \bar{z}_n)
\end{eqnarray}
whenever $|z_p^{(i)}| < |z_i-z_j|$ 
for all $j=1, \ldots, n$, $i\ne j$ and for $p=1, \dots, k$.

\item {\it Permutation property}: For $n\in \Z_{+}$ and 
$\sigma\in S_{n}$,
we have 
\begin{eqnarray}
\lefteqn{m_{cl}^{(n)}(u_{1}, \dots, u_{n}; z_{1}, \bar{z}_{1}, 
\dots, z_{n}, \bar{z}_{n})} \nn
&& = m_{cl}^{(n)}(u_{\sigma(1)}, \dots, u_{\sigma(n)}; 
z_{\sigma(1)}, \bar{z}_{\sigma(1)}, \dots, 
z_{\sigma(n)}, \bar{z}_{\sigma(n)})   \label{ffa-perm-axiom}
\end{eqnarray}
for $u_{1}, \dots, u_{n}\in V_{cl}$ and 
$(z_{1}, \dots, z_{n})\in M_{\C}^n$.

\item {\it $\mathbf{d}^L$ and $\mathbf{d}^R$ property}:
For $u_1, \dots, u_n \in V_{cl}$ and $a\in \C$, 
\bea
&&\hspace{-1cm} e^{a\mathbf{d}^L} e^{\bar{a}\mathbf{d}^R} 
m_{cl}^{(n)}(u_{1}, \dots, u_{n}; z_{1}, \bar{z}_{1}, 
\dots, z_{n}, \bar{z}_{n})  \nn
&&\hspace{-1cm} =m_{cl}^{(n)}(e^{a\mathbf{d}^L} e^{\bar{a}\mathbf{d}^R}u_1, \dots, 
e^{a\mathbf{d}^L} e^{\bar{a}\mathbf{d}^R}u_n, e^{a}z_1, e^{\bar{a}}\bar{z}_1, 
\dots, e^{a}z_n, e^{\bar{a}}\bar{z}_n).
\eea

\item {\it $D^{L}$ and $D^{R}$ property}: $[D^L, D^R]=0$ and
for $u_1, \dots, u_n \in V_{cl}$ and $a\in \C$, 
\bea  \label{m-cl-D}
&&e^{aD^L} e^{\bar{a}D^R} 
m_{cl}^{(n)}(u_{1}, \dots, u_{n}; z_{1}, \bar{z}_{1}, 
\dots, z_{n}, \bar{z}_{n})  \nn
&&\hspace{0.5cm} =m_{cl}^{(n)}(u_1, \dots, u_n, 
z_1+a, \bar{z}_1+\bar{a}, \dots, z_n+a, \bar{z}_n+\bar{a}).
\eea

\end{enumerate}
}

\end{defn}

We denote the $\R \times \R$-graded 
full field algebra defined above by 
$$(V_{cl}, m_{cl}, \mathbf{d}^L, \mathbf{d}^R, D^{L}, D^{R})$$ 
or simply by $V_{cl}$.

Subalgebra, homomorphism, monomorphism, epimorphism and isomorphism 
of full field algebra can be naturally defined.

Let $\mathbb{Y}: V_{cl}^{\otimes 2} \times \C^{\times} \rightarrow 
\overline{V}_{cl}$ be so that 
\beq  \label{Y-ffa}
\mathbb{Y}: (u\otimes v, z) \mapsto
\mathbb{Y}(u; z, \bar{z})v := m_2(u\otimes v; z, \bar{z}, 0,0).
\eeq
Then by the convergence property, it is easy to see that 
\bea 
\mathbb{Y}(\one_{cl}; z, \bar{z}) &=& \id_{V_{cl}}, \label{id-Y-ffa-1}\\
\lim_{z\rightarrow 0} \mathbb{Y}(u; z, \bar{z}) \one_{cl} &=& u, \quad\quad
\forall u\in V_{cl}. \label{id-Y-ffa-2}
\eea
Moreover, it is also not hard to show the following two properties of 
$\mathbb{Y}$:
\bnu
\item 
{\it The $\mathbf{d}^L$- and $\mathbf{d}^R$-bracket properties}: 
\begin{eqnarray}
\left[ \mathbf{d}^L, \mathbb{Y}(u; z,\bar{z})\right] 
&=& z\frac{\partial}{\partial z} \mathbb{Y}(u; z, \bar{z})
+\mathbb{Y}(\mathbf{d}^L u; z,\bar{z}),     \label{d-l-Y-ffa} \\
\left[ \mathbf{d}^R, \mathbb{Y}(u; z, \bar{z})\right]
&=& \bar{z}\frac{\partial}{\partial \bar{z}} \mathbb{Y}(u; z,\bar{z})
+ \mathbb{Y}(\mathbf{d}^R u; z, \bar{z}). \label{d-r-Y-ffa}
\end{eqnarray}

\item The {\it $D^{L}$- and $D^{R}$-derivative property}: 
\begin{eqnarray}
&&\left[ D^L, \mathbb{Y}(u; z,\bar{z})\right] = 
\mathbb{Y}(D^{L}u; z,\bar{z}) =
\frac{\partial}{\partial z}\mathbb{Y}(u; z, \bar{z}), 
\label{D-l-Y-ffa}\\
&&\left[ D^R, \mathbb{Y}(u; z, \bar{z})\right] = 
\mathbb{Y}(D^{R}u; z, \bar{z}) =
\frac{\partial}{\partial \bar{z}}\mathbb{Y}(u; z, \bar{z}). 
\label{D-r-Y-ffa}
\end{eqnarray}

\enu

It was shown in \cite{HKo2} that we have the
following expansion: 
\begin{equation}  \label{Y-exp-ffa}
\mathbb{Y}(u; z,\bar{z}) =\sum_{r, s\in \R} \mathbb{Y}_{l, r}(u)
e^{(-l-1)\log z} e^{(-r-1)\;\overline{\log
z}}
\end{equation}
where $\mathbb{Y}_{l, r}(u)\in \edo F$ with
$\text{\rm wt}^L \mathbb{Y}_{l, r}(u) = \text{\rm wt}^L u -l-1$ and 
$\text{\rm wt}^R \mathbb{Y}_{l, r}(u) = \text{\rm wt}^R u -r-1$.
Moreover, the expansion above is unique. 
Let $x$ and $\bar{x}$ be independent and commuting formal 
variables. We define the {\it formal full vertex operator}
 $\mathbb{Y}_f$ associated to $u\in V_{cl}$ by
\begin{equation} \label{formal-Y-ffa}
\mathbb{Y}_f(u; x,\bar{x}) = 
\sum_{l, r\in \R}\mathbb{Y}_{l, r}(u) x^{-l-1}\bar{x}^{-r-1}.
\end{equation}
For nonzero complex numbers $z$ and $\zeta$, 
we can substitute $e^{r\log z}$ and $e^{s\; \overline{\log \zeta}}$
for $x^{r}$ and $\bar{x}^{s}$ respectively 
in $\mathbb{Y}_f(u; x,\bar{x})$ to obtain an operator 
\beq
\mathbb{Y}_{\rm an}(u; z, \zeta): V_{cl} \times 
(\C^{\times})^2 \to \overline{V}_{cl}  \nonumber
\eeq
called {\it analytic full vertex operator}.

\begin{defn} \label{g-r-RR-g-ffa}
{\rm
An $\R \times \R$-graded full field algebra 
$(V_{cl}, m_{cl}, \mathbf{d}^L, \mathbf{d}^R, D^{L}, D^{R})$ 
is called {\it grading-restricted} 
if it satisfies the following grading-restriction conditions: 
\begin{enumerate}

\item There exists $M\in \R$ such that 
$(V_{cl})_{(m,n)}=0$ if $n<M$ or $m<M$.

\item $\dim (V_{cl})_{(m,n)}<\infty$ for $m, n\in \R$. \

\end{enumerate}
We  say that $V_{cl}$ is {\it lower-truncated} if $V_{cl}$ satisfies
the first grading restriction condition. }
\end{defn} 

In this case, for $u\in V_{cl}$ and $k\in \R$, we have  
$$
\sum_{l+r=k} \mathbb{Y}_{l, r}(u) \in \edo V_{cl}
$$ 
with total weight $\wt u -k-2$. 
We denote $\sum_{l+r=k}  \mathbb{Y}_{l, r}(u) $ 
by $\mathbb{Y}_{k-1}(u)$. Then we have the  expansion
\begin{equation}  \label{Y-x-x}
\mathbb{Y}_f(u; x, x) = 
\sum_{k\in \R} \mathbb{Y}_k(u) x^{-k-1},
\end{equation}
where $\wt \mathbb{Y}_k(u) = \wt u - k -1$. For given $u, v\in v_{cl}$,
 we have 
$\mathbb{Y}_k(u) w=0$ for sufficiently large $k$.

Let $(V^{L}, Y^{L}, \one^{L}, \omega^{L})$
and $(V^{R}, Y^{R}, \one^{R}, \omega^{R})$ be vertex operator
algebras. It was pointed out in \cite{HKo2} that 
$V^L\otimes V^R$ has a natural structure of full field algebra. 
Let $\iota_{cl}$ be an injective homomorphism from
the full field algebra $V^{L}\otimes V^{R}$ to $V_{cl}$. 
Then we have $\one_{cl}=\iota_{cl} (\one^{L}\otimes \one^{R})$,
$\mathbf{d}^{L}\circ \iota_{cl} =\iota_{cl}\circ (L^{L}(0)\otimes I_{V^{R}})$,
$\mathbf{d}^{R}\circ \iota_{cl}=
\iota_{cl} \circ (I_{V^{L}}\otimes  L^{R}(0)))$,
$D^{L}\circ \iota_{cl}=\iota_{cl}\circ (L^{L}(-1)\otimes I_{V^{R}})$ and
$D^{R}\circ \iota_{cl}=\iota_{cl}\circ (I_{V^{L}}\otimes  L^{R}(-1))$.
Moreover, $V_{cl}$ has a {\it left conformal element} 
$\iota_{cl}(\omega^{L}\otimes \one^{R})$ and an {\it right 
conformal element} $\iota_{cl}(\one^{L}\otimes \omega^{R})$. 
We have the following operators on $V_{cl}$: 
\begin{eqnarray*}
L^{L}(0)&=&\res_{x}\res_{\bar{x}}\bar{x}^{-1}\mathbb{Y}_{f}(
\iota_{cl}(\omega^{L}\otimes
\mathbf{1}^{R}); x, \bar{x}),\\
L^{R}(0)&=&\res_{x}\res_{\bar{x}}x^{-1}\mathbb{Y}_{f}(
\iota_{cl}(\mathbf{1}^{L}\otimes 
\omega^{L}); x, \bar{x}),\\
L^{L}(-1)&=&\res_{x}\res_{\bar{x}}x\bar{x}^{-1}\mathbb{Y}_{f}(
\iota_{cl}(\omega^{L}\otimes
\mathbf{1}^{R}); x, \bar{x}),\\
L^{R}(-1)&=&\res_{x}\res_{\bar{x}}x^{-1}\bar{x}\mathbb{Y}_{f}(
\iota_{cl}(\mathbf{1}^{L}\otimes
\omega^{L}); x, \bar{x}).
\end{eqnarray*}
Since these
operators are operators on $V_{cl}$, it should be easy to distinguish 
them from those operators with the same notation but 
acting on $V^{L}$ or $V^{R}$.

\begin{defn} 
{\rm Let $(V^{L}, Y^{L}, \one^{L}, \omega^{L})$ 
and $(V^{R}, Y^{R}, \one^{R}, \omega^{R})$ be vertex operator algebras. 
A {\it full field algebra over $V^{L}\otimes V^{R}$} 
is a grading-restricted $\R\times \R$-graded full field 
algebra $(V_{cl}, m_{cl}, \mathbf{d}^L,\mathbf{d}^R,
D^{L}, D^{R})$ equipped with an injective homomorphism $\iota_{cl}$ from 
the full field algebra $V^{L}\otimes V^{R}$ to $V_{cl}$ such that 
$\mathbf{d}^{L}=L^{L}(0)$, $\mathbf{d}^{R}=L^{R}(0)$, 
$D^{L}=L^{L}(-1)$ and $D^{R}=L^{R}(-1)$. A full field algebra over
$V^L\otimes V^R$ equipped with left and right conformal 
elements $\iota_{cl}(\omega^L\otimes \one^R)$ and 
$\iota_{cl}(\one^L \otimes \omega^R)$ is called {\em 
conformal full field algebra over $V^L\otimes V^R$}. 
}
\end{defn}

We shall denote the (conformal) full field algebra over 
$V^{L}\otimes V^{R}$
defined above by $(V_{cl}, m_{cl}, \iota_{cl})$ or simply by $V_{cl}$.
From now on, we will not distinguish $V^L$ with 
$V^L\otimes \one^R$ and $\iota_{cl}(V^L\otimes \one^R)$. Similarly for
$V^R$ with $\one^L \otimes V^R$ and $\iota_{cl}(\one^L \otimes V^R)$.
For $u^L\in V^L$, $\mathbb{Y}_{an}(u^L; z, \zeta)$ is
independent of $\zeta$. So we simply denote it as 
$\mathbb{Y}_{an}(u^L, z)$. Similarly, we denote 
$\mathbb{Y}_{an}(u^R; z, \zeta)$ as $\mathbb{Y}_{an}(u^R, \zeta)$
for $u^R\in V^R$.

\subsection{Open-closed field algebras}

\begin{defn}{\rm
An open-closed field algebra consists of a $\R\times \R$-graded full field algebra 
$$(V_{cl}, m_{cl}, \mathbf{d}_{cl}^L, \mathbf{d}_{cl}^R, D_{cl}^L, D_{cl}^R)$$
and a $\R$-graded vector space $V_{op}$ 
with the grading operator $\mathbf{d}_{op}$ and an additional operator 
$D_{op}\in \edo V_{op}$, 
together with a map for each pair of $n, l\in \N$: 
\bea
m_{cl-op}^{(l;n)}: V_{cl}^{\otimes l} \otimes V_{op}^{\otimes n} \times 
M_{\HH}^{l} \times \Lambda^{n}  &\rightarrow& \overline{V}_{op}   \nn
(u_1\otimes \dots \otimes u_l\otimes 
v_1 \otimes \dots \otimes v_n, (z_1, \dots, z_l ; r_1, \dots, r_n)) 
&\mapsto& \nn
&&\hspace{-5cm}
m_{cl-op}^{(l;n)}(u_1, \dots, u_l; v_1, \dots, v_n; 
z_1, \bar{z}_1 \dots, z_l, \bar{z}_l; r_1, \dots, r_n),
\nonumber
\eea
satisfying the following axioms:
\bnu
\item $m_{cl-op}^{(l;n)}
(u_1, \dots, v_n; 
z_1, \bar{z}_1, \dots, z_l, \bar{z}_l;r_1, \dots, r_n)$
is linear in $u_1, \dots, v_n$ and smooth in 
$r_1, \dots, r_n, z_1, \dots, z_l$. 

\item {\it Identity properties}: 
$m_{cl-op}^{(0;1)}(v;0) = v, \forall v\in V_{op}$ and 
$\one_{op}:=m_{cl-op}^{(0;0)}(1)\in (V_{op})_{(0)}$.

\item {\it Convergence properties}: 
\bnu
\item 
For $u_1, \dots, u_l, \tilde{u}_1, \dots, \tilde{u}_k \in V_{cl}$,
$v_1, \dots, v_n, \tilde{v}_1, \dots, \tilde{v}_m \in V_{op}$ 
and $i=1, \dots, n$,  the following series
\bea  \label{ocfa-conv-1}
&&\sum_{n_1\in \R} m_{cl-op}^{(l;n)}(u_1, \dots, u_l; v_1, \dots, v_{i-1},
P_{n_1} m^{(k;m)}_{cl-op}(\tilde{u}_1, \dots, \tilde{u}_{k}; \nn
&&\hspace{2cm} 
\tilde{v}_1, \dots, \tilde{v}_{m};  z_1^{(i)}, \overline{z_1^{(i)}}, 
\dots, z_{k}^{(i)}, \overline{z_k^{(i)}}; r_1^{(i)}, \dots, r_m^{(i)}); \nn
&&\hspace{1.5cm} v_{i+1}, \dots v_n;  
z_1, \bar{z}_1, \dots, z_l, \bar{z}_l; r_1, \dots, r_n)   
\eea
converges absolutely when 
$|z_s-r_i|, |r_t-r_i|>|z_p^{(i)}|, r_q^{(i)}\geq 0$ 
for all $s=1, \dots, l$,
$t=1, \dots, n$, $t\neq i$, $p=1, \dots, k$ and $q=1, \dots, m$ to 
\bea \label{ocfa-conv-2}
&&\hspace{-0.5cm}
m^{(l+k;n+m-1)}_{cl-op}(u_1, \dots, u_l, \tilde{u}_1, \dots,  \tilde{u}_k; 
v_1, \dots, v_{i-1}, \tilde{v}_1, \dots, \tilde{v}_m,  \nn
&&\hspace{0.5cm} 
v_{i+1}, \dots, v_n; 
z_1, \bar{z}_1, \dots, z_l, \bar{z}_l, 
z_1^{(i)}, \overline{z_1^{(i)}}, \dots, z_{k}^{(i)}, 
\overline{z_{k}^{(i)}};   \nn
&&\hspace{2cm} r_1, \dots, r_{i-1}, r_i+r_1^{(i)}, \dots, 
r_i+r_m^{(i)}, r_{i+1}, \dots, r_n). 
\eea

\item 
For $u_1, \dots, u_l, \tilde{u}_1, \dots, \tilde{u}_k \in V_{cl}$,
$v_1, \dots, v_n \in V_{op}$ 
and $i=1, \dots, n$,  the following series
\bea  \label{ocfa-conv-3}
&&\hspace{0cm}\sum_{n_1, n_2\in \R} m_{cl-op}^{(l;n)}( 
u_1, \dots, u_{i-1}, 
P_{(n_1, n_2)} m_{cl}^{(k)} (\tilde{u}_1, \dots, \tilde{u}_{k};
 z_1^{(i)}, \overline{z_1^{(i)}}, \dots, z_{k}^{(i)}, 
\overline{z_{k}^{(i)}}),  \nn
&&\hspace{2cm} u_{i+1}, \dots, u_{l}; v_{1}, \dots v_n;  
z_1, \bar{z}_1, \dots, z_l, \bar{z}_l; r_1, \dots, r_n)   
\eea
converges absolutely, when $|z_s-z_i|, |r_t-z_i|>|z_p^{(i)}|$
for all $s=1, \dots, l$, $s\neq i$, $t=1, \dots, n$, 
$p=1, \dots, k$, to 
\bea  \label{ocfa-conv-4}
&&\hspace{-1cm}
m_{cl-op}^{(l+k-1;n)}(
u_1, \dots, u_{i-1}, \tilde{u}_1, \dots, \tilde{u}_{k}, 
u_{i+1}, \dots, u_l;  v_{1}, \dots v_n; \nn
&&\hspace{1cm}
 z_1, \bar{z}_1, \dots, z_{i-1}, \bar{z}_{i-1}, 
 z_i+z_1^{(i)}, \overline{z_i+z_1^{(i)}}, \dots,  \nn
&&\hspace{2cm} z_i +z_k^{(i)}, \overline{z_i +z_k^{(i)}},
z_{i+1}, \bar{z}_{i+1}, \dots, z_l, \bar{z}_l; r_1, \dots, r_n). 
\eea
\enu

\item {\it Permutation axiom}: For $u_1, \dots, u_l\in V_{cl}$, 
$v_1, \dots, v_n\in V_{op}$ and $\sigma \in S_l$, 
\bea
&&\hspace{0cm}m_{cl-op}^{(l;n)} ( u_1, \dots, u_l;  v_1, \dots, v_n;
 z_1, \bar{z}_1, \dots, z_l, \bar{z}_l; 
r_1, \dots, r_n) \nn
&&\hspace{1cm}
= m_{cl-op}^{(l;n)} (u_{\sigma(1)}, \dots, u_{\sigma(l)}; v_1, \dots, v_n; \nn 
&&\hspace{3cm} 
z_{\sigma(1)}, \bar{z}_{\sigma(1)}, \dots, z_{\sigma(l)}, \bar{z}_{\sigma(l)};
r_1, \dots, r_n). 
\eea

\item {\it $\mathbf{d}_{op}-$, $\mathbf{d}_{cl}^L-$ 
and $\mathbf{d}_{cl}^R$-property}: 
For $u_1, \dots, u_l\in V_{cl}$, $v_1, \dots, v_n\in V_{op}$
and $a\in \R$, 
\bea   \label{d-m-co}
&&\hspace{-1cm}
e^{a\mathbf{d}_{op}} m_{cl-op}^{(l;n)} ( u_1, \dots, u_l; v_1, \dots, v_n; 
 z_1, \bar{z}_1, \dots, z_l, \bar{z}_l; r_1, \dots, r_n;) \nn
&&\hspace{-0.5cm}
= m_{cl-op}^{(l;n)} (e^{a(\mathbf{d}_{cl}^L + \mathbf{d}_{cl}^R)} u_1, \dots, 
e^{a(\mathbf{d}_{cl}^L + \mathbf{d}_{cl}^R)}u_l; 
e^{a\mathbf{d}_{op}}v_1, \dots, e^{a\mathbf{d}_{op}}v_n; \nn
&&\hspace{3cm} e^a z_1, e^a \bar{z}_1, \dots, e^az_l, e^a\bar{z}_l; 
e^ar_1, \dots, e^ar_n).
\eea

\item {\it $D_{op}$-property}: 
For $u_1, \dots, u_l\in V_{cl}$, $v_1, \dots, v_n\in V_{op}$ and
$r_n+a \geq 0$,
\bea \label{D-m-co}
&&\hspace{-1cm}
e^{aD_{op}} m_{cl-op}^{(l;n)} ( u_1, \dots, u_l; v_1, \dots, v_n;
 z_1, \bar{z}_1, \dots, z_l, \bar{z}_l; r_1, \dots, r_n) \nn
&&\hspace{-0.5cm}
= m_{cl-op}^{(l;n)} ( u_1, \dots, u_l;  v_1, \dots, v_n;\nn
&&\hspace{0.5cm}  
z_1+a, \bar{z}_1+a, \dots, z_l+a, \bar{z}_l+a;
r_1+a, \dots, r_n+a).
\eea

\enu

We denote such algebra 
by $(V_{cl}, V_{op}, m_{cl-op})$ for simplicity. 
{\it Homomorphisms}, {\it isomorphisms}, {\it subalgebras}
of open-closed field algebras are defined in the obvious way. 
}
\end{defn}

\begin{rema}  {\rm 
From above definition, it is clear that an open-closed field
algebra automatically includes a boundary field algebra as a substructure. }
\end{rema}

We discuss a few results which follow immediately from the 
definition. By the identity properties and (\ref{D-m-co}), we
have, for $a\geq 0$,
\beq  \label{m-1-0-D-op}
m_{cl-op}^{(0;1)}(v; a) = e^{aD_{op}} m_{cl-op}^{(0;1)}(v;0) = e^{aD_{op}}v. 
\eeq
By (\ref{m-1-0-D-op}) and the convergence property, 
for $i=1, \dots, n$ and $a \in \R$, we obtain 
\bea  \label{m-D-op-it}
&&m_{cl-op}^{(l;n)} ( u_1, \dots, u_l; v_1, \dots, v_{i-1}, 
e^{aD_{op}} v_i, v_{i+1}, \dots,  v_n;   \nn 
&&\hspace{5cm} 
z_1, \bar{z}_1, \dots, z_l, \bar{z}_l; r_1, \dots, r_n) \nn
&&\hspace{1cm}= m_{cl-op}^{(l;n)} (u_1, \dots, u_l; v_1, \dots, v_{i-1}, 
m_{cl-op}^{(0;1)}(v_i;a),  \nn
&&\hspace{4cm}   v_{i+1}, \dots,  v_n; 
z_1, \bar{z}_1, \dots, z_l, \bar{z}_l; r_1, \dots, r_n) \nn
&&\hspace{1cm}
=m_{cl-op}^{(l;n)} (u_1, \dots, u_l; v_1, \dots, v_n; 
 z_1, \bar{z}_1, \dots, z_l, \bar{z}_l  \nn
&&\hspace{5cm} 
r_1, \dots, r_{i-1}, r_i+a, r_{i+1}, \dots, r_n)
\eea
when $|r_j-r_i|,|z_k-r_i| > |a|$ for $j\neq i$ and $k=1,\dots, l$. 
By (\ref{m-cl-D}) and the convergence property, 
we also have, for $j=1, \dots, l$ and $b\in \C$,
\bea  \label{m-D-cl-it}
&&m_{cl-op}^{(l;n)} ( u_1, \dots, u_{j-1}, 
e^{bD_{cl}^L+\bar{b}D_{cl}^R}u_j, u_{j+1}, \dots, u_l; v_1, \dots,  v_n;\nn 
&&\hspace{6cm}
z_1, \bar{z}_1, \dots, z_l, \bar{z}_l;  r_1, \dots, r_n) \nn
&&\hspace{1cm} =
m_{cl-op}^{(l;n)} (
u_1, \dots, u_{j-1}, m_{cl}^{(1)}(u_j;b,\bar{b}), u_{j+1}, \dots, u_l; \nn 
&&\hspace{4cm}   v_1, \dots,  v_n;  
z_1, \bar{z}_1, \dots, z_l, \bar{z}_l; r_1, \dots, r_n) \nn
&&\hspace{1cm}= 
m_{cl-op}^{(l;n)} ( u_1, \dots, u_l; v_1, \dots, v_n; 
z_1, \bar{z}_1, \dots, z_{j-1}, \bar{z}_{j-1},\nn
&&\hspace{3cm}  
z_j+b, \overline{z_j+b}, z_{j+1}, \bar{z}_{j+1}, 
\dots, z_l, \bar{z}_l; r_1, \dots, r_n)
\eea
when
$|z_i-z_j|, |r_k-z_j|>|b|$ for $i=1, \dots, l, i\neq j$ 
and $k=1,\dots,n$.

Let $m_{op}^{(n)}:= m_{cl-op}^{(0;n)}$. The definition of 
open-closed field algebra immediately implies that 
$(V_{op}, m_{op}, \mathbf{d}_{op}, D_{op})$ is a boundary field algebra.
In particular, the map
$Y_{op}$ defined in (\ref{Y-op-bfa}) satisfies the equations 
(\ref{id-Y-op-bfa-1}), (\ref{id-Y-op-bfa-2}), 
(\ref{d-Y-op-bfa}) and (\ref{D-Y-op-bfa}). 
Similarly, we know that the map $\mathbb{Y}$ defined in (\ref{Y-ffa}) satisfies
the equations
(\ref{id-Y-ffa-1}), (\ref{id-Y-ffa-2}), 
(\ref{d-l-Y-ffa}), (\ref{d-r-Y-ffa}),
(\ref{D-l-Y-ffa}), (\ref{D-r-Y-ffa}), (\ref{Y-exp-ffa}) and
(\ref{formal-Y-ffa}). 

In an open-closed field algebra $(V_{cl}, V_{op}, m_{cl-op})$, 
there is an additional vertex operator map:  
\begin{eqnarray*}
\mathbb{Y}_{cl-op}: (V_{cl}\otimes V_{op}) \times \HH
&\to& \overline{V}_{op} \\
(u \otimes v, (z, \bar{z})) &\mapsto &
\mathbb{Y}_{cl-op}(u; z, \bar{z})v, 
\end{eqnarray*}
defined by 
\beq   \label{Y-co-ocfa}
\mathbb{Y}_{cl-op}(u; z, \bar{z})v := m_{cl-op}^{(1;1)}(u;v;z,\bar{z};0).
\eeq

By the convergence property, we have the following 
identity property:
\bea
\mathbb{Y}_{cl-op}(\one_{cl}; z, \bar{z})v &=&
m_{cl-op}^{(1;1)}(m_{cl}^{(0)}(1);v;z,\bar{z};0) \nn
&=& m_{cl-op}^{(0;1)}(v;0) \nn
&=& v.
\eea

\begin{prop}
For $u\in V_{cl}$, we have
\bea  \label{d-op-brack-0}
[\mathbf{d}_{op}, \mathbb{Y}_{cl-op}(u; z, \bar{z})] &=& 
\mathbb{Y}_{cl-op}((\mathbf{d}_{cl}^L + \mathbf{d}_{cl}^R)u; z, \bar{z}) 
  \nn
&&\hspace{0.5cm} + \left( z\frac{\partial}{\partial z}
+ \bar{z} \frac{\partial}{\partial \bar{z}} \right) 
\mathbb{Y}_{cl-op}(u; z, \bar{z}).
\eea
\end{prop}
\pf
Applying $\frac{\partial}{\partial a}|_{a=0}$ to both sides of 
(\ref{d-m-co}) when $n=1,l=1$ and $r_1=0$, we obtain 
(\ref{d-op-brack-0}) immediately. 
\epf

\begin{prop}
For $u\in V_{cl}$, we have
\bea  
[D_{op}, \mathbb{Y}_{cl-op}(u; z, \bar{z})] &=& 
\mathbb{Y}_{cl-op}((D_{cl}^L + D_{cl}^R)u; z, \bar{z}), 
\label{D-op-brack-1}  \\
\mathbb{Y}_{cl-op}(D_{cl}^Lu; z, \bar{z})
&=& \frac{\partial}{\partial z}
\mathbb{Y}_{cl-op}(u; z, \bar{z}), \label{D-op-brack-2} \\
\mathbb{Y}_{cl-op}(D_{cl}^Ru; z, \bar{z})
&=& \frac{\partial}{\partial \bar{z}}
\mathbb{Y}_{cl-op}(u; z, \bar{z}). \label{D-op-brack-3}
\eea
\end{prop}
\pf
Applying $\frac{\partial}{\partial b}|_{b=0}$ 
($\frac{\partial}{\partial \bar{b}}|_{\bar{b}=0}$) to both sides of 
(\ref{m-D-cl-it}) when $n=1,l=1$ and $r_1=0$, we obtain 
(\ref{D-op-brack-2}) and (\ref{D-op-brack-3}).

Applying $\frac{\partial}{\partial a}|_{a=0}$ to both sides of 
(\ref{D-m-co}) when $n=1,l=1$ and $r_1=0$ and using 
(\ref{D-op-brack-2}) and (\ref{D-op-brack-3})
we obtain the first identity in (\ref{D-op-brack-1}) immediately. 
\epf

\subsection{Analytic open-closed field algebras}

The notion of open-closed field algebra introduced in the last
subsection is very general. 
There are not much to say about open-closed field algebras 
in such generality. In this subsection, we study those
open-closed field algebras satisfying some nice analytic properties. 
In this case, the whole structure can be reconstructed from 
some simple ingredients.

\begin{defn} \label{an-ext-ocfa-def}
{\rm
An open-closed field algebra $(V_{cl}, V_{op}, m_{cl-op})$ 
is called {\it analytic} if it satisfies the following 
conditions:
\bnu

\item $\mathbb{Y}_{cl-op}$ can be extended to a map
$V_{cl}\otimes V_{op} \times \HH \times \overline{\HH} \rightarrow
\overline{V}_{op}$ such that for $z\in \HH, \zeta\in \overline{\HH}$
\beq
\mathbb{Y}_{cl-op}(u; z, \bar{z}) = \mathbb{Y}_{cl-op}(u;z,\zeta) 
|_{\zeta=\bar{z}}. 
\eeq

\item For $n\in \N$, $v_1,\cdots, v_{n+1}\in V_{op}, v'\in (V_{op})'$ 
and $u_1,\cdots, u_n\in V_{cl}$, the series
$$
\langle v', \mathbb{Y}_{cl-op}(u_1, z_1,\zeta_1)Y_{op}(v_1, r_1) \cdots 
\mathbb{Y}_{cl-op}(u_n, z_n,\zeta_n)Y_{op}(v_n, r_n)v_{n+1}\rangle
$$
is absolutely convergent when 
$|z_1|, |\zeta_1| > r_1> \cdots >|z_n|, |\zeta_n|>r_n>0$ and 
can be extended to a (possibly multivalued) analytic function on 
$$
\{ (z_1, \zeta_1, r_1, \dots, z_n, \zeta_n, r_n) \in M_{\C}^{3n} \}.
$$ 

\item For $n\in \N$, $v', u_1, \dots, u_{n+1}\in V_{cl}$, 
the series
$$
\langle v', \mathbb{Y}_{an}(u_1; z_1, \zeta_1)\dots 
\mathbb{Y}_{an}(u_n; z_n, \zeta_n)u_{n+1}\rangle 
$$
is absolutely convergent when 
$|z_1|>\dots >|z_{n}|>0$ and $|\zeta_1|>\dots >|\zeta_n|>0$
and can be extended to an analytic function on $M_{\C}^{2n}$.

\item For $v'\in V_{op}', v_1, v_2\in V_{op}, u\in V_{cl}$, the series
$$
\langle v', Y_{op}(\mathbb{Y}_{cl-op}(u; z, \zeta)v_1, r)v_2\rangle
$$
is absolutely convergent when $r>|z|, |\zeta|>0$.

\item For $v'\in V_{op}', v\in V_{op}, u_1, u_2\in V_{cl}$, the series
$$
\langle v', \mathbb{Y}_{cl-op}(\mathbb{Y}_{an}(u_1; z_1, \zeta_1)u_1; 
z_2, \zeta_2)v\rangle,
$$
converges absolutely when $|z_2|>|z_1|>0$, $|\zeta_2|>|\zeta_1|>0$
and $|z_1|+|\zeta_1|<|z_2-\zeta_2|$. 

\enu
}
\end{defn}

By the convergence properties of open-closed field algebra, 
$m_{cl-op}^{(l;n)}$ can be expressed as products of $\mathbb{Y}_{cl-op}, 
Y_{op}$ on a dense subdomain of $M_{\mathbb{H}}^l\times \Lambda^n$. 
$m_{cl-op}^{(l;n)}$ on the
complement of this dense subdomain is completely determined by
analytic extension. Therefore, $m_{cl-op}^{(l;n)}, l,n\in \N$ 
are completely determined by $\one_{op}$, $Y_{op}$ 
and $\mathbb{Y}_{cl-op}$. Similarly, $m_{cl}^{(n)}, n\in \N$ is 
completely determined by $\one_{cl}$ and $\mathbb{Y}_{an}$. 
Therefore, we also denote an analytic open-closed conformal 
field algebra by
$$
((V_{cl}, \mathbb{Y}_{an}, \one_{cl}), (V_{op}, Y_{op}, \one_{op}), 
\mathbb{Y}_{cl-op}). 
$$
or $(V_{cl}, V_{op}, \mathbb{Y}_{cl-op})$ for simplicity.

Two immediate consequences of the 
definition of analytic open-closed field algebra are given in 
the following two Lemmas.
\begin{lemma}   \label{lemma-D-brack}
{\em $D_{op}$-bracket properties}: $u\in V_{cl}$, we have
\bea  \label{D-op-brack}
[D_{op}, \mathbb{Y}_{cl-op}(u; z, \zeta)] &=& 
\mathbb{Y}_{cl-op}((D_{cl}^L + D_{cl}^R)u; z, \zeta), \nn
\mathbb{Y}_{cl-op}(D_{cl}^Lu; z, \zeta)
&=& \frac{\partial}{\partial z}
\mathbb{Y}_{cl-op}(u; z, \zeta), \nn
\mathbb{Y}_{cl-op}(D_{cl}^Ru; z, \zeta)
&=& \frac{\partial}{\partial \zeta}
\mathbb{Y}_{cl-op}(u; z, \zeta).
\eea
\end{lemma}
\pf
By (\ref{D-op-brack-1}),(\ref{D-op-brack-2}) and (\ref{D-op-brack-3}), 
for any fixed $z\in \HH$, 
all three identities holds when $\zeta=\bar{z}$. Replacing
$u$ by $(D_{cl}^R)^ku$ for any $k\in \N$ in 
(\ref{D-op-brack-1}), (\ref{D-op-brack-2}) and (\ref{D-op-brack-3}),  
we obtained that
all derivatives $(\frac{\partial}{\partial \zeta})^k |_{\zeta=\bar{z}}, 
\forall k\in \N$
of both sides of above three identities equal respectively. 
By the property of analytic function
and the fact that $\overline{\HH}$ is a simply connected domain, 
we obtain that both sides of above three identities, viewed as 
analytic functions for a fixed $z\in \HH$ and any $\zeta\in \overline{\HH}$,
must equal respectively.  
\epf
\begin{lemma}
{\it $\mathbf{d}_{op}$-bracket properties }:
$u\in V_{cl}$, we have 
\bea  \label{d-op-brack}
[\mathbf{d}_{op}, \mathbb{Y}_{cl-op}(u; z, \zeta)] &=& 
\mathbb{Y}_{cl-op}((\mathbf{d}_{cl}^L + \mathbf{d}_{cl}^R)u; z, \zeta) 
  \nn
&&\hspace{0.5cm} + \left( z\frac{\partial}{\partial z}
+ \zeta \frac{\partial}{\partial \zeta} \right) 
\mathbb{Y}_{cl-op}(u; z, \zeta).
\eea
\end{lemma}
\pf
The proof is the same as that of Lemma \ref{lemma-D-brack}. 
\epf


An analytic open-closed field algebra
$(V_{cl}, V_{op}, \mathbb{Y}_{cl-op})$ satisfies two  
associativities and two commutativities. These properties 
are very important for categorical formulation.

\begin{prop} [Associativity I]  \label{prop-asso-co-op}
For $u\in V_{cl}, v_1,v_2\in V_{op}, v'\in V_{op}'$ and
$z\in \HH,\zeta\in \overline{\HH}$,  we have
\beq \label{asso-co-op}
\langle v', \mathbb{Y}_{cl-op}(u; z, \zeta)
Y_{op}(v_1, r)v_2\rangle  
=\langle v', Y_{op}
( \mathbb{Y}_{cl-op}(u; z-r, \zeta-r)v_1, r)v_2\rangle 
\eeq
when $|z|,|\zeta|>r>0$ and $r>|r-z|,|r-\zeta|>0$. 
\end{prop}
\pf
We abbreviate the ``left (right) hand side'' as ``LHS (RHS)'' . 
For $z\in \HH$,
\bea
\mbox{ LHS of (\ref{asso-co-op}) } |_{\zeta=\bar{z}}
&=& m_{cl-op}^{(1;1)}(u; m_{op}^{(2)}(v_1,v_2; r,0); z,
\bar{z};0) \nn
&=& m_{cl-op}^{(1;2)}(u; v_1,v_2; z, \bar{z}; r,0) 
\eea
when $|z|>r>0$, and 
\bea
\mbox{ RHS of (\ref{asso-co-op}) } |_{\zeta=\bar{z}}
&=& m_{cl-op}^{(0;2)}(m_{cl-op}^{(1;1)}(u; v_1; z-r, \bar{z}-r;0),
v_2; r,0) \nn
&=& m_{cl-op}^{(1;2)}(u;v_1,v_2; z, \bar{z}; r,0) 
\eea
when $r>|r-z|>0$. Therefore (\ref{asso-co-op})
holds when $\zeta=\bar{z}$ and 
$|z|>r>|r-z|>0, z\in\HH$. Now replace $u$ in (\ref{asso-co-op})
by $(D_{cl}^{R})^k u$, $k\in \N$,  we obtain: 
\beq
\frac{\partial^k}{\partial \zeta^k}\lbar_{\zeta=\bar{z}}
\mbox{ LHS of (\ref{asso-co-op})}  
= \frac{\partial^k}{\partial \zeta^k}\lbar_{\zeta=\bar{z}}
\mbox{ RHS of (\ref{asso-co-op})} 
\eeq
when $z\in \HH$ and $|z|>r>|r-z|>0$. 
Then by the properties of analytic function, it is clear that 
(\ref{asso-co-op}) holds for all 
$z\in \HH, \zeta\in \overline{\HH}$ and 
$|z|,|\zeta|>r>0$ and $r>|r-z|,|r-\zeta|>0$. 
\epf

\begin{prop}[Associativity II] \label{prop-asso-co-co}
For $u_1, u_2 \in V_{cl}, v_1,v_2\in V_{op}, v'\in V_{op}'$ and
$z_1,z_2\in \HH,\zeta_1, \zeta_2\in \overline{\HH}$, we have 
\bea  \label{asso-co-co}
&& \langle w', \mathbb{Y}_{cl-op}(u_1; z_1, \zeta_1)
\mathbb{Y}_{cl-op}(u_2; z_2, \zeta_2)v_2\rangle  \nn
&& \hspace{1cm} =\langle v', \mathbb{Y}_{cl-op}
( \mathbb{Y}_{an}(u_1; z_1-z_2, \zeta_1-\zeta_2)
u_2; z_2, \zeta_2)v_2\rangle 
\eea
when $|z_1|, |\zeta_1|> |z_2|, |\zeta_2|$ and 
$|z_2| > |z_1-z_2|>0, |\zeta_2|>|\zeta_1-\zeta_2|>0$ and
$|z_2-\zeta_2|>|z_1-z_2|+|\zeta_1-\zeta_2|$.  
\end{prop}
\pf
The proof is similar to
that of (\ref{asso-co-op}). So we will only sketch the difference 
here. First, it is easy to show that 
(\ref{asso-co-co}) is true for $z_i\in \HH$, $\zeta_i=\bar{z}_i$, 
$i=1,2$ in a proper domain. Then replacing $u_1, u_2$ by 
$(D_{cl}^R)^p u_1$, $(D_{cl}^R)^qu_2$ 
for $p,q\in \N$ and using 
(\ref{D-r-Y-ffa}) and (\ref{D-op-brack}),  
we obtain that 
\beq
\frac{\partial^p}{\partial \zeta_1^p} 
\frac{\partial^q}{\partial \zeta_2^q}\lbar_{\zeta_i=\bar{z}_i}
\mbox{ LHS of (\ref{asso-co-co})}  
= \frac{\partial^p}{\partial \zeta_1^p} 
\frac{\partial^q}{\partial \zeta_2^q}\lbar_{\zeta_i=\bar{z}_i}
\mbox{ RHS of (\ref{asso-co-co})} 
\eeq
in a proper domain. 
By the property of analytic function again, 
the analytic extension of both sides of 
(\ref{asso-co-co}) to the following 
simply connected domain
$$
D_0 = \{ (z_1, \zeta_1, z_2, \zeta_2) |
z_i\in \HH, \zeta_i\in \overline{\HH},
|z_1|, |\zeta_1|> |z_2|, |\zeta_2| \}
$$
must be identical. Moreover, the additional 
restrictions of domain in the statement of
Associativity II guarantee the absolute 
convergence of both sides of (\ref{asso-co-co}). 
\epf

For an analytic open-closed field algebra, 
the map $\mathbb{Y}_{cl-op}$ can be uniquely extended to 
$V_{cl}\otimes V_{op}\times R$ where 
$$
R : = \{ (z, \zeta) \in \C^2 | z\in \HH \cup \R_+, \zeta\in \overline{\HH}\cup \R_+, z\neq \zeta \}. 
$$

\begin{prop}[Commutativity I] \label{comm-prop-1}
For $u\in V_{cl}, v_1,v_2\in V_{op}$ and $v'\in V_{op}'$,  
\beq  \label{comm-ocfa-coop-1}
\langle v', \mathbb{Y}_{cl-op}(u; z, \zeta)
Y_{op}(v_1, r)v_2\rangle,
\eeq
which is absolutely convergent when $z>\zeta > r>0$ (recall Definition 
\ref{an-ext-ocfa-def}), and 
\beq  \label{comm-ocfa-coop-2}
\langle v', Y_{op}(v_1, r)\mathbb{Y}_{cl-op}(u; z, \zeta)v_2\rangle,
\eeq
which is absolutely convergent when $r>z>\zeta>0$ (recall Definition
\ref{an-ext-ocfa-def}), are 
analytic continuation  of each other along the following path 
\footnote{The extended domain $R\backslash \{z=r, \text{or}, \zeta =r\}$ 
is simply connected for fixed $r>0$, 
all possible paths of analytic continuation  
are homotopically equivalent.}. 
\beq  \label{comm-bcft-1}
\epsfxsize 0.4\textwidth
\epsfysize 0.15\textwidth
\epsfbox{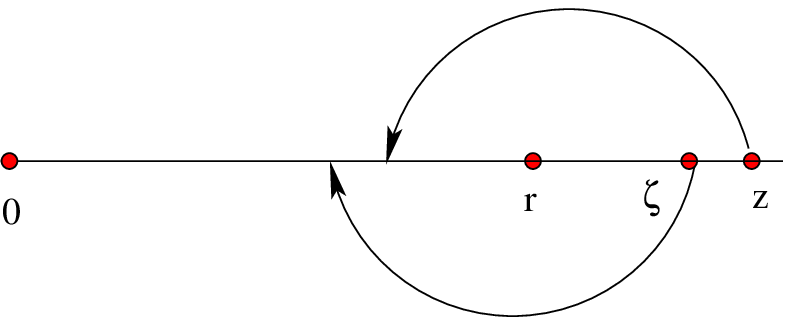}
\eeq
\end{prop}
\pf 
When $z\in \HH$, $\zeta=\bar{z}$ and $|z|>r>0$, 
(\ref{comm-ocfa-coop-1}) is absolutely convergent to
$m_{cl-op}^{(1;2)}(u;v_1,v_2; z,\bar{z}; r,0)$. 
Hence, if we analytic extend the analytic function 
(\ref{comm-ocfa-coop-1}) to 
$$
D_1:=\{ (z, \zeta) | z\in \HH, \zeta=\bar{z}, |z|=r \},
$$
then its value on $D_1$ must equal to 
$m_{cl-op}^{(1;2)}(u;v_1,v_2; z,\bar{z}; r,0)$ when $z\in \HH,|z|=r$ 
by the continuity of $m_{cl-op}^{(1;2)}$. 

Similarly, the unique extension of (\ref{comm-ocfa-coop-2})
from $r>|z|, |\zeta|>0$ to $D_1$ also equal to 
$m_{cl-op}^{(1;2)}(u;v_1,v_2; z, \bar{z}; r,0)$ when $z\in \HH, |z|=r$. 

Both (\ref{comm-ocfa-coop-1}) and (\ref{comm-ocfa-coop-2})
can be uniquely extended to two analytic functions on 
$\HH \times \overline{\HH}$. By the discussion above,  
these two extended analytic functions of $(z,\zeta)$ take
the same value on $D_1$ which is of lower dimension.

Now we replace $u$ by $(D_{cl}^L)^m(D_{cl}^R)^nu, m,n\in \N$ and
repeat the arguments in the proof of 
Proposition \ref{prop-asso-co-op}.
We obtain that all the derivatives of 
the two extended analytic functions also match on $D_1$. By the 
properties of analytic functions, these two extended functions
are identical on $\HH\times \overline{\HH}$. 
Then their extensions to $(z,\zeta)\in 
R\backslash \{z=r, \text{or}, \zeta =r\}$ are
also identical. Thus we have proved the first commutativity.  
\epf

\begin{prop}[Commutativity II] \label{comm-prop-2}
For $u_1,u_2\in V_{cl}, v\in V_{op}$ and $v'\in V_{op}'$,
\beq   \label{comm-ocfa-coco-1}
\langle v', \mathbb{Y}_{cl-op}(u_1; z_1, \zeta_1)
\mathbb{Y}_{cl-op}(u_2; z_2, \zeta_2)v\rangle,
\eeq
which is absolutely convergent when $z_1>\zeta_1>z_2>\zeta_2>0$, and
\beq  \label{comm-ocfa-coco-2}
\langle v', \mathbb{Y}_{cl-op}(u_2; z_2, \zeta_2)
\mathbb{Y}_{cl-op}(u_1; z_1, \zeta_1)v\rangle, 
\eeq
which is absolutely convergent when $z_2>\zeta_2>z_1>\zeta_1>0$, 
are analytic continuation of each other along the following paths.
\beq \label{comm-path-2}
\epsfxsize 0.8\textwidth
\epsfysize 0.15\textwidth
\epsfbox{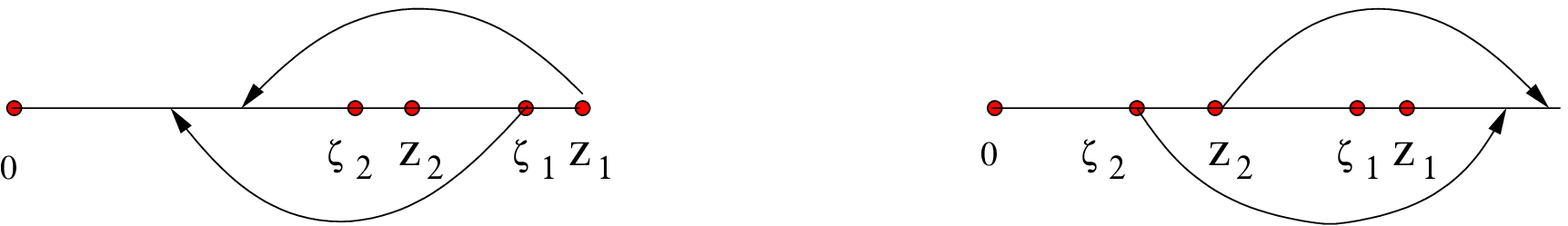}
\eeq
\end{prop}
\pf
Commutativity II follows directly from 
the Associativity II (\ref{asso-co-co})
and the skew-symmetry of the full field algebra $V_{cl}$ \cite{HKo2}. 
Here, we give a more direct proof which is
similar to that of Proposition \ref{comm-prop-1}.  
The unique extension of (\ref{comm-ocfa-coco-1}) from 
$|z_1|,|\zeta_1|>|z_2|,|\zeta_2|>0$ to 
$$
D_2: = \{ (z_1,\zeta_1, z_2, \zeta_2)| z_1,z_2\in M_{\HH}^2, |z_1|=|z_2|,
\zeta_i=\bar{z}_i, i=1,2 \}
$$
equal to $m_{cl-op}^{(2;1)}(u_1,u_2;v;z_1,\bar{z}_1,z_2,
\bar{z}_2;0)$, 
the unique extension of (\ref{comm-ocfa-coco-2})
from $|z_2|,|\zeta_2|>|z_1|,|\zeta_1|>0$ to $D_2$ also 
match with $m_{cl-op}^{(2;1)}(u_1,u_2;v;z_1,\bar{z}_1,
z_2,\bar{z}_2;0)$. 
By the similar argument as the proof of the first commutativity, 
we see that (\ref{comm-ocfa-coco-1}) in 
$|z_1|,|\zeta_1|>|z_2|,|\zeta_2|>0$ and
(\ref{comm-ocfa-coco-1}) in $|z_2|,|\zeta_2|>|z_1|,|\zeta_1|>0$
are analytic continuation of each other along the following paths.
\beq
\begin{picture}(14,2)
\put(1, 0){\resizebox{12cm}{2cm}
{\includegraphics{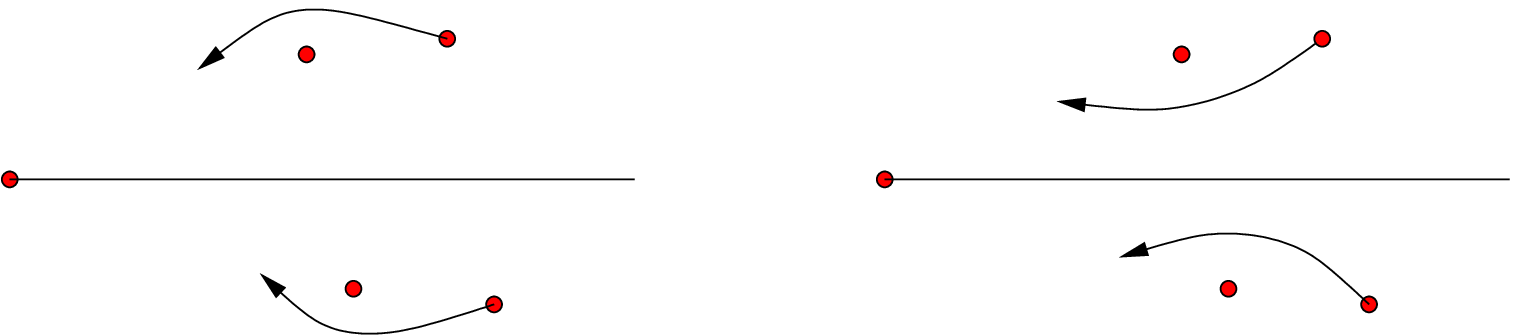}}}
\put(0.9, 0.5){$0$}
\put(3.6, 1.5){$z_2$}\put(4.6, 1.5){$z_1$}
\put(4, 0.4){$\zeta_2$}\put(5.1, 0.4){$\zeta_1$}

\put(7.8,0.5){$0$}
\put(10.5, 1.8){$z_2$}\put(11.6,1.8){$z_1$}
\put(10.8, 0){$\zeta_2$}\put(12, 0){$\zeta_1$}
\end{picture}
\eeq
Then it is obvious to see that the unique extension of 
both (\ref{comm-ocfa-coco-1}) and (\ref{comm-ocfa-coco-2}) to 
the subdomain of $\{(z_1,\zeta_1, z_2, \zeta_2)\in \R_+^4\}$, 
where $z_1, \zeta_1, z_2, \zeta_2$ have distinct values,
are analytic extension of each other 
along the paths (\ref{comm-path-2}). 
\epf

\subsection{Open-closed field algebras over $U$}
We gradually add more structures on an analytic
open-closed field algebra. At the end of this process, we will arrive 
at a notion called open-closed field algebra over 
a vertex operator algebra $U$. 
Let $(V^L, Y_{V^L}, \one^L)$ and $(V^R, Y_{V^R}, \one^R)$ 
be two vertex algebras.
We now consider an analytic open-closed field algebra  
$(V_{cl}, V_{op}, \mathbb{Y}_{cl-op})$
such that $V_{cl}$ is a full field algebra over $V^L\otimes V^R$. 
For such full field algebra $V_{cl}$,
we will not distinguish $V^L$ with $V^L\otimes \one^R \subset V_{cl}$ 
and $V^R$ with $\one^L \otimes V^R \subset V_{cl}$.

\begin{lemma} 
For $u^L\in V^L$ and $u^R\in V^R$,
$\mathbb{Y}_{cl-op}(u^L; z, \zeta)$ is independent of $\zeta$
and $\mathbb{Y}_{cl-op}(u^R; z, \zeta)$ is independent of $z$. 
\end{lemma}
\pf
For $u^L\in V^L, w\in V_{op}, w'\in (V_{op})'$,
using the associativity (\ref{asso-co-co}), we have
\bea  \label{u-l}
\langle w', \mathbb{Y}_{cl-op}(u^L; z, \zeta)w\rangle
&=& 
\langle w', \mathbb{Y}_{cl-op}(u^L; z, \zeta)
\mathbb{Y}_{cl-op}(\one_{cl}; z_1,\zeta_1)w\rangle   \nn
&=&   \langle w', \mathbb{Y}_{cl-op}(
\mathbb{Y}_{an}(u^L; z-z_1, \zeta-\zeta_1)\one_{cl}; 
z_1,\zeta_1)w\rangle
\eea
when $|z_1|>|z-z_1|>0$, $|\zeta_1|>|\zeta-\zeta_1|>0$
and $|z-z_1|+|\zeta-\zeta_1|<|z_1-\zeta_1|$.
The right hand side of (\ref{u-l}) is independent of 
$\zeta$ and the left hand side of (\ref{u-l}) is analytic in $\zeta$. 
Hence $\mathbb{Y}_{cl-op}(u^L; z, \zeta)$ 
is independent of $\zeta$ for all $z\in \HH$. 
Similarly, $\mathbb{Y}_{cl-op}(u^R; z, \zeta)$
is independent of $z$ for all $\zeta\in \overline{\HH}$
and $u^R\in V^R$. 
\epf  

In order to emphasis these $\zeta$- 
or $z$-independence properties, 
we denote them simply as $\mathbb{Y}_{cl-op}(u^L, z)$ 
and $\mathbb{Y}_{cl-op}(u^R, \zeta)$ 
for $u^L\in V^L$ and $u^R \in V^R$ respectively.

Replace $u$ in (\ref{d-op-brack}) by $u^L\in V^L$ and $u^R\in V^R$ 
respectively, we obtain 
\bea   \label{d-op-brack-u}
[\mathbf{d}_{op}, \mathbb{Y}_{cl-op}(u^L, z)] &=& 
\mathbb{Y}_{cl-op}(\mathbf{d}_{cl}^L u^L, z) 
+ z\frac{\partial}{\partial z} 
\mathbb{Y}_{cl-op}(u^L, z), \nn
\left[ \mathbf{d}_{op}, \mathbb{Y}_{cl-op}(u^R, \zeta) \right] &=& 
\mathbb{Y}_{cl-op}(\mathbf{d}_{cl}^R u^R, \zeta) 
+ \zeta \frac{\partial}{\partial \zeta} 
\mathbb{Y}_{cl-op}(u^R, \zeta).
\eea
As a consequence, we have 
\bea \label{u-l-expansion}
\mathbb{Y}_{cl-op}(u^L, z) &=& \sum_{n\in \R} u_n^L z^{-n-1},   \nn
\mathbb{Y}_{cl-op}(u^R, \zeta) &=& \sum_{n\in \R} u_n^R \zeta^{-n-1},
\eea
where $u_n^L, u_n^R \in \edo V_{op}$ 
and $\wt u_n^L = \lwt u^L -n-1$ and $\wt u_n^R = \rwt u^R -n-1$,
and $z^n=e^{n\log z}$ and $\zeta^n=e^{n\overline{\log \zeta}}$.

Moreover, we have 
\bea
[D_{op}, \mathbb{Y}_{cl-op}(u^L, z)] &=& \mathbb{Y}_{cl-op}(D_{cl}^L u^L, z)
= \frac{\partial}{\partial z}\mathbb{Y}_{cl-op}(u^L, z),  
\label{D-Y-clop-L}  \\
\left[ D_{op}, \mathbb{Y}_{cl-op}(u^R, \zeta)\right] &=& 
\mathbb{Y}_{cl-op}(D_{cl}^R u^R, \zeta)
= \frac{\partial}{\partial \zeta}\mathbb{Y}_{cl-op}(u^R, \zeta), 
\nonumber
\eea
which further implies that 
\bea  \label{conj-D-op}
e^{aD_{op}}\mathbb{Y}_{cl-op}(u^L, z)e^{-aD_{op}} &=& 
\mathbb{Y}_{cl-op}(u^L, z+a) \nn
e^{aD_{op}}\mathbb{Y}_{cl-op}(u^L, \zeta)e^{-aD_{op}} &=& 
\mathbb{Y}_{cl-op}(u^R, \zeta+a)
\eea
for $|z|>|a|, z+a\in \HH$ 
and $|\zeta|>|a|, \zeta+a\in \overline{\HH}$ respectively. 

\begin{lemma}
For any $u^L\in V^L\otimes \one^R$ and $u^R\in \one^L \otimes V^R$, 
the following two limits:
$$\lim_{z\rightarrow 0} \mathbb{Y}_{cl-op}(u^L, z) \one_{op}, 
\quad \quad \lim_{\zeta \rightarrow 0} \mathbb{Y}_{cl-op}(u^R, \zeta)\one_{op}$$
exist in $V_{op}$. 
\end{lemma}
\pf
By the associativity and the creation
property of open-string vertex algebra, we have
\bea  \label{u-l-lim-1}
\mathbb{Y}_{cl-op}(u^L, z+r)\one_{op} 
&=&\mathbb{Y}_{cl-op}(u^L, z+r)Y_{op}(\one_{op}, r)\one_{op} \nn
&=& Y_{op}(\mathbb{Y}_{cl-op}(u^L, z)\one_{op}, r)\one_{op}  \nn
&=& e^{rD_{op}} \mathbb{Y}_{cl-op}(u, z) \one_{op}
\eea
when $|z+r|> r >|z|>0$. For fixed $z\in \HH$, 
the left hand side of (\ref{u-l-lim-1})
is an analytic function valued in $\overline{V}_{op}$ 
on the domain $\{ r\in \C | z+r\neq 0 \}$. 
The right hand side of (\ref{u-l-lim-1}), as a power series
of $r$, is absolutely convergent when
$|r|>|z|>0$. By the general property of power series, 
the right hand side of (\ref{u-l-lim-1}) 
is absolutely convergent for all $r\in \C$ to a singlevalued
analytic function. Because both sides of (\ref{u-l-lim-1})
are analytic functions, the equality (\ref{u-l-lim-1}) must
hold for all $r\in \C$. In particular, 
$\lim_{r\rightarrow -z} \mathbb{Y}_{cl-op}(u^L, z+r)\one_{op}$  
exists. Equivalently, 
$\lim_{z\rightarrow 0} \mathbb{Y}_{cl-op}(u^L, z)\one_{op}$ exists. 
By the expansion (\ref{u-l-expansion}),
we must have $u_n^L\one_{op}=0$ for all $n> -1$. 
Moreover, it is also easy to see that $u_n^L\one_{op}=0$ for
any $n\notin -\Z_+$ by (\ref{D-Y-clop-L})
(see the proof of Proposition 1.8 in \cite{HKo1}). Therefore, we have
$$
\lim_{z\rightarrow 0} \mathbb{Y}_{cl-op}(u^L, z)\one_{op} = u_{-1}^L \one_{op}
\in V_{op}. 
$$
The proof of the existence of the second limit is entirely the same.
\epf

By the Lemma above, we can define two maps
$h^L : V^L \rightarrow V_{op}$ and
$h^R:  V^R \rightarrow V_{op}$ as follows: for all $u^L\in V^L$ and 
$u^R\in V^R$, 
\bea  \label{h-L-R}
&&h^L: u^L\mapsto \lim_{z\rightarrow 0}
\mathbb{Y}_{cl-op}(u^L, z)\one_{op}, \nn
&&h^R: u^R\mapsto \lim_{\zeta\rightarrow 0}
\mathbb{Y}_{cl-op}(u^R, \zeta)\one_{op}.
\eea
Notice also that $h^L, h^R$ preserve the weights. Namely
\bea
\wt h^L(u^L) &=& \wt (u_{-1}^L\one_{op}) = \lwt u^L,    \nn
\wt h^R(u^R) &=& \wt (u_{-1}^R\one_{op}) = \rwt u^R.   \nonumber
\eea
Therefore both $h^L$ and $h^R$ can be naturally extended
to maps $\overline{V^L} 
\rightarrow \overline{V_{op}}$ and
$\overline{V^R} \rightarrow \overline{V_{op}}$.
We still denote the extended maps as $h^L$ and $h^R$ respectively.

\begin{lemma}
For $u^L\in V^L, u^R\in V^R$, we have
\bea  \label{Y-co-creat}
\mathbb{Y}_{cl-op}(u^L, z) \one_{op} &=& e^{zD_{op}} h^L(u^L),  \nn
\mathbb{Y}_{cl-op}(u^R, \zeta) \one_{op} &=& e^{\zeta D_{op}} h^R(u^R). 
\eea
\end{lemma}
\pf
Since we have shown that (\ref{u-l-lim-1}) holds
for all $r\in \C, z\in \HH$ and both sides of 
(\ref{u-l-lim-1}) are analytic for all $r,z\in \C$,
if we take the limit 
$\lim_{z\rightarrow 0}$ on both sides of (\ref{u-l-lim-1}), the 
equality should still hold. Thus we obtain 
the first identity in (\ref{Y-co-creat}). 

The proof of the second identity is entirely the same. 
\epf

\begin{prop}
$h^L$ and $h^R$ are homomorphisms from $V^L$ and $V^R$ respectively
to their images, viewed as graded vertex algebras.  
\end{prop}
\pf
We have shown that $h^L, h^R$ preserve gradings. 
By the identity property of open-closed field algebra, 
we have
$$
h^L(\one_{cl})=\one_{op}.
$$
Next, for $u^L \in V^L$ and $r>0$,  we define
\bea \label{prop-h-l-r-equ-1}
\mathbb{Y}_{cl-op}(u^L, r)w &:=& \lim_{z\rightarrow r} 
\mathbb{Y}_{cl-op}(u^L, z)w \nn
&=& \lim_{z\rightarrow r} \mathbb{Y}_{cl-op}(u^L, z)Y_{op}(\one_{op}, r)w
\eea
where the limit is taken along a path from a fixed initial 
point in $\HH$ to $r>0$. 
Since $\mathbb{Y}_{cl-op}(u^L, z)w$ is analytic in $\C^{\times}$, 
the limit is independent of the path we choose. 
So we choose a path in the domain $\{ z\in \HH| |z|>r>|z-r|>0\}$.
In this domain, we can apply the associativity (\ref{asso-co-op}) to 
the right hand side of (\ref{prop-h-l-r-equ-1}). We obtain
\bea
\mathbb{Y}_{cl-op}(u^L, r)w
&=& \lim_{z\rightarrow r} Y_{op}(\mathbb{Y}_{cl-op}(u^L, z-r)\one_{op}, r)w \nn
&=& Y_{op}(h^L(u^L), r)w.   \nonumber
\eea

For $|z|>r>0$, by (\ref{Y-co-creat}), we have, 
\beq
\mathbb{Y}_{cl-op}(\mathbb{Y}(u^L; r, r)v^L, z)\one_{op} =
e^{zD_{op}} h^L(\mathbb{Y}(u^L; r, r)v^L), \label{Y-clop-Y-one}
\eeq
the right hand side of which is absolutely convergent for
all $z\in \C$. And both sides are analytic in $z$. Therefore 
$\mathbb{Y}_{cl-op}(\mathbb{Y}(u^L; r, r)v^L, z)\one_{op}$
is absolutely convergent for all $z\in \C$ and 
(\ref{Y-clop-Y-one}) holds for all $z\in \C$. 
By the associativity, we have
\beq  \label{prop-h-l-r-equ-2}
\mathbb{Y}_{cl-op}(\mathbb{Y}(u^L; r, r)v^L, z)\one_{op} = 
\mathbb{Y}_{cl-op}(u^L, r+z)\mathbb{Y}_{cl-op}(v^L, z)\one_{op}
\eeq
for $|r+z|>|z|>r>0$. Again both sides of (\ref{prop-h-l-r-equ-2})
are analytic in $z$. Hence the left hand side of 
(\ref{prop-h-l-r-equ-2}) defined for all $z\in \C$ is the 
analytic extension of the right hand side of 
(\ref{prop-h-l-r-equ-2}), which is defined on $\{ |r+z|>|z| \}$.
Since the extension is free of singularity on entire $\C$, 
the right hand side of (\ref{prop-h-l-r-equ-2}) must be well-defined
on entire $\C$. Therefore, we must have 
\beq \label{prop-h-l-r-equ-3}
\lim_{z\rightarrow 0} 
\mathbb{Y}_{cl-op}(\mathbb{Y}(u^L,r)v^L, z)\one_{op}  
= \lim_{z\rightarrow 0} 
\mathbb{Y}_{cl-op}(u^L, r+z)\mathbb{Y}_{cl-op}(v^L, z)\one_{op}. 
\eeq

Combining above results, we have 
\bea
h^L(\mathbb{Y}(u^L, r)v^L) &=& 
\lim_{z\rightarrow 0} e^{zD_{op}} h^L(\mathbb{Y}(u^L, r)v^L)  \nn
&=& \lim_{z\rightarrow 0} 
\mathbb{Y}_{cl-op}(\mathbb{Y}(u^L,r)v^L, z)\one_{op}  \nn
&=& \lim_{z\rightarrow 0} 
\mathbb{Y}_{cl-op}(u^L, r+z)\mathbb{Y}_{cl-op}(v^L, z)\one_{op} \nn
&=& \mathbb{Y}_{cl-op}(u^L, r)h^L(v^L) \nn
&=& Y_{op}(h^L(u^L), r) h^L(v^L).   \nonumber
\eea
Thus $h^L$ is a vertex algebra homomorphism. 
The proof for $h^R$ is entirely the same. 
\epf


Let $(U, Y, \one, \omega)$ be a vertex operator algebra with 
central charge $c$. 
$U$ and $U\otimes U$ naturally give an analytic 
open-closed field algebra, in which 
$\mathbb{Y}_{cl-op}(\cdot; z,\zeta)\cdot$ is
given by  
\bea
\mathbb{Y}_{cl-op}(u\otimes v; z, \zeta)w &=& Y(u,z)Y(v,\zeta)w,
\quad\quad |z|>|\zeta|>0,  \nn
&=& Y(v,\zeta)Y(u,z)w, \quad\quad |\zeta|>|z|>0 \label{U-UU-act}
\eea
for $u,v,w\in U$.
In this case, $h^L: u\otimes \one \mapsto u$ and 
$h^R: \one\otimes u \mapsto u$. We denote this open-closed
field algebra as $(U\otimes U, U)$. 

In general, let $\rho^L, \rho^R \in \text{Aut}(U)$ where
$\text{Aut}(U)$ is the set of automorphisms
of $U$ as vertex operator algebra. We can obtain 
a new action of $U\otimes U$ on $U$ by 
composing (\ref{U-UU-act}) with the automorphism 
$\rho^L\otimes \rho^R: U\otimes U\rightarrow U\otimes U$. 
Namely, there exists another
open-closed field algebra structure on $U$ and $U\otimes U$, 
in which 
$\mathbb{Y}_{cl-op}(u\otimes v; z, \zeta)w$, for $u,v,w\in U$,
is given by 
\bea
&&Y(\rho^L(u),z)Y(\rho^R(v),\zeta)w,
\quad\quad \mbox{ for $|z|>|\zeta|>0$,}  \nn
&& Y(\rho^R(v),\zeta)Y(\rho^L(u),z)w, \quad\quad 
\mbox{ for $|\zeta|>|z|>0$. }
\eea 
In this case, $h^L: u\otimes \one \mapsto \rho^L(u)$ and 
$h^R: \one\otimes u \mapsto \rho^R(u)$. We denote such 
open-closed field algebra as $(U\otimes U, U, \rho^L, \rho^R)$. 
In particular, $(U\otimes U, U, \id_U, \id_U)$
is just $(U\otimes U, U)$.

\begin{rema}  {\rm
$(U\otimes U, U, \rho^L, \rho^R)$ for general automorphisms
$\rho^L$ and $\rho^R$ is very interesting in physics. 
But it adds some technical subtleties in later formulations. 
So we postpone its study to future publications. In this
work, we focus on $(U\otimes U, U)$.
}
\end{rema}

\begin{defn} \label{def-ocfa-U}
{\rm 
Let $(U, Y, \one, \omega)$ be a vertex operator algebra.
An {\it open-closed field algebra over $U$} is an 
analytic open-closed field algebra 
$$
( (V_{cl}, m_{cl}, \iota_{cl}), (V_{op}, Y_{op}, \iota_{op}),  
\mathbb{Y}_{cl-op}),
$$ 
where $(V_{cl}, m_{cl}, \iota_{cl})$ 
is a conformal full field algebra over 
$U\otimes U$ and $(V_{op}, Y_{op}, \iota_{op})$ 
is an open-string vertex operator algebra
over $U$, satisfying the following conditions:
\bnu
\item  {\em $U$-invariant boundary condition}: 
$h^L=h^R=\iota_{op}$.

\item {\em Chirality splitting property}: 
$\forall u\in V_{cl}$,  
$u=u^L\otimes u^R\in W^L\otimes W^R \subset V_{cl}$ 
for some $U$-modules $W^L, W^R$. 
There exist $U$-modules $W_1, W_2$ and intertwining operators 
$\Y^{(1)}, \Y^{(2)}, \Y^{(3)}, \Y^{(4)} $ of type 
$\binom{V_{op}}{W^LW_1}$, 
$\binom{W_1}{W^RV_{op}}$, $\binom{V_{op}}{W^RW_2}$,
$\binom{W_2}{W^LV_{op}}$ respectively\footnote{It was proved in \cite{HKo1} that 
$V_{op}$ is a $U$-modules, and in \cite{HKo2}
that $V_{cl}$ is a $U\otimes U$-modules.}, such that
\beq \label{ch-ach-Y-co-1}
\langle w', \mathbb{Y}_{cl-op}(u; z, \zeta)w\rangle
=  \langle w', \Y^{(1)}(u^L, z)\Y^{(2)}(u^R, \zeta) w\rangle
\eeq
when $|z|>|\zeta|>0$ (recall the convention (\ref{branch-cut}) 
and (\ref{power-fun})), and 
\beq \label{ch-ach-Y-co-2}
\langle w', \mathbb{Y}_{cl-op}(u; z, \zeta)w\rangle
= \langle w', \Y^{(3)}(u^R, \zeta)\Y^{(4)}(u^L, z) w\rangle
\eeq
when $|\zeta|>|z|>0$
for all $u\in V_{cl}, w\in V_{op}, w'\in V_{op}$.  
\enu    

In the case that $U$ is generated by $\omega$, i.e.
$U=\langle \omega \rangle$, the $\langle \omega\rangle$-invariant
boundary condition is simply called {\it conformal invariant boundary 
condition}. We also call open-closed
field algebra over $\langle \omega\rangle$ 
{\it open-closed conformal field algebra}. 
}
\end{defn}

\begin{rema}  {\rm
The $U$-invariant boundary condition actually says that 
the open-closed field algebra over $U$ contains $(U\otimes U, U)$ as
a subalgebra. 
If we only want to construct 
algebras over Swiss-cheese partial operad,
the $U$-invariant boundary condition in Definition \ref{def-ocfa-U}
can be weaken to the {\it conformal invariant boundary condition}: 
\beq  \label{conf-inv-bc}
h^L|_{\langle \omega \rangle} = h^R|_{\langle \omega \rangle} 
= \iota_{op}|_{\langle \omega \rangle}.
\eeq
These situations appear in physics in 
the study of the so-called symmetry breaking boundary conditions
(see for example \cite{FS1}\cite{FS2} and references therein). 
All examples studied in this work and \cite{Ko2} 
satisfy the $U$-invariant boundary condition. 
We leave the study of general symmetry-broken 
situations to the future.
}
\end{rema}

\begin{rema}
{\rm 
The chirality splitting property is a 
very natural condition because the interior sewing operation of 
$\tilde{\mathfrak{S}}^c$ is defined by a double sewing 
operation as given in (\ref{I-sew-ext}). 
Unfortunately, we do not know whether this chirality
splitting property is necessary for general constructions of 
algebras over $\tilde{\mathfrak{S}}^c$.
}
\end{rema}


For an open-closed field algebra over $V$, 
there are three Virasoro elements, 
$\omega_{op}:=\iota_{op}(\omega)$ and 
$\omega^L:=\iota_{cl}(\omega\otimes \one)$ and 
$\omega^R:=\iota_{cl}(\one \otimes \omega)$, and we have
\bea
Y_{op}(\omega_{op}, r) &=& \sum_{n\in \Z} L(n) r^{-n-2},  \nn
\mathbb{Y}(\omega^L, z) &=& \sum_{n\in \Z} L^L(n) z^{-n-2},  
\nn
\mathbb{Y}(\omega^R, \zeta) &=& \sum_{n\in \Z} L^R(n) 
\zeta^{-n-2}. \nonumber
\eea
where $L^L(n)=L(n)\otimes 1$ and $L^R(n)=1\otimes L(n)$ for
$n\in \Z$.

When $U=V$ a vertex operator algebra 
satisfying the conditions in 
Theorem \ref{ioa}, we have a very simple description of 
open-closed field algebra over $V$ given in the following Theorem. 
\begin{thm} \label{prop-op-cl-V}
An open-closed field algebra over $V$ is equivalent to 
the following structure: $(V_{op}, Y_{op}, \iota_{op})$
an open-string vertex operator algebra over
$V$ and $(V_{cl}, m_{cl}, \iota_{cl})$ 
a conformal full field algebra over $V\otimes V$, 
together with a vertex map $\mathbb{Y}_{cl-op}(\cdot; z, \zeta)\cdot$ 
given by intertwining operators 
$\Y^{(i)}, i=1,2,3,4$ as in (\ref{ch-ach-Y-co-1}) and 
(\ref{ch-ach-Y-co-2}),  satisfying the unit property:
\beq  \label{Y-co-unit}
\mathbb{Y}_{cl-op}(\one_{cl}; z, \bar{z}) v := 
\Y^{(1)}(\one, z)\Y^{(2)}(\one, \bar{z})v = v, 
\quad\quad \forall v\in V_{op}, 
\eeq
Associativity I (\ref{asso-co-op}), Associativity II (\ref{asso-co-co}) 
and Commutativity I given in Proposition \ref{comm-prop-1}. 
\end{thm}
\pf
It is clear that an open-closed field algebra over $V$
gives the data and properties included in the statement of the 
Theorem. We only need to show that
such data is sufficient to reconstruct an open-closed
field algebra over $V$. Moreover, such open-closed
field algebra over $V$ with the given data is unique.

Since $V$ satisfies the condition in Theorem \ref{ioa}, 
the conditions listed in 
the Definition \ref{an-ext-ocfa-def} are all 
automatically satisfied. 
In particular, for $v_i\in V_{op}, u_i\in V_{cl}, i=1,\ldots, n$
and $v'\in V_{op}'$, the following series: 
\bea \label{gen-corr-bcft}
\langle v', \mathbb{Y}_{cl-op}(u_1; z_1,\zeta_1)Y_{op}(v_1, r_1) \cdots 
\mathbb{Y}_{cl-op}(u_n; z_n,\zeta_n)Y_{op}(v_n, r_n) \one_{op} \rangle.
\eea
is absolutely convergent when 
$z_1 >\zeta_1>r_1>\cdots >z_n>\zeta_n>r_n > 0$, and 
can be extended by analytic continuation
to a multi-valued analytic function for variables 
$z_i, \zeta_i \in \C^{\times}, i=1, \ldots, n$ 
with possible singularities only when 
two of $z_i, \zeta_j, r_k$ are equal. 
Using this property of (\ref{gen-corr-bcft}),
we can define a non-analytic but single-valued smooth 
function in $\Lambda^n \times M_{\HH}^n$ as follows.

Let $(r_1, \dots, r_n) \in \Lambda^n$ 
and $(\xi_1, \dots, \xi_n)\in M_{\HH}^n$, $n\in \Z_+$.
Let $\gamma_1$ be a smooth path from $(3n, 3(n-1),\dots, 1)$ to  
$(\xi_1,\ldots, \xi_n)$ such that $\gamma_1((0,1)) \subset M_{\HH}^n$. 
Let $\gamma_0: [0,1]\rightarrow \R_+^n$ be so that
$$
\gamma_0(t) = ((1-t) (3n-1) + t 3n, \dots, (1-t)2 + t 3).
$$
Clearly, $\gamma_0$ is the straight line from 
$(3n-1, 3(n-1)-1,\dots, 2)$ 
to $(3n, 3(n-1), \dots, 3)$ in $\R_+^n$.  Then we define 
a path $\gamma_2: [0,1]\in (\overline{\HH} \cup \R_+)^n$ to
be the composition $\gamma_2= \bar{\gamma}_1 \circ \gamma_0$ where
$\bar{\gamma}_1$ is the complex conjugation of $\gamma_1$. 
It is clear that $\gamma_2$ is a path from 
$(3n-1, 3(n-1)-1,\dots, 2)$ to
$(\bar{\xi}_1,\ldots, \bar{\xi}_n)$. 
Let $\gamma_3: [0,1]\rightarrow \R_+^n$ be so that 
$$
\gamma_3(t) = ((1-t) (3n-2) + t r_1, \dots, (1-t)1 + t r_n).
$$
So $\gamma_3$ is the straight line from 
$(3n-2, 3(n-1)-2, \dots, 1)$ to $(r_1, \dots, r_n)\in \Lambda^n$. 

Combining $\gamma_1, \gamma_2, \gamma_3$, we obtain a path 
$\gamma$ in $\C^{3n}$ from the initial point 
$(3n, 3n-1, 3n-2, \dots, 3, 2, 1)\in \R_+^{2n}$
to the final point 
$(\xi_1, \bar{\xi}_1, r_1, \dots, \xi_n, \bar{\xi}_n, r_n)$ 
in the obvious way.

Then we define 
$$
\langle v', m_{cl-op}^{(n;n)}(u_1, \dots, u_n, v_1, \dots ,v_n;
\xi_1,\bar{\xi}_1, \ldots, \xi_n, \bar{\xi}_n, r_1, \dots, r_n)\rangle
$$ 
to be the value obtained from the value of 
(\ref{gen-corr-bcft}) at the initial point 
$(3n, 3n-1, 3n-2, \dots, 3, 2, 1)$ 
by analytic continuation along the path $\gamma$. 
Following a similar argument as in the proof of 
Theorem 2.11 in \cite{HKo2}, it is easy to show, 
that such defined 
\beq  \label{m-def-0}
m_{cl-op}^{(n;n)}(u_1,\dots, u_n; v_1, \dots, v_{n};
\xi_1,\bar{\xi}_1, \ldots, \xi_n, \bar{\xi}_n, r_1, \dots, r_n), 
\eeq
is independent of the choice of $\gamma_1$ and its initial points. 
Moreover, such defined (\ref{m-def-0}) is single-valued and smooth
in $\HH \times \overline{\HH}$.

For $n>l\geq 0$, we define
\bea
&&\hspace{-0.2cm}m_{cl-op}^{(l;n)}(u_1,\dots, u_l; v_1, \dots, v_n;
\xi_1,\bar{\xi}_1, \ldots, \xi_l, \bar{\xi}_l; r_1, \dots, r_n), \nn
&&\hspace{0cm} 
:= m_{cl-op}^{(n;n)}(u_1,\dots,u_l, \one_{cl}, \dots, \one_{cl}; 
v_1, \dots, v_n;
\xi_1,\bar{\xi}_1, \ldots, \xi_n, \bar{\xi}_n; r_1, \dots, r_n),
\nonumber
\eea
and for $l>n\geq 0$, we define
\bea
&&\hspace{-0.2cm}m_{cl-op}^{(l;n)}(v_1,\dots, v_l; w_1, \dots, w_n; 
\xi_1,\bar{\xi}_1, \ldots, \xi_l, \bar{\xi}_l; r_1, \dots, r_n),  \nn
&&\hspace{0cm} := 
m_{cl-op}^{(l;l)}(u_1,\dots,u_l; v_1, \dots, v_{n}, 
\one_{op},\dots, \one_{op};  
\xi_1,\bar{\xi}_1, \ldots, \xi_l, \bar{\xi}_l; r_1, \dots, r_l),
\eea
and for $n=l=0$, we define $m_{cl-op}^{(0;0)}(1)=\one_{op} \in V_{op}$.

Immediately following from the construction of $m_{cl-op}^{(l;n)}$,
we have, for all $v\in V_{op}$, 
$$
m_{cl-op}^{(0;1)}(v;0) = m_{cl-op}^{(1;1)}(\one_{cl};v;z;0) 
= \mathbb{Y}_{cl-op}(\one_{cl}; z,\bar{z})v 
= v.
$$

Now we show the permutation axiom for $m_{cl-op}^{(l;n)}$.  
This is enough to just consider adjacent permutations 
$(ii+1), i=1,\dots, l-1$ 
because they generate the whole permutation group. 
We can just consider $(12)$ because all the other cases
are exactly the same. This is amount to show that 
\bea  \label{perm-axiom-equ}
&&\hspace{-0.5cm} m_{cl-op}^{(l;n)}(v_1, v_2, \dots, v_l; 
w_1, \dots, w_n; 
\xi_1,\bar{\xi}_1, \xi_2, \bar{\xi}_2, \dots, 
\xi_l, \bar{\xi}_l; r_1, \dots, r_n),  \nn
&&\hspace{-0.2cm}=m_{cl-op}^{(l;n)}(v_2, v_1,\dots, v_l; 
w_1, \dots, w_n; \xi_2, \bar{\xi}_2, 
\xi_1,\bar{\xi}_1, \ldots, \xi_l, \bar{\xi}_l; r_1, \dots, r_n). 
\eea
By our construction, the only difference of the two sides of 
(\ref{perm-axiom-equ}) is that they are obtained by analytic
continuation along paths with different initial points. 
The initial points of the path for the left hand side of 
(\ref{perm-axiom-equ}) is $z_1>\zeta_1>r_1>z_2>\zeta_2>r_2>\dots$, 
that for the right hand side of  
(\ref{perm-axiom-equ}) is $z_2>\zeta_2>r_1>z_1>\zeta_1>r_2>\dots$. 
But by the commutativity, 
the value of (\ref{gen-corr-bcft}) at these two initial points
are analytic continuation of each other along the 
paths given in the commutativity axiom of open-closed field algebra. 
Hence the equality (\ref{perm-axiom-equ}) follows.

((\ref{d-m-co}) and (\ref{D-m-co})) can be proved by 
first proving similar properties of (\ref{gen-corr-bcft}),
which is obvious by the properties of intertwining operators.  
Then those properties
(\ref{d-m-co}) and (\ref{D-m-co}) of $m_{cl-op}^{(l;n)}$
follow from analytic continuations.

Since $\Y^{(i)}, i=1,2,3,4$ are intertwining operators of $V$, 
we have, for $u\in V$,  
\bea
h^L(u) &=& \lim_{z\rightarrow 0} \mathbb{Y}_{cl-op}(u, z)\one_{op} \nn
&=& \lim_{z\rightarrow 0} \lim_{z_1\rightarrow 0} 
\Y^{(1)}(Y(u, z_1)\one, z)\Y^{(2)}(\one, \bar{z}) \iota_{op}(\one) \nn
&=& \lim_{z\rightarrow 0} \lim_{z_1\rightarrow 0}
Y_{V_{op}}(u,z+z_1)\Y^{(1)}(\one, z)\Y^{(2)}(\one, \bar{z}) 
\iota_{op}(\one) \nn
&=& \lim_{z\rightarrow 0} \lim_{z_1\rightarrow 0}
Y_{V_{op}}(u,z+z_1)\mathbb{Y}_{cl-op}(\one_{cl}; z, \bar{z})
\iota_{op}(\one) \nn
&=& \lim_{z\rightarrow 0} \lim_{z_1\rightarrow 0}
Y_{V_{op}}(u,z+z_1)\iota_{op}(\one)    \nn
&=& \lim_{z\rightarrow 0} \lim_{z_1\rightarrow 0} \iota_{op}(Y(u, z+z_1)\one)\nn
&=& \iota_{op}(u). 
\eea
Similarly, one can show that $h^R(u)=\iota_{op}(u), \forall u\in V$. 
Thus we have proved the 
$V$-invariant boundary condition $h^L=h^R=\iota_{op}$.

It remains to show the convergence properties
of open-closed field algebra. 
For the first convergence property ((\ref{ocfa-conv-1}) and 
(\ref{ocfa-conv-2})), one first consider
cases when $z_j^{(i)}, r_p^{(i)}, j=1,\dots, k, p=1, \dots, m$ 
in (\ref{ocfa-conv-1}) have distinct absolute values
and $z_j, r_p, j=1,\dots, l, p=1,\dots, n$ in (\ref{ocfa-conv-1})
have distinct absolute values. 
In these cases, one can express $m_{cl-op}^{(l;n)}$
and $m_{cl-op}^{(k;m)}$ as products of $\mathbb{Y}_{cl-op}$ and $Y_{op}$. 
Then by using the associativity (\ref{asso-co-op}) and that of 
open-string vertex operator algebra, it is easy to show that 
(\ref{ocfa-conv-1}) converges absolutely to 
(\ref{ocfa-conv-2}) in the required domain. 

The rest cases can all be reduced to above cases 
by using (\ref{D-m-co}) (see the proof of Theorem 2.11 in 
\cite{HKo2} for reference). 
More precisely, for $z_1^{(i)}, \dots, z_k^{(i)}$,
$z_1, \dots, z_l\in \HH$; $r_1^{(i)}>\dots >r_m^{(i)}\geq 0$, 
$r_1> \dots >r_n > 0$\footnote{If $r_n=0$, we can further introduce
another small real variable $b$, using (\ref{D-m-co}) to
move $r_n=0$ to some $r_n'>0$. We omit
the detail.}, there always exists $a\in \R_+$ small enough so that 
both of the following sets 
\bea
&\{ z_1^{(i)}+a, \dots, z_k^{(i)}+a, r_1^{(i)}+a, \dots, r_k^{(i)}+a \}, & \nn 
&\{ z_1-a, \dots, z_l-a, r_1-a,\dots, r_n-a \} &  \nonumber
\eea
are sets whose elements have distinct absolute values. 
Then (\ref{ocfa-conv-1}) equals to 
the following iterate series: 
\bea  \label{series-1}
&&\sum_{n_2} \sum_{n_1} m_{cl-op}^{(l;n)}(u_1, \dots, u_l; v_1, \dots, v_{i-1},
P_{n_2} e^{-aD_{op}} 
P_{n_1} m^{(m;k)}_{cl-op}(\tilde{u}_1, \dots, \tilde{u}_{k}; \nn
&&\hspace{2cm} 
\tilde{v}_1, \dots, \tilde{v}_{m}; z_1^{(i)}+a, \overline{z_1^{(i)}}+a, 
\dots, z_{k}^{(i)}+a, \overline{z_k^{(i)}}+a; 
r_1^{(i)}+a, \dots, r_m^{(i)}+a), \nn
&&\hspace{2cm} v_{i+1}, \dots v_n; z_1, \bar{z}_1, 
\dots, z_l, \bar{z}_l; r_1, \dots, r_n)   a_1^{n_1}a_2^{n_2}
\eea 
when $a_1=a_2=1$. 
Hence we first switch the order of above iterate sum.
Then using (\ref{m-D-op-it}) and 
the analytic extension properties of (\ref{gen-corr-bcft}),
we can easily show that the iterate series (\ref{series-1}) 
with opposite summing order is absolutely convergent
when $1\geq |a_1|, |a_2|$. Then we can switch the order of 
iterate sum in (\ref{series-1}) freely without changing 
the value of the double sum (\ref{series-1}). By (\ref{m-D-op-it}), 
we have reduced all the remaining cases to the previous cases. 

The proof of the second convergence 
property ((\ref{ocfa-conv-3}) and (\ref{ocfa-conv-4})) is 
entirely the same except that one use the associativity
(\ref{asso-co-co}) and (\ref{m-D-cl-it}). We omit the details. 

It is also clear that such open-closed field algebra over
$V$ is unique because products of $\mathbb{Y}_{cl-op}$ and $Y_{op}$ 
determine a dense subset of the set of all cases. The rest cases are
uniquely determined by continuity. 
\epf

\renewcommand{\theequation}{\thesection.\arabic{equation}}
\renewcommand{\thethm}{\thesection.\arabic{thm}}
\setcounter{equation}{0}
\setcounter{thm}{0}

\section{Operadic formulation}

In this section, we first review the notion of
$2$-colored partial operad and algebra over it. 
Then we recall the notion of Swiss-cheese partial operad
$\mathfrak{S}$, its relation to sphere partial operad $K$
and its $\C$-extension  $\tilde{\mathfrak{S}}^c$ \cite{HKo1}. 
In the end, we prove that an open-closed field algebra
over $V$ canonically gives an algebra over $\tilde{\mathfrak{S}}^c$.

\subsection{$2$-colored partial operads}

We recall the notion of $2$-colored (partial) operad 
\cite{V}\cite{Kont} and algebras over it.  
The $2$-colored operad is called {\it relative operad} 
in \cite{V} and {\it colored operad} in \cite{Kont}.

We first recall some basic notions from \cite{H4}. 
The notion of (partial) operad
can be defined in any symmetric monoidal category \cite{MSS}. In this
work, we only work in the category of sets. 
We will use the definition of (partial) operad given in \cite{H4}. 
We denote a (partial) operad as a triple 
$(\mathcal{P}, I_{\mathcal{P}}, \gamma_{\mathcal{P}})$, where
$\mathcal{P}= \{ \mathcal{P}(n) \}_{n\in \N}$ is a family of sets, 
$I_{\mathcal{P}}$ the identity element and 
$\gamma_{\mathcal{P}}$ the substitution map. Note that 
the definition of (partial) operad in \cite{H4} 
is slightly different from 
that in \cite{MSS} for the appearance of $\mathcal{P}(0)$ which 
is very important for the study of conformal field theory.   
If a triple $(\mathcal{P}, I_{\mathcal{P}}, \gamma_{\mathcal{P}})$
satisfies all the axioms of a (partial) operad except the 
associativity of $\gamma_{\mathcal{P}}$, then it is called
{\it (partial) nonassociative operad}.

We consider an important example of partial nonassociative 
operad. Let $U=\oplus_{n\in J} U_{(n)}$ be a vector space graded 
by an index set $J$
and $E_U = \{ E_U(n) \}_{n\in \N}$ a family of vector spaces, where
$$
E_{U}(n) = \hom_{\C}(U^{\otimes n}, \overline{U}).
$$
For $k,n_1,\dots, n_k\in \N$, 
$f\in E_U(k), g_i\in E_U(n_i), i=1,\dots, k$ and 
$v_j\in U, j=1,\dots, n_1+\dots +n_k$, 
\bea
&& \hspace{-1cm}
\gamma_{E_U}(f; g_1,\dots, g_k)(v_1\otimes \dots \otimes v_{n_1+\dots+n_k})
\nn
&&:= \sum_{s_1,\dots, s_k \in J} 
f(P_{s_1} g_1(v_1\otimes \dots \otimes v_{n_1}) \otimes 
\cdots  \nn
&&\hspace{3cm} \otimes 
P_{s_k} g_k(v_{n_1+\dots + n_{k-1}+1}\otimes \dots \otimes
v_{n_1+\dots + n_k}) )
\eea
is well-defined if the sum is countable and absolutely convergent. 
$\gamma_{E_U}$ is not associative in general because 
an iterate series may converge in both order but may not
converge to the same value. It is clear that
$(E_U, \id_U, \gamma_{E_U})$ is a partial nonassociative operad. 
We sometimes denote it simply as $E_U$. If $J$ is the set of equivalence
classes of irreducible modules over a group $G$ and 
$U_{(n)}$ is a direct sum of irreducible $G$-modules of the equivalence class
$n\in J$, we denote this partial nonassociative operad as $E_U^G$. 

\begin{defn}  {\rm
An {\it algebra over a partial operad 
$(\mathcal{P}, I_{\mathcal{P}}, \gamma_{\mathcal{P}})$}, 
or simply a $\mathcal{P}$-algebra, is a graded vector space $U$,
together with a partial nonassociative operad homomorphism 
$\nu: \mathcal{P}\rightarrow E_U$. 
}
\end{defn}

\begin{defn}  {\rm
Given a partial operad 
$(\mathcal{P},I_{\mathcal{P}}, \gamma_{\mathcal{P}})$, a subset $G$ of
$\mathcal{P}(1)$ is called {\it rescaling group} for $\mathcal{P}$
if
\bnu
\item For any $n\in \N$, $P_i\in G, i=0,\dots,n$ 
and $P\in \mathcal{P}(n)$, $\gamma_{\mathcal{P}}(P; P_1, \dots, P_n)$ and
$\gamma_{\mathcal{P}}(P_0;P)$ are well-defined.  
\item $I_{\mathcal{P}}\in G$ and $G$ together with the identity 
$I_{\mathcal{P}}$ and
multiplication map $G\times G \xrightarrow{\gamma_{\mathcal{P}}} G$
is a group. 
\enu
}
\end{defn}

\begin{defn}  {\rm 
A partial operad 
$(\mathcal{P}, I_{\mathcal{P}}, \gamma_{\mathcal{P}})$ is called
{\it $G$-rescalable}, if for $P_i \in \mathcal{P}(n_i), i=1,\dots,k$
and $P_0\in \mathcal{P}(k)$, then there exists 
$g_i\in G, i=1,\dots, k$ such that 
$$
\gamma_{\mathcal{P}}(\gamma_{\mathcal{P}}(P_0; g_1,\dots, g_k); 
P_1, \dots, P_k)
$$ is well-defined. 
}
\end{defn}

\begin{defn}  {\rm
An {\it algebra over a $G$-rescalable partial operad
$(\mathcal{P}, I_{\mathcal{P}}, \gamma_{\mathcal{P}})$}, 
or a {\it $G$-rescalable $\mathcal{P}$-algebra},
consists of a completely reducible $G$-modules
$U = \oplus_{n\in J} U_{(n)}$, where $J$ is the set of 
equivalence classes of irreducible $G$-modules
and $U_{(n)}$ is a direct sum of irreducible
$G$-modules of the equivalence class $n$, and a partial 
nonassociative operad homomorphism 
$\nu: \mathcal{P} \rightarrow E_U^G$ such that 
$\nu: G \rightarrow \edo U_{(n)}$ coincides 
with $G$-module structure on $U_{(n)}$. 
}
\end{defn}
We denote such algebra as $(U, \nu)$.

For $m\in \Z_+$, let $m=m_1+\dots +m_n$ be an ordered 
partition and $\sigma\in S_n$. 
The block permutation $\sigma_{(m_1,\dots, m_n)} \in S_m$ is
the permutation acting on $\{1, \dots, m\}$ by permuting $n$ intervals
of lengths $m_1, \dots, m_n$ in the same way that $\sigma$ permute
$1, \dots, n$. Let $\sigma_i\in S_{m_i}, i=1,\dots, n$, we
view the element 
$(\sigma_1, \dots, \sigma_n)\in S_{m_1}\times \dots \times S_{m_n}$
naturally as an element in $S_m$ by the canonical embedding
$S_{m_1}\times \dots \times S_{m_n} \hookrightarrow S_m$.

\begin{defn}   {\rm 
Let $(\mathcal{Q}, I_{\mathcal{Q}}, \gamma_{\mathcal{Q}})$ be an operad. 
A right module over $\mathcal{Q}$,
or a right $\mathcal{Q}$-module,  is a family of 
sets $\mathcal{P}= \{ \mathcal{P}(n) \}_{n\in \N}$ 
with actions of permutation groups, equipped with maps: 
$$
\mathcal{P}(k) \times 
\mathcal{Q}(n_1) \times \dots \times \mathcal{Q}(n_l) 
\xrightarrow{\gamma} 
\mathcal{P}(n_1+ \dots + n_k)
$$
such that 
\bnu
\item For $c\in \mathcal{P}(k)$, we have
\beq
\gamma(c; I_{\mathcal{Q}}, \dots, I_{\mathcal{Q}}) = c.
\eeq

\item $\gamma$ is associative. Namely, for $c\in \mathcal{P}(k)$, 
$d_i\in \mathcal{Q}(p_i),i=1,\dots, k$, 
$e_j\in \mathcal{Q}(q_j), j=1,\dots, p_1+\dots+p_k$, we have
\beq   \label{rel-op-asso}
\gamma(\gamma(c;d_1,\dots, d_k); e_1, \dots, e_{p_1+\dots+p_k})
= \gamma(c; f_1, \dots, f_k),
\eeq
where 
$$
f_s= \gamma_{\mathcal{Q}}(d_s; e_{p_1+\dots+p_{s-1}+1}, \dots, 
e_{p_1+\dots+p_s}).
$$

\item For $c\in \mathcal{P}(k;l)$, 
$d_i\in \mathcal{Q}(p_i),i=1,\dots, l$, $\sigma\in S_l$ and
$\tau_j\in S_{p_j}, j=1,\dots, l$, 
\bea
\gamma(\sigma(c); d_1, \dots, d_l) &=& 
\sigma_{(p_1,\dots, p_l)}( \gamma(c; d_1, \dots, d_l) ),  
\label{perm-axiom-1}\\
\gamma(c; \tau_1(d_1), \dots, \tau_l(d_l)) &=&
(\tau_1, \dots, \tau_l)(\gamma(c; d_1,\dots, d_l)).
\label{perm-axiom-2} 
\eea
\enu
Homomorphisms and isomorphisms between right $\mathcal{Q}$-modules
are naturally defined. 
}
\end{defn}

The left module over a partial operad can be similarly defined. 

\begin{defn}  {\rm
A right module $\mathcal{P}$ over a partial operad $\mathcal{Q}$, 
or a right $\mathcal{Q}$-module,
is called {\it $G$-rescalable} if $\mathcal{Q}$ is $G$-rescalable
and for any $c\in \mathcal{P}(k)$, 
$d_i\in \mathcal{Q}(n_i),i=1,\dots, k$, there exist 
$g_i\in G$, $i=1,\dots,k$ such that 
$$
\gamma(\gamma(c;g_1,\dots, g_k);  d_1, \dots, d_k )
$$
is well-defined.
} 
\end{defn}

\begin{defn} {\rm
A {\em $2$-colored operad } consists of an operad
$(\mathcal{Q}, I_{\mathcal{Q}}, \gamma_{\mathcal{Q}})$,
a family of sets $\mathcal{P}(m;n)$ equipped with an
$S_m\times S_n$-action for $m,n\in \N$, a distinguished
element $I_{\mathcal{P}}\in \mathcal{P}(1,1)$
and substitution maps $\gamma_1, \gamma_2$ given as follows: 
\bea
&&\mathcal{P}(k; l) \times \mathcal{P}(m_1; n_1) \times \dots \times
\mathcal{P}(m_k; n_k)  \nn
&& \hspace{4cm} \xrightarrow{\gamma_1}
\mathcal{P}(m_1+\ldots +m_k; l+n_1+\ldots + n_k ), \nn
&& \mathcal{P}(k; l) \times 
\mathcal{Q}(p_1) \times \dots \times \mathcal{Q}(p_l) 
\xrightarrow{\gamma_2} \mathcal{P}(k; p_1+ \dots + p_l),
\eea 
satisfying the following axioms: 
\bnu
\item The family of sets 
$\mathcal{P} := \{ \cup_{n\in \N} \mathcal{P}(m;n) \}_{m\in \N}$ 
equipped with the natural $S_m$-action on 
$\mathcal{P}(m) = \cup_{n\in \N} \mathcal{P}(m;n)$, together
with identity element $I_{\mathcal{P}}$ and substitution maps $\gamma_1$
is an operad. 

\item $\gamma_2$ gives each $\mathcal{P}(k)$ a right 
$\mathcal{Q}$-module structure for $k\in \N$. 
\enu
We denote it as $(\mathcal{P}|\mathcal{Q}, (\gamma_1,\gamma_2))$ or
$\mathcal{P}$ for simplicity.  
}
\end{defn}

\begin{rema} {\rm 
The substitution map $\gamma_1, \gamma_2$ can be combined into
a single substitution map $\gamma = (\gamma_1, \gamma_2)$:
\bea
&&\mathcal{P}(k; l) \times \mathcal{P}(m_1; n_1) \times \ldots \times
\mathcal{P}(m_k; n_k)  \times 
\mathcal{Q}(p_1) \times \dots \times \mathcal{Q}(p_l) \nn
&& \hspace{2cm} \xrightarrow{\gamma}
\mathcal{P}(m_1+\ldots +m_k; n_1+\ldots + n_k+ p_1 +\dots + p_l). 
\eea
For this reason, we also denote 
$(\mathcal{P}|\mathcal{Q}, (\gamma_1,\gamma_2))$ as
$(\mathcal{P}|\mathcal{Q}, \gamma)$.
For some examples we encounter later, $\gamma_1,\gamma_2$ are
often defined all together in terms of $\gamma$. 
}
\end{rema}

\begin{defn}  {\rm
$2$-colored partial operad is defined similarly 
as that of $2$-colored operad
except that $\gamma_1,\gamma_2$ are only partially defined and 
$\mathcal{Q},\mathcal{P}$ are partial operads and 
(\ref{rel-op-asso}) holds whenever both sides exist. 
If only the associativities of $\gamma_{\mathcal{Q}},\gamma_1,\gamma_2$
do not hold, then it is called 
{\it $2$-colored nonassociative (partial) operad}. 
}
\end{defn}

We give an important example of $2$-colored partial nonassociative 
operad. Let $J_1, J_2$ be two index sets. 
$U_1=\oplus_{n\in J_1} (U_1)_{(n)}, U_2=\oplus_{n\in J_2} (U_2)_{(n)}$ be
two graded vector spaces. Consider two families of vector spaces,
\bea
E_{U_2}(n) &=& \hom_{\C}(U_2^{\otimes n}, \overline{U}_2), \nn
E_{U_1|U_2}(m;n) &=& \hom_{\C}(U_1^{\otimes m} \otimes U_2^{\otimes n}, 
\overline{U}_1).
\eea
We denote both of the projection operators 
$U_1 \rightarrow (U_1)_{(n)}$, $U_2\rightarrow (U_2)_{(n)}$ as $P_n$ 
for $n\in J_1$ or $J_2$.  For any $f\in E_{U_1|U_2}(k;l)$,
$g_i \in E_{U_1|U_2}(m_i;n_i), i=1,\dots, k$, 
and $h_j \in E_{U_2}(p_j)$, $j=1,\dots,l$, we say that 
\bea
&&\hspace{-1cm}\Gamma(f; g_1, \dots, g_n; h_1, \dots, h_l):=\nn
&&\hspace{0.5cm} \sum_{s_1, \dots, s_k\in J_1; t_1, \dots, t_l\in J_2} 
f(P_{s_1}g_1(u_{1}^{(1)}, 
\dots, u_{m_1}^{(1)}, v_{1}^{(1)},\dots, v_{n_1}^{(1)}), \nn
&&\hspace{4cm} \dots, P_{s_k} g_k(u_{1}^{(k)}, \dots, u_{m_k}^{(k)}, 
v_{1}^{(k)}, \dots, v_{n_k}^{(k)}); \nn
&&\hspace{2cm} P_{t_1}h_1(w_1^{(1)}, \dots, w_{n_1}^{(1)}), \dots, 
P_{t_l}h_l(w_1^{(l)}, \dots, w_{n_l}^{(l)}) ) \nonumber
\eea
where $u_j^{(i)}\in U_1, v_j^{(i)}\in U_2$, 
is well-defined if 
the multiple sum is absolutely convergent. This gives arise 
to a partially defined substitution map $\Gamma$: 
\bea
&& E_{U_1|U_2}(k;l) \times E_{U_1|U_2}(m_1;n_1)\times \dots \times 
E_{U_1|U_2}(m_k;n_k)\times
E_{U_2}(p_1)\times \dots \times E_{U_2}(p_l)  \nn
&&\hspace{2cm} \xrightarrow{\Gamma}
E_{U_1|U_2}(m_1+\dots + m_k, n_1 +\dots
+n_k + p_1 +\dots +p_l).  \nonumber
\eea
$\Gamma$ does not satisfy the associativity in general. 
Let $E_{U_1|U_2}=\{ E_{U_1|U_2}(m;n) \}_{m,n\in \N}$ and
$E_{U_2} = \{ E_{U_2}(n) \}_{n\in \N}$. It is obvious that 
$(E_{U_1|U_2}|E_{U_2}, \Gamma)$
is a $2$-colored nonassociative partial operad.

Let $U_1$ be a completely reducible $G_1$-modules and 
$U_2$ a completely reducible $G_2$-modules. Namely, 
$U_1=\oplus_{n_1\in J_1} (U_1)_{(n_1)}, U_2= \oplus_{n_2\in J_2} (U_2)_{(n_2)}$
where $J_i$ is the set of equivalence classes of 
irreducible $G_i$-modules and $(U_i)_{(n_i)}$ 
is a direct sum of irreducible $G_i$-modules of the 
equivalence class $n_i$ for $i=1,2$. 
In this case, we denote $E_{U_1|U_2}$ by $E_{U_1|U_2}^{G_1|G_2}$.

\begin{defn} {\rm
A homomorphism between two $2$-colored (partial) operads 
$(\mathcal{P}_i|\mathcal{Q}_i, \gamma_i), i=1,2$ consists of 
two (partial) operad homomorphisms: 
\beq
\nu_{\mathcal{P}_1|\mathcal{Q}_1}: \mathcal{P}_1 \rightarrow \mathcal{P}_2,  
\quad\quad \mbox{and} \quad\quad
\nu_{\mathcal{Q}_1}: \mathcal{Q}_1 \rightarrow \mathcal{Q}_2  \nonumber
\eeq 
such that $\nu_{\mathcal{P}_1|\mathcal{Q}_1}: \mathcal{P}_1 \rightarrow 
\mathcal{P}_2 $, where $\mathcal{P}_2$ is a
right $\mathcal{Q}_1$-module by $\nu_{\mathcal{Q}_1}$, 
is also a right $\mathcal{Q}_1$-module homomorphism.
}
\end{defn}

\begin{defn} {\rm
An {\it algebra over a $2$-colored partial operad 
$(\mathcal{P}|\mathcal{Q}, \gamma)$}, or a 
{\it $\mathcal{P}|\mathcal{Q}$-algebra} 
consists of 
two graded vector spaces $U_1, U_2$ and 
a homomorphism $(\nu_{\mathcal{P}|\mathcal{Q}}, \nu_{\mathcal{Q}})$ 
from $(\mathcal{P}| \mathcal{Q}, \gamma)$ to 
$(E_{U_1|U_2}|E_{U_2}, \Gamma)$. We denote this algebra as 
$(U_1|U_2, \nu_{\mathcal{P}|\mathcal{Q}}, \nu_{\mathcal{Q}})$. 
}
\end{defn}

\begin{defn}  {\rm
If a $2$-colored partial operad $(\mathcal{P}|\mathcal{Q}, \gamma)$
is so that $\mathcal{P}$ 
is a $G_1$-rescalable partial operad
and a $G_2$-rescalable right $\mathcal{Q}$-module, then 
it is called {\it $G_1|G_2$-rescalable}. 
}
\end{defn}

\begin{defn} {\rm
A {\it $G_1|G_2$-rescalable $\mathcal{P}|\mathcal{Q}$-algebra}
$(U_1|U_2,\nu_{\mathcal{P}|\mathcal{Q}}, \nu_{\mathcal{Q}})$ 
is a $\mathcal{P}|\mathcal{Q}$-algebra so that 
$\nu_{\mathcal{P}|\mathcal{Q}}: \mathcal{P}\rightarrow E_{U_1|U_2}^{G_1|G_2}$
and $\nu_{\mathcal{Q}}: \mathcal{Q} \rightarrow E_{U_2}^{G_2}$; 
moreover, $\nu_{\mathcal{P}|\mathcal{Q}}: G_1\rightarrow \edo U_1$
coincides with the $G_1$-module structure on $U_1$ 
and $\nu_{\mathcal{Q}}: G_2\rightarrow \edo U_2$ coincides with
the $G_2$-module structure on $U_2$. 
}
\end{defn}

\subsection{Swiss-cheese partial operad $\mathfrak{S}$}

A {\it disk with 
strips and tubes of type $(m_-, m_+ ; n_-, n_+)$} 
($m_-, m_+, n_-, n_+\in \N$) 
is a disk $S$ with the following additional data:
\bnu
\item $m_-+m_+$ distinct ordered punctures 
$p^{B}_{-m_-}, \dots , p^{B}_{-1}, p^{B}_{1}, \dots, p^{B}_{m_+}$ 
(called {\it boundary punctures}) 
on $\partial S$ (the boundary of $S$), where
$p^{B}_{-m_-}, \dots, p^{B}_{-1}$ are negatively oriented and 
$p^{B}_{1}, \dots, p^{B}_{m_+}$ are positively oriented, 
together with local coordinates:
$$
(U^{B}_{-m_-}, \varphi^{B}_{-m_-}), \dots, (U^{B}_{-1}, \varphi^{B}_{-1});
\quad (U^{B}_{1}, \varphi^{B}_{1}), \dots, (U^{B}_{m_+}, \varphi^{B}_{m_+}),
$$
where $U^{B}_{i}$ is a neighborhood of $p_i^{B}$ and 
$\varphi^{B}_{i}: U^{B}_{i} \to \tilde{\HH}$ is an 
analytic map which vanishes at $p^{B}_{i}$ and maps
$U^{B}_{i} \cap \partial S$ to $\R$,
for each $i=-m_-, \dots, -1$, $1, \dots, m_+$.

\item $n_-+n_+$ distinct ordered points 
$p^{I}_{-n_-}, \dots , p^{I}_{-1}, p^{I}_{1}, \dots , p^{I}_{n_+}$ 
(called {\it interior punctures}) 
in the interior of $S$, where
$p^{I}_{-n_-}, \dots, p^{I}_{-1}$ are negatively oriented and 
$p^{I}_{1}, \dots, p^{I}_{n_+}$ are positively oriented,
together with local coordinates: 
$$
(U^{I}_{-n_-}, \varphi^{I}_{-n_-}), \dots, (U^{I}_{-1}, \varphi^{I}_{-1});
\quad (U^{I}_{1}, \varphi^{I}_{1}), \dots, (U^{I}_{1}, \varphi^{I}_{n_+}),
$$
where $U^{I}_{j}$ is a local neighborhood of $p^{I}_{j}$ and 
$\varphi^{I}_{j}: U^{I}_{j} \to \C$ is an 
analytic map which vanishes at $p^{I}_{i}$
for each $j=-n_-, \dots, -1, 1, \dots, n_+$.  

\enu

Two disks with strips and tubes are conformal equivalent
if there exists between them a biholomorphic map
which maps punctures to punctures and preserves the order of
punctures and the germs of local
coordinate maps. We denote the moduli space of 
the conformal equivalence classes of disks with strips and tubes 
of type $(m_-, m_+; n_-, n_+)$ as
$\mathbb{S}(m_-, m_+ | n_-, n_+)$. 
The structure on $\mathbb{S}(m_-, m_+ | n_-, n_+)$ will be discussed
in \cite{Ko2}. 
In this work, we are only interested in disks with strips and tubes 
of types $(1, m_+; 0, n_+)$ for $n, l\in \N$. For simplicity, we
denote $\mathbb{S}(1, m_+ | 0, n_+)$ by $\Upsilon(m_+; n_+)$. 
For such disks, we label the only negatively 
oriented boundary puncture as the $0$-th boundary puncture.


We can choose a 
canonical representative for each conformal equivalence class
in $\Upsilon(m_+; n_+)$
just as we did for disks with strips \cite{HKo1} 
and sphere with tubes \cite{H4}. 
More precisely, for a disk with strips and tubes of type
$(1, m_+; 0, n_+)$ where $m_+>0$, 
we first use a conformal map $f$ to map 
the disk to $\hat{\HH}$. Then 
we use an automorphism of $\hat{\HH}$ to move the only negatively 
oriented puncture (the $0$-th puncture) to $\infty$ and
the smallest $r_{k}$ to $0$, 
and fix the local coordinate map $f_{0}$ at $\infty$
to be so that $\lim_{w\rightarrow \infty} wf_{0}(w)=-1$.
As a consequence, the canonical representative 
of a generic conformal equivalent
class of disk with strips and tubes $Q\in \Upsilon(m_+; n_+)$ is
a disk $\hat{\HH}$, together with 
a negatively oriented boundary puncture at $\infty\in \hat{\HH}$
and local coordinate map given by
$$
f_0^B(w) = - e^{\sum_{j=1}^{\infty} -B_j^{(0)} w^{-j+1}\frac{d}{dw} } 
\, \frac{1}{w},
$$
where $B_j^{(0)}\in \R$, 
and positively oriented boundary punctures at 
$r_1,\dots r_{m_+} \in \R_+\cup \{ 0\}$ ($r_k=0$) and
local coordinate maps given by
$$
f_i^B(w) = e^{\sum_{j=1}^{\infty} B_j^{(i)} x^{j+1}\frac{d}{dx} } 
(b_0^{(i)})^{x\frac{d}{dx}} x \lbar_{x=w-r_i}, i=1, \dots, m_+, 
$$
where $B_j^{(i)}\in \R, b_0^{(i)} \in \R_+$, 
and positively oriented interior punctures at 
$z_1, \dots, z_{n_+}\in \HH$ with local coordinate maps given by 
$$
f_i^{I}(w) = e^{\sum_{j=1}^{\infty} A_j^{(i)} x^{j+1}\frac{d}{dx} } 
(a_0^{(i)})^{x\frac{d}{dx}} x \lbar_{x=w-z_i}, i=1, \dots, n_+.
$$
where $A_j^{(i)}\in \C$ and $a_0^{(i)}\in \C^{\times}$.

Let $\Pi_{\R}$ ($\Pi_{\C}$) be the set of sequences of real numbers 
(complex numbers) $\{ C_j\}_{j=1}^{\infty}$ such that 
$e^{\sum_{j>0}C_j x^{j+1}\frac{d}{dx}}x$ as a power series
converges in some neighborhood of $0$. We define
$$
\tilde{\Lambda}^n:=\{(r_1,\dots, r_n)|\,\, \exists \sigma\in S_n,
\,\, r_{\sigma(1)}>\dots >r_{\sigma(n)}=0\}.
$$ 
Using the data on the canonical representative of $Q$,
we denote $Q\in \Upsilon(m_+; n_+)$ as follows
\bea \label{gen-elm-moduli}
& & \hspace{-1cm} 
[r_1, \ldots, r_{m_+-1}; B^{(0)}, (b_{0}^{(1)},B^{(1)}), \ldots, 
 (b_{0}^{(m_+)},B^{(m_+)}) | \nn
& & \hspace{2cm} z_1,\ldots, z_{n_+}; (a_{0}^{(1)},A^{(1)}), \ldots, 
 (a_{0}^{(n_+)},A^{(n_+)}) ] ,
\eea
where $(r_1, \dots, r_{m_+})\in \tilde{\Lambda}^{m_+}$
and $b_0^{(i)}\in \R_+, a_0^{(j)}\in \C^{\times}$,
and $B^{(i)}= \{ B_j^{(i)} \}_{j=1}^{\infty} \in \Pi_{\R}$
and $A^{(p)}=\{ A_q^{(p)} \}_{q=1}^{\infty} \in \Pi_{\C}$
for all $i=1, \ldots, m_+, p=1,\ldots, n_+$.

Using notations (\ref{Lambda}) and (\ref{M-HH}), 
for $m_+ > 0$ and $n_+\in \N$, we can express the moduli space of 
disks with strips and tubes of type $(1, m_+; 0, n_+)$ as follows:
$$
\Upsilon(m_+; n_+)= \tilde{\Lambda}^{m_+-1}\times \Pi_{\R} \times 
(\R_+ \times \Pi_{\R})^{m_+} \times M_{\HH}^{n_+}
\times (\C^{\times} \times \Pi_{\C})^{n_+}.
$$ 
For $m_+=0, n_+\in \N$, we used automorphism of $\hat{\HH}$ 
to fix $B_1^{(0)}=0$. Hence, we have
$$
\Upsilon(0; n_+)=\{ B^{(0)} \in \Pi\;|\; B_1^{(0)}=0\}\times 
(\C^{\times} \times \Pi_{\C})^{n_+}.
$$ 

Note that 
$\Upsilon(m_+; 0)$ is nothing but $\Upsilon(m_+)$ 
introduced and studied in \cite{HKo1}. 
$\Upsilon:=\{ \Upsilon(n) \}_{n\in \N}$ is a partial operad
of disks with strips \cite{HKo1}. The identity $I_{\Upsilon}$ 
is an element of $\Upsilon(1; 0)$. Also, for $m_+, n_+\in \N$,
$S_{n_+}$ acts on $\Upsilon(m_+; n_+)$ 
in the obvious way.
Let $\mathfrak{S}(m)=\cup_{n_+\in \N}\Upsilon(m_+; n_+)$ for $m\in \N$,
and $\mathfrak{S} = \cup_{m\in \N} \mathfrak{S}(m)$. 

There are two kinds of sewing operations on $\mathfrak{S}$,
The first kind is called {\it boundary sewing operation} which
sews the positively oriented boundary puncture in the first
disk with a negatively oriented
boundary puncture in the second disk. The second is 
called {\it interior sewing operation} 
which sews a positively oriented interior 
puncture in a disk with a negatively 
oriented puncture in a sphere with tubes. 
We describe these sewing operations more precisely. 
Let us consider
$P\in \Upsilon(m_+; n_+)$ and $Q\in \Upsilon(p_+; q_+)$ 
($Q \in K(p_+)$) for boundary sewing operations 
(for interior sewing operations).  
Let $B^r$ ($\bar{B}^r$) denote the open  
(closed) ball in $\C$ center at $0$ with radius 
$r$, $\varphi_i$ the germs of local coordinate map at $i$-th boundary
(interior) puncture $p$ of $P$, 
and $\psi_0$ the germs of local coordinate map
at $0$-th puncture $q$ of $Q$. Then we say that the $i$-th strip 
of $P$ can be sewn with the $0$-th strip (tube) of $Q$ if there is a 
$r\in \R_+$ such that $p$ and $q$ are the only punctures in 
$\varphi_i^{-1}(\bar{B}^r)$ and $\psi_0^{-1}(\bar{B}^{1/r})$ respectively. 
A new disk with strips and tubes in $\Upsilon(m_++p_+-1, n_++q_+)$
($\Upsilon(m_+, n_++p_+-1)$),
denoted as $ P_{^{i}}\infty_{^{0}}^{B} Q$ ($ P _{^i} \infty_{^0}^I Q$), 
is obtained by cutting out $\varphi_i^{-1}(B^r)$ and 
$\psi_0^{-1}(B^{1/r})$ from $P$ and $Q$ respectively, and then 
identifying the boundary of 
$\varphi_i^{-1}(\bar{B}^r)$ 
and $\psi_0^{-1}(\bar{B}^{1/r})$ via the map 
$$
\psi^{-1} \circ J_{\hat{\HH}} \circ \varphi_i
$$
where $J_{\hat{\HH}}: w \rightarrow -\frac{1}{w}$.

The boundary sewing operations and interior sewing operations
induce the following partially defined substitution maps: 
\bea
&&\Upsilon(k; l) \times \Upsilon(m_1; n_1) \times \ldots \times
\Upsilon(m_k; n_k)  \times K(p_1) \times \dots \times K(p_l) \nn
&& \hspace{2cm} \xrightarrow{\gamma}
\Upsilon(m_1+\ldots +m_k; n_1+\ldots + n_k+ p_1 +\dots + p_l). 
\eea

The following Proposition is clear. 
\begin{prop} {\rm
$(\mathfrak{S}|K, \gamma)$ is a 
$\R_+|\C^{\times}$-rescalable $2$-colored partial operad. 
} \end{prop}

This $\R_+|\C^{\times}$-rescalable $2$-colored partial operad
$(\mathfrak{S}|K, \gamma)$
is a generalization of Voronov's 
Swiss-cheese operad \cite{V}. So we will call it
{\it Swiss-cheese partial operad} and sometimes denote it by
$\mathfrak{S}$ for simplicity. The relation between 
Swiss-cheese operad and Swiss-cheese partial operad is an analogue
of that between 
little disk operad and sphere partial operad \cite{H4}.

\subsection{Sewing equations and the doubling map $\delta$}

We are interested in finding the canonical representative of 
disk with strips and tubes obtained by  
sewing two such disks or sewing a disk with a sphere. 
In the case of sphere partial operad $K$, 
such canonical representatives
were obtained by Huang \cite{H1}\cite{H2}\cite{H4} 
by solving the so-called sewing equation. 
Similarly, the canonical representatives of disks 
with strips and tubes
obtained by two types of sewing operations can also be determined 
by solving two types of sewing equations.

We start with the boundary sewing operations.  
For $Q\in \mathfrak{S}$, we denoted the canonical representative
of $Q$ as $\Sigma_Q$. 
Let $P\in \Upsilon(m;n)$ and $Q\in \Upsilon(p;q)$. 
Let $g_0$ be the local coordinate map at $\infty \in \Sigma_Q$ 
and $f_i$ be that at $z_i\in \Sigma_P, 1\leq i\leq m$. 
We assume that $P _{^i}\infty_{^0}^B Q$ exists. Then 
the canonical disk $\Sigma_{P _{^i}\infty_{^0}^B Q}$ 
can be obtained by solving the following sewing equation:
\begin{equation} \label{sew-equ-B}
F_{(1)}^B(w) = F_{(2)}^B 
\left( g_0^{-1}\left( \frac{-1}{f_i(w)} \right) \right)
\end{equation}
where $F_{(1)}^B$ is a conformal map from an open neighborhood of
$\infty\in \Sigma_P$
to an open neighborhood of 
$\infty \in \Sigma_{P_{^i}\infty_{^0}^BQ}$, and 
$F_{(2)}^B$ is a conformal map from an open neighborhood of 
$0\in \Sigma_Q$ to an open neighborhood of 
$0\in \Sigma_{P_{^i}\infty_{^0}^BQ}$, with 
the following normalization conditions: 
\beq  \label{norm-equ-B} 
F_{(1)}^B(\infty) = \infty,  \hspace{1cm}
F_{(2)}^B(0)  = 0,  \hspace{1cm}
\lim_{w\rightarrow \infty} \frac{F_{(1)}^B (w)}{w} = 1.
\eeq
It is easy to see that 
the solution of (\ref{sew-equ-B}) and (\ref{norm-equ-B}) is
unique. Notice that $F_{(1)}^B$, $F_{(2)}^B$, $f_0$ and $f_i$ in 
(\ref{sew-equ-B}) are all real analytic.

Similarly, let $P\in \Upsilon(m;n)$ and $Q\in K(p)$. 
We denote the canonical sphere with tubes of $Q$ as $\Sigma_Q$. 
Let $g_0$ be the local coordinate map at 
$\infty \in Q$ and $f_i$ be that at $z_i\in P, 1\leq i\leq n$. 
We assume that $P _{^i}\infty_{^0}^I Q$ exists. Then 
$\Sigma_{P _{^i}\infty_{^0}^I Q}$ can be obtained by solving 
the following sewing equation:
\begin{equation} \label{sew-equ-I}
F_{(1)}^I(w) = F_{(2)}^I \left( g_0^{-1}\left( 
\frac{-1}{f_i(w)} \right) \right)
\end{equation}
where $F_{(1)}^I$ is a conformal map from an open neighborhood of
$\hat{\R}$ in  $\Sigma_P$
to an open neighborhood of $\hat{\R}$ in 
$\Sigma_{P _{^i}\infty_{^0}^I Q}$, and 
$F_{(2)}^I$ is a conformal map from an open neighborhood of $0$
in $\Sigma_Q$ to an open subset of 
$\HH \subset \Sigma_{P _{^i}\infty_{^0}^I Q}$, with 
the following normalization conditions:
\beq  \label{norm-equ-I}
F_{(1)}^I(\infty) = \infty,   \hspace{1cm}
F_{(1)}^I(0) = 0,  \hspace{1cm}
\lim_{w\rightarrow \infty} \frac{F_{(1)}^I(w)}{w} = 1.
\eeq
It is easy to see that the solution of (\ref{sew-equ-I}) and
(\ref{norm-equ-I}) is unique as well.

Notice that $F_{(1)}^I$ is real analytic because it 
maps $\R$ to $\R$. Hence, $F_{(1)}^I$
is also the unique solution for the following equation 
\beq   \label{sew-equ-1}
F_{(1)}^I(w) = \overline{F_{(2)}^I} 
\left( \overline{g_0^{-1} } \left( 
\frac{-1}{\overline{f_i}(w)} \right) \right)
\eeq
with the same normalization condition (\ref{norm-equ-I}).

In sphere partial operad, 
we define a complex conjugation map $\text{Conj}: K \rightarrow K$ 
as follows
\bea
& & \hspace{-1cm} \text{Conj}: 
(z_1, \ldots, z_{n-1}; A^{(0)}, (a_0^{(1)}, A^{(1)}), \ldots, 
(a_0^{(n)}, A^{(n)}) ) \nn
& & \hspace{1cm} \longmapsto 
(\bar{z}_1, \ldots, \bar{z}_{n-1}; \bar{A}^{(0)}, 
(\bar{a}_0^{(1)}, \bar{A}^{(1)}), \ldots, 
(\bar{a}_0^{(n)}, \bar{A}^{(n)}) ).
\eea
For simplicity, we denote 
$\text{Conj}(Q)$ as $\bar{Q}$ for $Q\in K$. 

\begin{prop} \label{conj-auto-prop-K}
{\rm Conj} is an partial operad automorphism of $K$. 
\end{prop}
\pf
It is clear that 
$\text{Conj}((\mathbf{0},(1,\mathbf{0}))) = (\mathbf{0},(1,\mathbf{0}))$ and
$\text{Conj}$ is equivariant with respect to the action of 
permutation group. Moreover, $\text{Conj}$ is obviously bijective. 
It only remains to show that, for $1\leq i\leq m$,    
\begin{equation} \label{conj-auto-K}
\overline{ P _{^i}\infty_{^0} Q } = \bar{P} _{^i}\infty_{^0} \bar{Q}
\end{equation}
for any pair $P\in K(m), Q\in K(n)$ such that 
$P _{^i}\infty_{^0} Q$ exists.

Let $f_i$ be the local coordinate map at $i$-th
puncture in $P$ and $g_0$ be that at $\infty$ in $Q$. 
Then the local coordinate map at $i$-th puncture in $\bar{P}$ is
$\bar{f}_i$ and that at $\infty$ 
in $\bar{Q}$ is $\bar{g}_0$.  Also notice that 
the sewing equation and normalization equation 
for sphere partial operad (\cite{H4}) are 
the same as the equation (\ref{sew-equ-B}) and (\ref{norm-equ-B}). 
Let $F_{sphere}^{(1)}, F_{sphere}^{(2)}$ be the solution of
(\ref{sew-equ-B}) and (\ref{norm-equ-B}) for the sphere with tubes , 
then $\bar{F}_{sphere}^{(1)}, \bar{F}_{sphere}^{(2)}$ also satisfy
the same normalization condition and the following sewing equation
$$
\bar{F}_{sphere}^{(1)}(w) = \bar{F}_{sphere}^{(2)} 
\left( \overline{g_0^{-1}}\left( 
\frac{-1}{\bar{f}_i(w)} \right) \right),
$$
which is the sewing equation for
$\bar{P} _{^i}\infty_{^0} \bar{Q}$. 
Using the explicit formula ((A.6.1)-(A.6.5) in \cite{H4}) 
of the moduli $\bar{P} _{^i}\infty_{^0} \bar{Q}$ in terms of 
$\bar{F}_{sphere}^{(1)}, \bar{F}_{sphere}^{(2)}$, 
$\bar{f}_i, \bar{g}_0$, one can
easily see that (\ref{conj-auto-K}) is true. 
\epf

There is a canonical doubling map
$\delta: \mathfrak{S} \rightarrow K$ defined as follows. 
Let $Q\in \Upsilon(n;l)$ with form
\bea
Q&=& [\, r_1, \ldots, r_{n-1}; B^{(0)}, (b_{0}^{(1)},B^{(1)}), \ldots, 
 (b_{0}^{(n)},B^{(n)}) | \nn 
& & \hspace{2cm} z_1,\ldots, z_l; (a_{0}^{(1)},A^{(1)}), \ldots, 
 (a_{0}^{(l)},A^{(l)}) \,].  \nonumber
\eea 
Then 
\bea   \label{delta-map}
&&\hspace{-0.8cm} 
\delta(Q)= (z_1, \ldots, z_l, \bar{z}_1, \dots, \bar{z}_l, 
r_1, \ldots, r_{n-1} ;   
(a_{0}^{(1)},A^{(1)}),  \ldots, (a_{0}^{(l)},A^{(l)}),  \nn
& & \hspace{0.5cm} 
(\overline{a_{0}^{(1)}},\overline{A^{(1)}}), 
\dots, (\overline{a_{0}^{(l)}},\overline{A^{(l)}});
B^{(0)}, (b_{0}^{(1)},B^{(1)}), \ldots, (b_{0}^{(n)},B^{(n)}) ). 
\eea 

\begin{prop}  
Let $P\in \Upsilon(m;n)$ and $Q\in \Upsilon(p;q)$. 
Assume that $P _{^i}\infty_{^0}^B Q $ exists. We have
\beq
\delta( P _{^i}\infty_{^0}^B Q ) = \delta(P) _{^{2n+i}}\infty_{^0} 
\delta(Q). 
\eeq
\end{prop}
\pf
Since $F_{(1)}^B$, $F_{(2)}^B$, $g_0^{-1}$ and $f_i$ in 
(\ref{sew-equ-B}) are all real analytic, 
every solution of (\ref{sew-equ-B})
and (\ref{norm-equ-B}) for disk with strips and tubes
can be extended to a solution of the same sewing equation for
sphere with tubes by Schwarz's reflection principle. 
Then the Proposition follows immediately from this fact. 
\epf

Given a canonical disk with strips and tubes $\Sigma_Q$
corresponding to moduli $Q\in \mathfrak{S}$, we consider its
complex conjugation, denoted as $\overline{\Sigma_Q}$.
$\overline{\Sigma_Q}$ is the lower half plane $\hat{\overline{\HH}}$
together with 
the same boundary punctures and local coordinate maps as those in 
$\Sigma_Q$, and interior punctures which are the complex conjugation 
of the interior punctures in $\Sigma_Q$ with local coordinate maps
being the complex conjugation of those in $\Sigma_Q$. We can denote
it as 
\bea \label{Sigma-complex-conj}
& & [r_1, \ldots, r_{m_+-1}; B^{(0)}, (b_{0}^{(1)},B^{(1)}), \ldots, 
 (b_{0}^{(m_+)},B^{(m_+)}) | \nn
& & \hspace{2cm} \bar{z}_1,\ldots, \bar{z}_{n_+}; 
(\overline{a_{0}^{(1)}}, \overline{A^{(1)}}), \ldots, 
 (\overline{a_{0}^{(n_+)}},\overline{A^{(n_+)}}) ]
\eea
where $z_1, \dots, z_{n_+}\in \HH$.

\begin{lemma} \label{lemma-sew-I-conj}
Let $P\in \Upsilon(m;n)$ and $Q\in K(p)$. Assume that 
$P _{^i}\infty_{^0}^I Q$ exists for $1\leq i\leq m$. 
Then $\overline{\Sigma_P} _{^i}\infty_{^0}^I \bar{Q}$ also exists and 
\beq  \label{sew-I-conj}
\overline{\Sigma_{P _{^i}\infty_{^0}^I Q}} \cong 
\overline{\Sigma_P} \, _{^i}\infty_{^0}^I \, \bar{Q}. 
\eeq
\end{lemma}
\pf
We can choose a canonical representative of 
$\overline{\Sigma_P} \, _{^i}\infty_{^0}^I \, \bar{Q}$
as a lower half plane $\hat{\overline{\HH}}$
with two punctures at $\infty$ and 
$0$ and the local coordinate map $g_0$ at $\infty$ being 
so that $\lim_{w\rightarrow \infty} wg_0(w)=-1$. We denote 
such representative of 
$\overline{\Sigma_P} \, _{^i}\infty_{^0}^I \, \bar{Q}$ (a
lower half plane) as $\Sigma_1$. 
$\Sigma_1$ can be obtained by solving the sewing equation 
\beq  \label{sew-equ-2}
G_{(1)}^I(w) = G_{(2)}^I 
\left( \overline{g_0^{-1} } \left( 
\frac{-1}{\overline{f_i}(w)} \right) \right),
\eeq
where $G_{(1)}^{I}$ is a conformal map from a neighborhood of 
$\hat{\R} \subset \overline{\Sigma}_P$ to a neighborhood of 
$\hat{\R} \subset \Sigma_1$, and $G_{(2)}^I$ is a
conformal map from a neighborhood of $0\in \overline{\Sigma_Q}$
to an open subset of $\overline{\HH}\subset \Sigma_1$, 
with the normalization equation (\ref{norm-equ-I}). 
By comparing (\ref{sew-equ-2}) with (\ref{sew-equ-1}), 
we see that the unique solution 
$F_{(1)}^I$ (real analytic) and $\overline{F_{(2)}^I}$ 
of (\ref{sew-equ-1}) and 
(\ref{norm-equ-I}) exactly gives the unique solution 
$G_{(1)}^I$ and $G_{(2)}^I$of 
(\ref{sew-equ-2}) and (\ref{norm-equ-I}).
Hence we have $\Sigma_1 = \overline{\Sigma_{P _{^i}\infty_{^0}^I Q}}$. 
\epf

\begin{prop} \label{cardy-double-prop}
Let $P\in \Upsilon(n;l), Q\in K(m)$ and $P _{^i}\infty_{^0}^I Q$
exists for $1\leq i\leq l$. Then 
\begin{equation} \label{double-equ}
\delta (P _{^i}\infty_{^0}^I Q ) = 
(\delta(P) _{^i}\infty_{^0} Q) _{^{l+m-1+i}}\infty_{^0} \bar{Q}
\end{equation}
\end{prop}
\pf
Now we first consider the right hand side of (\ref{double-equ}).
We denote the canonical representative of 
any $R\in K(l)$ as $\Sigma_R$. Then $\Sigma_{\delta(P)}$ 
can be viewed as a union of the closure of
upper half plane, denoted as $U_+$, 
and the closure of lower half plane, denoted as $U_-$. 
Let $\Sigma_+$ be the Riemann surface obtained by sewing
$U_+$ with $Q$, and $\Sigma_-$ the Riemann surface obtained by
sewing $U_-$ with $\bar{Q}$.  
By identifying the real line in $U_+\subset \Sigma_+$ 
with the real line in $U_-\subset \Sigma_-$ using identity map,
we obtain the surface $\Sigma_+ \# \Sigma_-$ which is
isomorphic to the canonical sphere with tubes
$\Sigma_{(\delta(P) _{^i}\infty_{^0} Q) _{^{l+m-1+i}}\infty_{^0} \bar{Q}}$.

Both $\Sigma_+$ and $\Sigma_-$ are
disks with strips and tubes. Since $U_+ = \Sigma_P$, 
there is a unique biholomorphic map $f$ from 
$\Sigma_+$ to the canonical disk with strips and tubes 
$\Sigma_{Q _{^i}\infty_{^0}^I P }$. 
Similarily, because $U_-=\overline{\Sigma_P}$, 
by Lemma \ref{lemma-sew-I-conj}, 
there is a unique biholomorphic map $g$ from $\Sigma_-$ 
to the canonical disk with strips and tubes
$\overline{\Sigma_{P _{^i}\infty_{^0}^I Q}}$.  
The restriction of $f$ on the neighborhood of 
$\hat{\R}=\partial U_+$ is nothing but the unique  
real analytic map 
$F_{(1)}^I$ satisfying (\ref{sew-equ-I}) and (\ref{norm-equ-I}). 
Meanwhile, the restriction of $g$ on a neighborhood of 
$\hat{\R}=\partial U_-$ is nothing but the same 
$F_{(1)}^I$ satisfying (\ref{sew-equ-1}) and (\ref{norm-equ-I}). 
So we must have 
$f|_{\hat{\R}}=g|_{\hat{\R}}$, which further implies that 
$f^{-1}|_{\hat{\R}} = g^{-1}|_{\hat{\R}}$. 
Hence, $f^{-1}$ can be extended to a
biholomorphic map from $\Sigma_{\delta (P _{^i}\infty_{^0}^I Q)}$
to $\Sigma_+ \# \Sigma_-$. 

Therefore $\Sigma_{\delta (P  _{^i}\infty_{^0}^I Q )}$ 
must be biholomorphic to 
$\Sigma_{(\delta(P) _{^i}\infty_{^0} Q) _{^{l+m-1+i}}\infty_{^0} \bar{Q}}$. 
Since they are both canonical representatives of sphere with tubes, 
we must have the equality (\ref{double-equ}). 
\epf

\begin{rema}
{\rm 
Proposition \ref{cardy-double-prop} is nothing but 
the doubling trick \cite{A}\cite{C1} stated in our partial operad language. 
It also implies that the bulk theories in an open-closed conformal field theory
must contain both chiral parts and anti-chiral parts. 
}
\end{rema}

\begin{cor} 
For $P\in \delta(\Upsilon(n,l)) \subset K(n+2l)$ as in (\ref{delta-map}) 
and $Q\in K(m)$, the sewing operations:
$(P _{^i}\infty_{^0} Q) _{^{l+m-1+i}}\infty_{^0} \bar{Q}$ 
for $1\leq i \leq l$ define an action of the diagonal 
$\{ (Q, \bar{Q}) \in K\times K \}$ 
of $K\times K$ on $\delta(\mathfrak{S})$.    
\end{cor}




\subsection{The $\C$-extensions of $\mathfrak{S}$}

In order to study open-closed conformal field theories with
nontrivial central charges, we need study the $\C$-extensions of 
Swiss-cheese partial operad $\mathfrak{S}$.

For $c\in \C$, 
let $\tilde{K}^c$ be the $\frac{c}{2}$-th power of determinant line
bundle over $K$ \cite{H4}. 
We denote the pullback line bundle over $\mathfrak{S}$ 
through the doubling map $\delta$ as $\tilde{\mathfrak{S}}^{c}$. 
$\delta$ can certainly be extended to a map 
on $\tilde{\mathfrak{S}}^{c}$. We still denote it as $\delta$. 
For any $n\in \N$, the restrictions of the sections 
$\psi_{n+2l}$ of $\tilde{K}^{c}(n+2l)$ 
for $l\in \N$ to $\Upsilon(n; l)$
gives a section of $\tilde{\mathfrak{S}}^{c}(n)$ and 
we shall use $\psi_{n}^{\mathfrak{S}}$ to denote this section. 
It is clear that
$\tilde{\Upsilon}^c$, the $\C$-extension of the partial operad of
disks with strips, is the pullback bundle of the inclusion map 
$\Upsilon \hookrightarrow \mathfrak{S}$.

The boundary sewing operations in $\tilde{\mathfrak{S}}^{c}$ are
naturally induced from the sewing operations of $\tilde{K}^{c}$. 
We denote the boundary sewing operations in
$\tilde{\mathfrak{S}}^{c}$ 
as $_{^i}\widetilde{\infty}_{^0}^B$.\footnote{Since we always work with a
fixed $c\in \C$, it is convenient to make 
the dependence on $c$ implicit in some notations.} 
More explicitly, 
let $P \in \Upsilon(n;l)$ and 
$Q\in \Upsilon(m;k)$ be so that $P _{^i}\infty_{^0}^B Q$
exists. Let $\tilde{P}, \tilde{Q}$ be elements in the fiber over 
$P$ and $Q$ respectively. We define
\begin{equation} \label{B-sew-ext}
\tilde{P} _{^i}\widetilde{\infty}_{^0}^B \tilde{Q} := \delta^{-1} 
( \delta(\tilde{Q}) _{^{2l+i}}\widetilde{\infty}_{^0} \tilde{Q}),
\end{equation}
where $\delta^{-1}$ is defined on the image of $\delta$.

We would also like to lift an interior sewing operation 
$_{^i}\infty_{^0}^I$ to a sewing operation 
between an element in $\tilde{\mathfrak{S}}^c$ and an element in 
$\tilde{K}^c \otimes \overline{\tilde{K}^{\bar{c}}}$.
We still call it interior sewing operation 
and denote it as $_{^i}\widetilde{\infty}_{^0}^I$. 
Let $P \in \Upsilon(n;l), Q\in K(m)$ such that 
$P _{^i}\infty_{^0}^I Q$ exists. Let $\tilde{P},\tilde{Q}$ 
be elements in the fibers over $P$ and $Q$ respectively. 
Let $\psi_m\otimes \bar{\psi}_m$ be the canonical section on 
$\tilde{K}^c \otimes \overline{\tilde{K}^{\bar{c}}}(m)$. Then we have
$\tilde{Q}= \lambda \psi_m\otimes \bar{\psi}_m(Q) $ for
some $\lambda \in \C$. 
Then we define $\tilde{P} _{^i}\widetilde{\infty}_{^0}^I \tilde{Q}$ by
\begin{equation} \label{I-sew-ext}
\tilde{P} _{i}\widetilde{\infty}_{^0}^I \tilde{Q} := \delta^{-1}
( (\delta(\tilde{P}) _{^i}\widetilde{\infty}_{^0}  \lambda\psi_m(Q) )
_{^{l+m-1+i}}\widetilde{\infty}_{^0} \psi_m(\bar{Q}) )
\end{equation}

The following Lemma shows that the 
interior sewing operations are associative. 
\begin{lemma}
Let $P \in \Upsilon(n;l)$, $Q_1\in K(m_1), Q_2\in K(m_2)$ 
and $\tilde{P}, \tilde{Q}_1,
\tilde{Q}_2$ be elements in fibers of line bundles 
$\tilde{\mathfrak{S}}^c$ and 
$\tilde{K}^c \otimes \overline{\tilde{K}^{\bar{c}}}$ over the base points
$P, Q_1,Q_2$ respectively. Let $1\leq i \leq l$ and $1\leq j \leq m_1$. 
Then we have
\begin{equation} \label{double-ext-equ}
(\tilde{P} _{^i}\widetilde{\infty}_{^0}^I \tilde{Q}_1 ) _{^{i+j-1}}
\widetilde{\infty}_{^0}^I \tilde{Q}_2
= \tilde{P} _{^i}\widetilde{\infty}_{^0}^I (\tilde{Q}_1 \,
_{^{j}}\widetilde{\infty}_{^0} \tilde{Q}_2 )
\end{equation}
assuming that the sewing operations 
appeared in (\ref{double-ext-equ}) are all well-defined. 
\end{lemma}
\pf
Let $\lambda_i \in \C,  i=1,2$ be such that 
$\tilde{Q}_i = \lambda_i \psi_{m_i} \otimes 
\bar{\psi}_{m_i}  (Q_i), i=1,2$. 
Let $(a_0^{(i)}, A^{(i)})$ be the local coordinate map at $j$-th 
puncture of $Q_1$ and Let $B^{(0)}$ be the local
coordinate map at $\infty$ in $Q_2$. 
By (\ref{I-sew-ext}),  the $\delta$ image of 
the left hand side of (\ref{double-ext-equ}) equals to
\bea
& & ((( \delta(\tilde{P}) _{^i}\widetilde{\infty}_{^0}
\lambda_1 \psi_{m_1}(Q_1) ) _{^{l+m_1-1+i}}\widetilde{\infty}_{^0}
\psi_{m_1}(\bar{Q}_1) ) \nn
& & \hspace{4cm}  _{^{i+j-1}}\widetilde{\infty}_{^0} 
\lambda_2  \psi_{m_2}(Q_2) )
_{^{l+m_1+m_2+i+j-3}}\widetilde{\infty}_{^0} \psi_{m_2}(\bar{Q}_2)  \nonumber 
\eea
By the associativity of the partial operad $\tilde{K}^c$, 
the above formula equals to 
\bea
& & ( \delta(\tilde{P}) _{^i}\widetilde{\infty}_{^0}
(\lambda_1 \psi_{m_1}(Q_1) _{^j}\widetilde{\infty}_{^0} 
\lambda_2 \psi_{m_2}(Q_2)) )  \nn
& & \hspace{3cm}  _{^{l+m_1+m_2-2+i}}\widetilde{\infty}_{^0}
(\psi_{m_1}(\bar{Q}_1) _{^j}\widetilde{\infty}_{^0} \psi_{m_2}(\bar{Q}_2))   \nn
& & \hspace{1cm} = 
( \delta(\tilde{P}) _{^i}\widetilde{\infty}_{^0}
(\lambda_1 \lambda_2 e^{\Gamma(A^{(i)}, B^{(0)}, a_0^{(i)})c}
\psi_{m_1+m_2-1} (Q_1\, _{^j}\infty_{^0} Q_2) ) )   \nn
& & \hspace{3cm} _{^{l+m_1+m_2-2+i}}\widetilde{\infty}_{^0}
e^{\Gamma(\overline{A^{(i)}}, \overline{B^{(0)}}, \overline{a_0^{(i)}})c}  
\psi_{m_1+m_2-1} (\overline{Q_1\, _{^j}\infty_{^0} Q_2})  \nn
& & \hspace{1cm} = \delta \big( 
\tilde{P} _{^i}\widetilde{\infty}_{^0}^I (
\lambda_1 \lambda_2 e^{\Gamma(A^{(i)}, B^{(0)}, a_0^{(i)})c}
e^{\Gamma(\overline{A^{(i)}}, \overline{B^{(0)}}, \overline{a_0^{(i)}})c}  \nn
& & \hspace{4cm} 
\psi_{m_1+m_2-1} \otimes \bar{\psi}_{m_1+m_2-1} (Q_1\, _{^j}\infty_{^0} Q_2)
) \big)  
\nn
& & \hspace{1cm} = \delta \big( \tilde{P} 
_{^i}\widetilde{\infty}_{^0}^I  
(\lambda_1 \psi_{m_1}\otimes \bar{\psi}_{m_1}(Q_1)  
 _{^j}\widetilde{\infty}_{^0} \lambda_2 \psi_{m_2}\otimes
\bar{\psi}_{m_2}(Q_2) ) \big) \nn
& & \hspace{1cm} = \delta(\tilde{P} _{^i}\widetilde{\infty}_{^0}^I
(\tilde{Q}_1\, _{^i}\widetilde{\infty}_{^0} \tilde{Q}_2) ), \nonumber
\eea
which is nothing but the right hand side of (\ref{double-ext-equ}). 
\epf

Boundary sewing operations and interior sewing operations 
induce the following 
partially defined substitution maps $\tilde{\gamma}$: 
\bea
&&
\tilde{\Upsilon}^c(n; l) \times 
\tilde{\Upsilon}^c(m_1; k_1) \times \dots 
\times \tilde{\Upsilon}^c(m_n; k_n)
\times \tilde{K}^c(n_1) \times \dots \times \tilde{K}^c(n_l)   \nn
&&\hspace{2cm} \xrightarrow{\tilde{\gamma}} 
\tilde{\Upsilon}^c(m_1+\dots +m_n; k_1+\dots +k_n+n_1+\dots +n_l)
\nonumber
\eea

The following Proposition is clear. 
\begin{prop}
$(\tilde{\mathfrak{S}}^c|
\tilde{K}^c \otimes \overline{\tilde{K}^{\bar{c}}}, \tilde{\gamma})$ 
is a $\R_+|\C^{\times}$-rescalable $2$-colored partial operad. 
\end{prop}

We will call $(\tilde{\mathfrak{S}}^c|
\tilde{K}^c \otimes \overline{\tilde{K}^{\bar{c}}}, \tilde{\gamma})$
{\it Swiss-cheese partial 
operad with central charge $c$}. Note that 
$\tilde{\mathfrak{S}}^c$ restricted on $\Upsilon$ is just
$\tilde{\Upsilon}^c$ which was introduced in \cite{HKo1}.


\subsection{Smooth $\tilde{\mathfrak{S}}^{c}|\tilde{K}^c\otimes 
\overline{\tilde{K}^{\bar{c}}}$-algebras}

Let $V^{O}=\oplus_{n\in \R}V^{O}_{(n)}$, where $V^{O}_{(n)}$ has a
structure of irreducible $\R_+$-module given by 
$r \mapsto r^n \id_{V^{O}_{(n)}}$ for $r\in \R_+$.
Let $V^{C}=\oplus_{(m,n)\in \R\times \R}V^{C}_{(m,n)}$, where
$V^{C}_{(m,n)}$ has a structure of irreducible $\C^{\times}$-module
given by $z \mapsto z^m \bar{z}^n\, \id_{V^{C}_{(m,n)}}$ for 
$z\in \C^{\times}$ (recall (\ref{branch-cut}) and (\ref{power-fun}))
for all $m,n$. This is also implies that $V^{C}_{(m,n)}=0$ for 
all $m-n\notin \Z$. Let 
$$
(V^O|V^C,\nu_{\tilde{\mathfrak{S}}^{c}|\tilde{K}^c\otimes 
\overline{\tilde{K}^{\bar{c}}}}, 
\nu_{\tilde{K}^c\otimes \overline{\tilde{K}^{\bar{c}}}})
$$ 
be a $\R_+|\C^{\times}$-rescalable
$\tilde{\mathfrak{S}}^{c}|\tilde{K}^c\otimes 
\overline{\tilde{K}^{\bar{c}}}$-algebra.

\begin{defn}\label{sw-ch-alg}{\rm  
The $\R_+|\C^{\times}$-rescalable
$\tilde{\mathfrak{S}}^{c}|\tilde{K}^c\otimes 
\overline{\tilde{K}^{\bar{c}}}$-algebra
$(V^O|V^C,\nu_{\tilde{\mathfrak{S}}^{c}|\tilde{K}^c\otimes 
\overline{\tilde{K}^{\bar{c}}}}, 
\nu_{\tilde{K}^c\otimes \overline{\tilde{K}^{\bar{c}}}})$
is called {\it smooth} if it satisfies the following two conditions:
\bnu
\item $\dim V^{O}_{(s)}<\infty$ for $s\in \R$, $V_{(n)}^O=0$ for $n<<0$ 
and $\dim V^{C}_{(m,n)}<\infty$ for $m,n\in \Z$,
$V_{(m,n)}^{C}$ for $m<<0$ or $n<<0$. 

\item For $w\in (V^O)'$, $v_1, \dots, v_m\in V^O$,
$u_1, \dots, u_n\in V^C$, and 
$\tilde{P}\in \tilde{\Upsilon}^c(m;n)$, the following map
$$
\tilde{P} \mapsto 
\langle w, \nu_{\tilde{\mathfrak{S}}^{c}|\tilde{K}^c\otimes 
\overline{\tilde{K}^{\bar{c}}}}(\tilde{P})(v_1\otimes \dots \otimes v_m 
\otimes u_1 \otimes \dots \otimes u_n)\rangle 
$$
is linear on fiber and smooth on the base space $\Upsilon(m;n)$. 
\enu 
}
\end{defn}

By \cite{H4}, a vertex operator algebra $U$ 
canonically gives a $\tilde{K}^c$-algebra. 
Using the doubling map $\delta$, this $\tilde{K}^c$-algebra
naturally gives a smooth $\tilde{\mathfrak{S}}^c|\tilde{K}^c\otimes 
\overline{\tilde{K}^{\bar{c}}}$-algebra, in which
$V^O=U$ and $V^C=U\otimes U$.
It is also easy to see that this 
$\tilde{\mathfrak{S}}^c|\tilde{K}^c\otimes 
\overline{\tilde{K}^{\bar{c}}}$-algebra
canonically gives an analytic 
open-closed field algebra which is nothing but
$(U\otimes U, U, \id_U, \id_U)$ or
simply $(U\otimes U, U)$ discussed in Section 1. 
We still denote this smooth 
$\tilde{\mathfrak{S}}^c|\tilde{K}^c\otimes 
\overline{\tilde{K}^{\bar{c}}}$-algebra 
as $(U\otimes U, U)$.
A smooth $\tilde{\mathfrak{S}}^c|\tilde{K}^c\otimes 
\overline{\tilde{K}^{\bar{c}}}$-algebra 
containing $(U\otimes U, U)$ as a subalgebra is called a {\it 
smooth $\tilde{\mathfrak{S}}^c|\tilde{K}^c\otimes 
\overline{\tilde{K}^{\bar{c}}}$-algebra over $U$}.

Let $(V_{cl}, V_{op}, m_{cl-op})$ 
be an open-closed field algebra over $V$. 
By definition, $V_{cl}$ is a conformal full field algebra over
$V\otimes V$. By the results in \cite{Ko1}, $V_{cl}$
has a structure of smooth 
$\tilde{K}^c\otimes \overline{\tilde{K}^{\bar{c}}}$-algebra structure,
we denote it as 
$(V_{cl}, \nu_{\tilde{K}^c\otimes \overline{\tilde{K}^{\bar{c}}}})$.

Let $Q\in \Upsilon(n;l)$ of form (\ref{gen-elm-moduli}) 
such that $r_1>\ldots >r_{n-1}> r_n=0$. 
We define a map 
$\nu_{\tilde{\mathfrak{S}}^c|\tilde{K}^c\otimes \overline{\tilde{K}^{\bar{c}}}}: \tilde{\mathfrak{S}}^c \rightarrow 
E_{V_{op}|V_{cl}}^{\R_+|\C^{\times}}$ as follows:
\bea  \label{Phi-def}
&&\nu_{\tilde{\mathfrak{S}}^c|\tilde{K}^c\otimes \overline{\tilde{K}^{\bar{c}}}}(\lambda \psi_{n}^{\mathfrak{S}}(Q))
(v_1\otimes \dots  \otimes v_n \otimes u_1\otimes \dots \otimes u_l ) \nn
&&\hspace{0.3cm} := \lambda e^{-L_-(B^{(0)})} 
m_{cl-op}^{(l;n)}(e^{-L_+(A^{(1)})}(a_0^{(1)})^{-L(0)}\otimes
e^{-L_+(\overline{A^{(1)}})}(\overline{a_0^{(1)}})^{-L(0)}u_1,    \nn
&&\hspace{2cm} \dots 
e^{-L_+(A^{(l)})}(a_0^{(l)})^{-L(0)} \otimes 
e^{-L_+(\overline{A^{(l)}})}(\overline{a_0^{(l)}})^{-L(0)} u_l;    \nn
&&\hspace{2cm} e^{-L_+(B^{(1)})}(b_0^{(1)})^{-L(0)}v_1, \dots,  
e^{-L_+(B^{(n)})}(b_0^{(n)})^{-L(0)}v_n; \nn
&&\hspace{5cm} z_1, \bar{z}_1, \dots, z_l,\bar{z}_l ; r_1,\dots, r_{n}),
\eea 
where $L_{\pm}(A) = \sum_{j=1}^{\infty} L(\pm j)A_j$ for any 
$A=\{A_1, A_2, \dots\}, A_j\in \C$, 
for $u_1,\ldots, u_l\in V_{cl}, v_1, \ldots, v_n\in V_{op}$. 
Let $e_l$ be the identity element of $S_l$. 
$\forall \sigma \in S_n$, we define 
\beq  \label{perm-S-alg-def}
\nu_{\tilde{\mathfrak{S}}^c|\tilde{K}^c\otimes \overline{\tilde{K}^{\bar{c}}}} \big( (\sigma, e_l)(\lambda \psi_{n}^{\mathfrak{S}}(Q)) \big) =
 \nu_{\tilde{\mathfrak{S}}^c|\tilde{K}^c\otimes \overline{\tilde{K}^{\bar{c}}}}  \big( \lambda \psi_{n}^{\mathfrak{S}}(Q) \big) .
\eeq
We have finished the definition of $\nu_{\tilde{\mathfrak{S}}^c|\tilde{K}^c\otimes \overline{\tilde{K}^{\bar{c}}}}$ in all cases. 
By results in \cite{HKo1}, the restriction of 
$\nu_{\tilde{\mathfrak{S}}^c|\tilde{K}^c\otimes \overline{\tilde{K}^{\bar{c}}}}$ on $\tilde{\Upsilon}^c$ clearly gives a morphism of 
partial nonassociative operad from $\tilde{\Upsilon}^c$ 
to $E_{V_{op}}^{\R_+}$.


\begin{thm}
$(V_{op}|V_{cl}, \nu_{\tilde{\mathfrak{S}}^c|\tilde{K}^c\otimes \overline{\tilde{K}^{\bar{c}}}}, \nu_{\tilde{K}^c\otimes \overline{\tilde{K}^{\bar{c}}}})$ is 
a smooth $\tilde{\mathfrak{S}}^c|\tilde{K}^c\otimes \overline{\tilde{K}^{\bar{c}}}$-algebra.
\end{thm}
\pf
By the permutation property of open-closed field algebra
and (\ref{perm-S-alg-def}), it is clear that 
$\nu_{\tilde{\mathfrak{S}}^c|\tilde{K}^c\otimes \overline{\tilde{K}^{\bar{c}}}}$is equivariant with the actions of permutation groups. 

The conditions in Definition \ref{sw-ch-alg} 
are automatically satisfied. It remains to show that 
\beq   \label{sew-morph}
\nu_{\tilde{\mathfrak{S}}^c|\tilde{K}^c\otimes 
\overline{\tilde{K}^{\bar{c}}}}  \circ \tilde{\gamma} 
= \Gamma \circ ( \nu_{\tilde{\mathfrak{S}}^c|\tilde{K}^c\otimes 
\overline{\tilde{K}^{\bar{c}}}},
\dots, \nu_{\tilde{\mathfrak{S}}^c|\tilde{K}^c\otimes 
\overline{\tilde{K}^{\bar{c}}}}, 
\nu_{\tilde{K}^c\otimes \overline{\tilde{K}^{\bar{c}}}}, \dots,
\nu_{\tilde{K}^c\otimes \overline{\tilde{K}^{\bar{c}}}})
\eeq
as a map
\bea
&&\tilde{\Upsilon}^c(n; l) \times 
\tilde{\Upsilon}^c(m_1; k_1) \times \dots 
\times \tilde{\Upsilon}^c(m_n; k_n)
\times \tilde{K}^c(n_1) \times \dots \times \tilde{K}^c(n_l) \nn
&& \hspace{2cm}  \rightarrow 
\hom(V_{cl}^{\otimes k_1+\dots +k_n + n_1+\dots +n_l}\otimes 
V_{op}^{\otimes n+m_1+\dots +m_n} , \overline{V_{op}}).
\eea
Thanks to the doubling map $\delta$ and 
(\ref{double-ext-equ}), $\tilde{\mathfrak{S}}^c$ can be viewed as
a partial suboperad of $\tilde{K}^c$ with single-sewing operations
for the punctures in $\R_+$ and double-sewing operations
for mirror pairs of punctures in upper and lower half planes.
By using the $V$-invariant boundary condition,  
the chirality splitting property and the 
convergence and extension properties
of any products and iterates of intertwining operators
proved by Huang for any $V$ satisfying the condition in 
Theorem \ref{ioa}, 
it is easy to generalize the proof of 
Huang's fundamental result Proposition 5.4.1 in \cite{H4} 
and results in \cite{H5}\cite{H6} to 
shows that (\ref{sew-morph}) holds.  The arguements are 
standard but tedious. We omit the details.
\epf

\renewcommand{\theequation}{\thesection.\arabic{equation}}
\renewcommand{\thethm}{\thesection.\arabic{thm}}
\setcounter{equation}{0}
\setcounter{thm}{0}


\section{Categorical formulation}

In this section, we study 
open-closed field algebras over $V$ from a 
tensor-categorical point of view.  
In Section 3.1, we recall some basic ingredients of the 
vertex tensor categories. In Section 3.2, we reformulate 
the notion of open-closed field algebra over $V$ categorically
by a categorical notion called 
open-closed $\mathcal{C}_V|\mathcal{C}_{V\otimes V}$-algebra.

\subsection{Vertex tensor categories}

The theory of tensor products for modules over a
vertex operator algebra was developed by Huang and Lepowsky
\cite{HL2}-\cite{HL5}\cite{H3}. 
By Theorem \ref{ioa} and our assumption on $V$, 
the category of $V$-modules, denoted as $\mathcal{C}_V$, 
have a structure of vertex tensor category \cite{HL2}, 
In particular, it has a structure of semisimple 
braided tensor category.

We review some of the ingredients of vertex tensor 
category $\mathcal{C}_V$ and set our notations along the way.

There is a tensor product bifunctor 
$\boxtimes_{P(z)}: \mathcal{C}_V \times \mathcal{C}_V \rightarrow 
\mathcal{C}_V$ 
for each $P(z), z\in \C^{\times}$ in sphere partial operad $K$, 
where $P(z)$ is the conformal equivalence class of sphere with 
three punctures $0, z, \infty$ and standard local coordinates 
\cite{H4}. We denote $\boxtimes_{P(1)}$ simply as $\boxtimes$. 
For any pair of $V$-modules $W_1,W_2$, 
the module $W_1\boxtimes_{P(z)} W_2$ is spanned by 
the homogeneous components of 
$w_1\boxtimes_{P(z_1)} w_2 \in \overline{W_1\otimes W_2}, 
\forall w_1\in W_1, w_2\in W_2$.

For each $V$-module $W$, there is a left unit isomorphism 
$l_W: V\boxtimes W \rightarrow W$ defined by
\beq  \label{l-unit-cat}
\overline{l_W}(v\boxtimes w) = 
Y_{W}(v,1)w, \quad \quad \forall v\in V, w\in W,
\eeq
where $\overline{l_W}$ is the unique extension
of $l_W$ on $\overline{V\boxtimes W}$ and
$Y_W$ is the vertex operator which defines the module
structure on $W$,  and a right unit isomorphism 
$r_W: W\boxtimes V \rightarrow W$ defined by
\beq \label{r-unit-cat}
\overline{r_W}(w\boxtimes v) = e^{L(-1)}Y_{W}(v, -1)w, 
\quad \quad \forall v\in V, w\in W.
\eeq

\begin{rema} {\rm 
We have used ``overline'' for the extensions of maps, algebraic
completions of graded vector spaces and 
complex conjugations of complex variables.  
One shall not confuse them because they acts on different things. 
}
\end{rema}

Let $W_1$ and $W_2$ be $V$-modules. For a given path 
$\gamma \in \C^{\times}$ from a point $z_1$ to $z_2$, there is
a parallel isomorphism associated to this path 
$$
\mathcal{T}_{\gamma} : W_1 \boxtimes_{P(z_1)} W_2 \longrightarrow 
W_1 \boxtimes_{P(z_2)} W_2.  
$$
Let $\Y$ be the intertwining operator corresponding 
to the intertwining 
map $\boxtimes_{P(z_2)}$ and $l(z_1)$ the value of the logarithm of 
$z_1$ determined by $\log z_2$ and analytic continuation 
along the path $\gamma$. 
For $w_1\in W_1,w_2\in W_2$, $\mathcal{T}_{\gamma}$ is defined by
$$
\overline{\mathcal{T}_{\gamma}} (w_1\boxtimes_{P(z_1)} w_2)
=\Y(w_1, e^{l(z_1)})w_2,
$$
where $\overline{\mathcal{T}_{\gamma}}$ is the natural extension of 
$\mathcal{T}_{\gamma}$.
Moreover, the parallel isomorphism
depends only on the homotopy class of $\gamma$.

For $z_1>z_2>z_1-z_2>0$ and each triple of 
$V$-modules $W_1, W_2, W_3$,  
there is an associativity isomorphism: 
$$
\mathcal{A}^{P(z_{1}-z_{2}), P(z_{2})}_{P(z_{1}), P(z_{2})}: 
W_1\boxtimes_{P(z_{1})} (W_2\boxtimes_{P(z_{2})} W_3)\to 
(W_1\boxtimes_{P(z_{1}-z_{2})} W_2)\boxtimes_{P(z_{2})} W_3,
$$
which is characterized by 
\begin{equation}\label{assoc-iso}
\overline{\mathcal{A}^{P(z_{1}-z_{2}), P(z_{2})}_{P(z_{1}), P(z_{2})}}
(w_{(1)}\boxtimes_{P(z_{1})} (w_{(2)}\boxtimes_{P(z_{2})} w_{(3)}))
=(w_{(1)}\boxtimes_{P(z_{1}-z_{2})} w_{(2)})\boxtimes_{P(z_{2})}w_{(3)}
\end{equation}
for $w_{(i)} \in W_i, i=1,2,3$. The associativity isomorphism 
$\mathcal{A}$ of the braided tensor category 
is characterized by the following commutative diagram: 
\beq  \label{A-T-A-T}
\xymatrix{
W_1\boxtimes_{P(z_{1})} (W_2\boxtimes_{P(z_{2})} W_3)
\ar[rrr]^{\hspace{0.6cm} \mathcal{T}_{\gamma_1} \circ (\id_{W_1} \boxtimes_{P(z_1)} \mathcal{T}_{\gamma_2}) }
\ar[d]_{\mathcal{A}^{P(z_{1}-z_{2}), P(z_{2})}_{P(z_{1}), P(z_{2})}} & & &
W_1\boxtimes (W_2\boxtimes W_3)
\ar[d]^{\mathcal{A} }     \\
(W_1\boxtimes_{P(z_{1}-z_{2})} W_2)\boxtimes_{P(z_{2})} W_3
\ar[rrr]^{\hspace{0.7cm}\mathcal{T}_{\gamma_{2}}\circ (\mathcal{T}_{\gamma_{3}}
\boxtimes_{P(z_{2})} \id_{W_3})} & & & 
(W_1 \boxtimes W_2) \boxtimes W_3 \, ,
}   
\eeq
where $\gamma_1, \gamma_2, \gamma_3$ are paths in $\R_+$
from $z_1, z_2, z_1-z_2$ to $1$, respectively.

There is also a braiding isomorphism, for $z>0$,  
$\mathcal{R}_+^{P(z)}: W_1\boxtimes_{P(z)} W_2 
\rightarrow W_2\boxtimes_{P(z)} W_1$ 
for each pair of $V$-modules $W_1,W_2$, defined as 
\begin{equation}  \label{R-+-cat}
\overline{\mathcal{R}_+^{P(z)}}(w_1\boxtimes_{P(z)} w_2) = e^{zL(-1)} 
\overline{\mathcal{T}}_{\gamma_+} (w_2\boxtimes_{P(-z)} w_1),
\end{equation}
where $\gamma_+$ 
is a path from $-z$ to $z$ given by $\gamma_+(t) = -e^{i\pi t}z, t\in [0,1]$. 
The inverse of $\mathcal{R}_+^{P(z)}$ is denoted by 
$\mathcal{R}_-^{P(z)}$, which is characterized by 
\begin{equation}   \label{R---cat}
\overline{\mathcal{R}_-^{P(z)}}(w_2\boxtimes_{P(z)} w_1) = e^{zL(-1)} 
\overline{\mathcal{T}}_{\gamma_-} (w_1\boxtimes_{P(-z)} w_2),
\end{equation}
where $\gamma_-$ is a path given by $\gamma_-(t)= -e^{-i\pi t}z$. 
We denote $\mathcal{R}_{\pm}^{P(1)}$ simply as $\mathcal{R}_{\pm}$. 
Accordingly, we define the twist
by $\theta_W = e^{-2\pi iL(0)}$ for any $V$-module $W$.

For $z_1,z_2\in \R_+$, 
the naturalness of $\mathcal{T}$ implies the commutativity 
of the following diagram: 
\beq  \label{R-T-R-T}
\xymatrix{
W_1\boxtimes_{P(z_{1})} W_2
\ar[rrr]^{\hspace{0.6cm} \mathcal{T}_{\gamma}}
\ar[d]_{\mathcal{R}_{\pm}^{P(z_1)}} & & &
W_1\boxtimes_{P(z_2)} W_2
\ar[d]^{\mathcal{R}_{\pm}^{P(z_2)} }     \\
W_2\boxtimes_{P(z_{1})} W_1
\ar[rrr]^{\hspace{0.5cm} \mathcal{T}_{\gamma}} & & & 
W_2 \boxtimes_{P(z_2)} W_1 \, ,
}   
\eeq
where $\gamma$ is a path in $\R_+$ from $z_1$ to $z_2$.

Let $\mathcal{V}_{W_1W_2}^{W_3}$ denotes the 
space of intertwining operators of type $\binom{W_3}{W_1W_2}$. 
For the chosen branch cut (\ref{branch-cut}), 
the space $\mathcal{V}_{W_1W_2}^{W_3}$ is canonically isomorphic to 
the space of intertwining maps denoted as $\mathcal{M}[P(z)]_{W_1W_2}^{W_3}$. 
By the universal property of tensor product, the space
$\mathcal{M}[P(z)]_{W_1W_2}^{W_3}$ is canonically identified with
the space $\hom_V(W_1\boxtimes_{P(z)} W_2, W_3)$. 
Let $\Y\in \mathcal{V}_{W_1W_2}^{W_3}$. We denote the 
corresponding morphism in $\hom_V(W_1\boxtimes W_2, W_3)$ as
$m_{\Y}^{P(z)}$. 
For $r\in \Z$, $\Omega_r: \mathcal{V}_{W_1W_2}^{W_3}\rightarrow
\mathcal{V}_{W_2W_1}^{W_3}$ is an isomorphism defined as follows: 
$$
\Omega_r(\Y)( w_2, x) w_1 = e^{xL(-1)}\Y(w_1, e^{(2r+1)\pi i}x) w_2
$$
for all $w_1\in W_1, w_2\in W_2$. Then by Proposition 3.1 in \cite{Ko1}, 
where the definition of braiding is opposite to our choice here (recall 
(\ref{R-+-cat})(\ref{R---cat})), we obtain the following identities:
\beq \label{Omega-R}
m_{\Y}^{P(z)} = m_{\Omega_0(\Y)}^{P(z)} \circ \mathcal{R}_+^{P(z)} =
m_{\Omega_{-1}(\Y)}^{P(z)} \circ \mathcal{R}_-^{P(z)} .
\eeq

The tensor product $V\otimes V$ is also a 
vertex operator algebra \cite{FHL}. 
Moreover, it was shown in \cite{HKo2} that 
$V\otimes V$ also satisfies the conditions in
Theorem \ref{ioa}. Therefore, 
$\mathcal{C}_{V\otimes V}$, the category of 
$V\otimes V$-modules, also has a structure of 
semisimple braided tensor category. 
For $z,\zeta\in \C^{\times}$, 
let $\boxtimes_{P(z)P(\zeta)}$ be the tensor
product bifunctor in $\mathcal{C}_{V\otimes V}$ defined by
$$
(A\otimes B) \boxtimes_{P(z)P(\zeta)} (C\otimes D)
= (A\boxtimes_{P(z)}B)\otimes (C\boxtimes_{P(\zeta)}D),
$$ 
where $A, B, C, D$ are $V$-modules. 
We denote the bifunctor $\boxtimes_{P(1)P(1)}$ simply as $\boxtimes$.

There are a few different braiding structures on 
$\mathcal{C}_{V\otimes V}$ \cite{Ko1}. We choose the one given by 
$\mathcal{R}_+\otimes \mathcal{R}_-$. The twist $\theta_A: A\rightarrow A$,
for $A\in \mathcal{C}_{V\otimes V}$,  is given as follows
\beq  \label{twist-theta-2}
\theta_A = e^{-2\pi iL^L(0)} \otimes e^{2\pi iL^R(0)}.
\eeq
An object $A$ is said to have a trivial twist if $\theta_A =\id_A$.

Now we recall the definition of (commutative) associative
algebra in a braided tensor category $\mathcal{C}$
with tensor product $\otimes$, unit object $1_{\mathcal{C}}$,
left unit isomorphism $l_W$ and right unit isomorphism $r_W$
for any object $W$, the associativity $\mathcal{A}$ and 
the braiding $\mathcal{R}$. 
\begin{defn}   {\rm
An associative algebra in $\mathcal{C}$
(or associative $\mathcal{C}$-algebra) is an object $A$ in 
$\mathcal{C}$ together with $\mu_A: A\otimes A\rightarrow A$
and a monomorphism $\iota_A: 1_{\mathcal{C}} \rightarrow A$ satisfying
\bnu
\item {\it Associativity}: $\mu_A \circ (\mu_A \otimes \id_A)\circ 
\mathcal{A} = \mu_A \circ (\id_A \otimes \mu_A)$. 
\item {\it Unit properties}: $\mu_A \circ (\iota_A \otimes \id_A)
\circ l_A^{-1} = \mu_A \circ (\iota_A \otimes \id_A) \circ r_A^{-1}
=\id_A$. 
\enu
$A$ is called {\it commutative} if $\mu_A =\mu_A \circ \mathcal{R}$. 
}
\end{defn}

The following Theorem is proved in \cite{HKo1}.
\begin{thm}  \label{osva-thm-cat}
The category of open-string vertex operator algebras
over $V$ is isomorphic to the category of associative
$\mathcal{C}_V$-algebras. 
\end{thm}

The following Theorem is proved in \cite{Ko1}. 
\begin{thm}  \label{ffa-thm-cat}
The category of conformal full field algebras over $V\otimes V$
is isomorphic to the category of commutative associative 
algebras in $\mathcal{C}_{V\otimes V}$ with a trivial twist. 
\end{thm}

We are interested in studying the relation 
between above two algebras as ingredients of 
an open-closed field algebra over $V$. Notice that 
these two algebras live in different categories. 
So we will first 
discuss a functor between these two categories.

Recall a functor $T_{P(z)}: \mathcal{C}_{V\otimes V} \rightarrow
\mathcal{C}_V$ \cite{HL4}. 
In particular, for $W_1, W_2$ 
being $V$-module, $T_{P(z)}(W_1\otimes W_2) = W_1\boxtimes_{P(z)} W_2$. 
We will simply write $T_{P(1)}$ as $T$. Let
\beq \label{varphi-0}
\varphi_0:=l_{V}^{-1}=r_{V}^{-1}: V \rightarrow 
V \boxtimes V = T(V \otimes V).
\eeq
For each four $V$-modules $W_i^L, W_i^R, i=1,2$, notice that
\bea \label{varphi-2-domain}
T(W_1^L\otimes W_1^R) \boxtimes T(W_2^L\otimes
W_2^R) &=& (W_1^L\boxtimes W_1^R) \boxtimes (W_2^L\boxtimes W_2^R), 
\\
T((W_1^L\otimes W_2^R) \boxtimes (W_2^L\otimes
W_2^R)) &=& (W_1^L\boxtimes W_2^L) \boxtimes (W_1^R\boxtimes W_2^R). 
\label{varphi-2-codomain}
\eea
We define $\varphi_2: (W_1^L\boxtimes W_1^R) 
\boxtimes (W_2^L\boxtimes W_2^R) \rightarrow 
(W_1^L\boxtimes W_2^L) \boxtimes (W_1^R\boxtimes W_2^R)$ by 
\begin{equation} \label{varphi-2-defn}
\varphi_2:= 
\mathcal{A} \circ
(\id_{W_1} \boxtimes \mathcal{A}^{-1}) \circ 
(\id_{W_1} \boxtimes \mathcal{R}_- \boxtimes \id_{W_4}) \circ
(\id_{W_1} \boxtimes \mathcal{A}) \circ \mathcal{A}^{-1}.
\end{equation}
The above definition
of $\varphi_2$ can be naturally extended to a morphism 
$T(A)\boxtimes T(B) \rightarrow T(A\boxtimes B)$ for 
each pair of objects $A$ and $B$ in $\mathcal{C}_{V\otimes V}$. 
We still denote the extended morphism as $\varphi_2$.  
The following result is clear. 
\begin{lemma}   \label{lemma-T-ten-fun}
$T$ together with $\varphi_0, \varphi_2$ given in (\ref{varphi-0}) 
(\ref{varphi-2-defn}) is a monoidal functor. 
\end{lemma}

For any four objects  $W_i,i=1,2,3,4$ in $\mathcal{C}_V$, 
we define a morphism $\sigma$ in 
$\hom ( (W_1\boxtimes W_2)\boxtimes (W_3\boxtimes W_4), 
(W_3\boxtimes W_4)\boxtimes (W_1 \boxtimes W_2) )$ 
according to the following graph:
\beq  \label{sigma-fig}
\begin{picture}(14,2.5)
\put(5, 0.4){\resizebox{4cm}{1.8cm}{\includegraphics{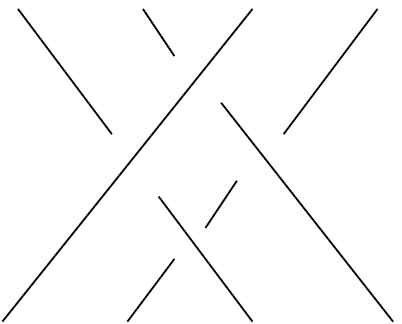}}}
\put(4.9, 0){$W_1$}\put(6.2,0){$W_2$}\put(7.4,0){$W_3$}
\put(8.9,0){$W_4$}
\put(4.9,2.3){$W_3$}\put(6.4,2.3){$W_4$}
\put(7.4,2.3){$W_1$}\put(8.9,2.3){$W_2$}
\end{picture}
\eeq
Clearly, for any $A\in \mathcal{C}_{V\otimes V}$, 
$\sigma$ can be extended to an automorphism on 
$T(A)\boxtimes T(A)$, denoted as $\sigma_A$.

\begin{prop} \label{prop-alg-alg-comm}
Let $(A, \mu_{A}, \iota_{A})$ 
be a commutative algebra in $\mathcal{C}_{V\otimes V}$. Let 
$\mu_{T(A)}:= T(\mu_A)\circ \varphi_2$ and 
$\iota_{T(A)}:=T(\iota_A) \circ \varphi_0$.
Then 
$(T(A), \mu_{T(A)}, \iota_{T(A)})$ is an associative $\mathcal{C}_V$-algebra
satisfying the following commutativity: 
\beq \label{comm-T-cl}
\mu_{T(A)} = \mu_{T(A)} \circ \sigma_A.
\eeq
\end{prop}
\pf
That $(T(A), \mu_{T(A)}, \iota_{T(A)})$ is an associative 
$\mathcal{C}_V$-algebra follows from Lemma \ref{lemma-T-ten-fun}. 
We only prove (\ref{comm-T-cl}).
Let us assume $A=W_1\otimes W_2$ with $W_1,W_2$ being 
$V$-modules (the proof for general $A$ is the same). 
Consider the following diagram:
\beq  \label{T-comm-diag}
\xymatrix{ (W_1\boxtimes W_2)\boxtimes (W_1\boxtimes W_2) 
\ar[d]^{\sigma} \ar[r]^{\varphi_2} & (W_1\boxtimes W_1) 
\boxtimes (W_2\boxtimes W_2) 
\ar[d]^{\mathcal{R}_+ \boxtimes \mathcal{R}_-} 
\ar[r]^{\hspace{1cm} T(\mu_{cl})} & W_1\boxtimes W_2 \ar[d]^{\id}  \\
(W_1\boxtimes W_2)\boxtimes (W_1\boxtimes W_2)
\ar[r]^{\varphi_2} & (W_1\boxtimes W_1) \boxtimes (W_2\boxtimes W_2)
\ar[r]^{\hspace{1cm} T(\mu_{cl})} & W_1\boxtimes W_2 \,\, .
}
\eeq
The left subdiagram is commutative because two paths
corresponding to the same braiding; the right subdiagram is 
commutative because of the commutativity of $A$. Hence 
the above diagram is commutative. (\ref{comm-T-cl}) follows 
from the commutativity of (\ref{T-comm-diag}).
\epf

\subsection{Open-closed $\mathcal{C}_V|\mathcal{C}_{V\otimes V}$-algebras}

Let $((V_{cl}, m_{cl}, \iota_{cl}), 
(V_{op}, Y_{op}, \iota_{op}), \mathbb{Y}_{cl-op})$ be an open-closed
field algebra over $V$ throughout this subsection.

By Theorem \ref{ffa-thm-cat},  $(V_{cl}, m_{cl}, \iota_{cl})$, 
a conformal full field algebra over $V\otimes V$
is equivalent to a commutative associative algebra with a trivial twist in
$\mathcal{C}_{V\otimes V}$ with braiding $\mathcal{R}_{+-}$. 
We denoted this $\mathcal{C}_{V\otimes V}$-algebra as a triple
$(V_{cl}, \mu_{cl}, \iota_{cl})$,
where $\mu_{cl}=m_{\mathbb{Y}_f}$ (recall (\ref{formal-Y-ffa})). 
By Proposition \ref{prop-alg-alg-comm},
$(T(V_{cl}), \nu_{T(V_{cl})}, \iota_{T(V_{cl})})$ 
is an associative $\mathcal{C}_V$-algebra.  
The property of algebra $T(V_{cl})$ can be expressed in the 
following graphic equations: 
\beq  \label{T-cl-1-fig}
\epsfxsize  0.8\textwidth
\epsfysize  0.15\textwidth
\epsfbox{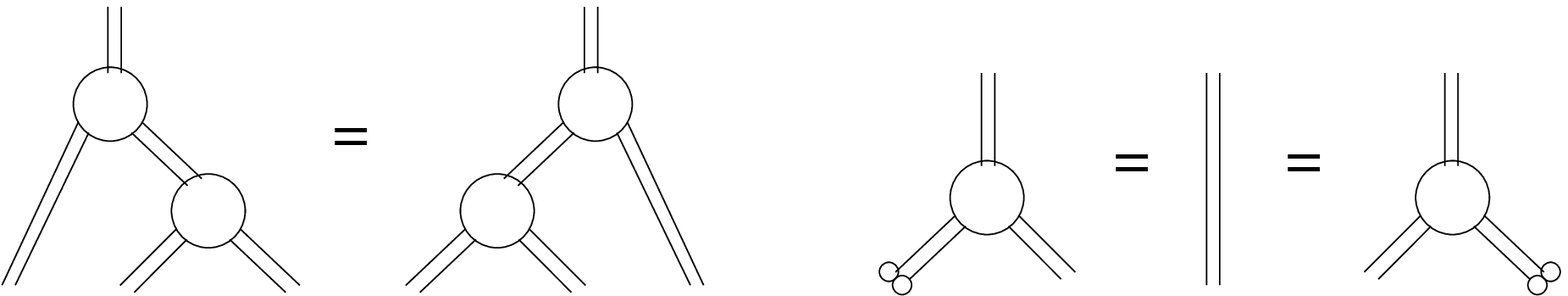}\quad ,
\eeq
\beq  \label{T-cl-2-fig}
\epsfxsize  0.4\textwidth
\epsfysize  0.15\textwidth
\epsfbox{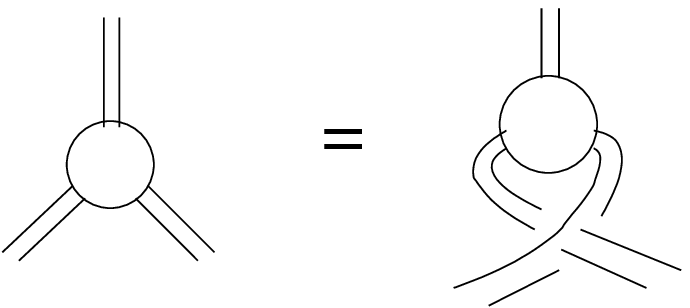}\quad\quad .
\eeq

By Theorem \ref{osva-thm-cat}, $(V_{op},Y_{op}, \iota_{op})$,
an open-string vertex operator algebra over $V$ is equivalent to 
an associative $\mathcal{C}_V$-algebra, denoted as a triple
$(V_{op}, \mu_{op}, \iota_{op})$,
where $\mu_{op} = m_{Y_{op}^f}$
(recall (\ref{formal-Y-osva-1}) and (\ref{formal-Y-osva-2}))., 
The defining properties of $(V_{op}, \mu_{op}, \iota_{op})$
can be expressed in the following graphic equations: 
\beq  \label{T-op-1-fig}
\epsfxsize  0.8\textwidth
\epsfysize  0.1\textwidth
\epsfbox{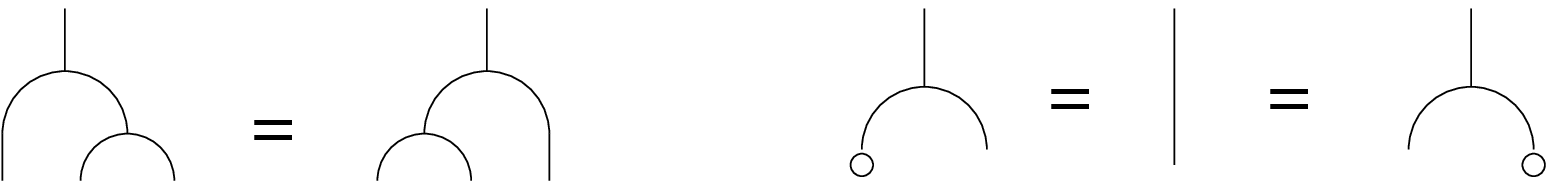}\quad .
\eeq

It remains to study the categorical formulation of the only 
remaining data $\mathbb{Y}_{cl-op}$ and its properties.

By the chirality splitting properties and the 
associativity of intertwining operator algebra \cite{H7}, there exist
intertwining operators $\Y^{(5)}, \Y^{(6)}$ such that 
\beq  \label{Y-hat-tilde}
\langle v', \mathbb{Y}_{cl-op}(u\otimes \bar{u}, z, \zeta)v\rangle =
\langle v', \Y^{(5)}(\Y^{(6)}(u^L, z-\zeta)u^R, \zeta)v \rangle
\eeq
for $v'\in V_{op}', v\in V_{op}, u^L\otimes u^R \in V_{cl}$ 
and $\zeta>z-\zeta>0$. 
For $u^L\otimes u^R \in V_{cl}$, we have
$u^L \boxtimes_{P(z-\zeta)} u^R \in
\overline{T_{P(z-\zeta)}(V_{cl})}$.

We define a map $\mu_{cl-op}^{P(z-\zeta)P(\zeta)}: 
T_{P(z-\zeta)}(V_{cl}) \boxtimes_{P(\zeta)} V_{op} \rightarrow V_{op}$ 
as follows:
\beq  \label{m-co-Y-co-0}
\mu_{cl-op}^{P(z-\zeta)P(\zeta)}: = m_{\Y^{(5)}}^{P(\zeta)} \circ 
(m_{\Y^{(6)}}^{P(z-\zeta)} \boxtimes_{P(\zeta)} \id_{V_{op}}).
\eeq
By the convergence and analytic properties of the right hand side of 
(\ref{Y-hat-tilde}) and 
the fact that module maps are continuous with respect to the topology on
graded vector spaces defined in Section 1, we obtain,
for $z>\zeta>z-\zeta>0$, $u^L\otimes u^R \in V_{cl}, v\in V_{op}$, 
the following identity:
\begin{equation} \label{m-co-Y-co}
\overline{\mu_{cl-op}^{P(z-\zeta)P(\zeta)}}
((u^L\boxtimes_{P(z-\zeta)} u^R) \boxtimes_{P(\zeta)} v)
= \Y^{(5)}(\Y^{(6)}(u^L, z-\zeta) u^R, \zeta)v,
\end{equation}
where $\overline{\mu_{cl-op}^{P(z-\zeta)P(\zeta)}}$ is the unique extension
of $\mu_{cl-op}^{P(z-\zeta)P(\zeta)}$ to the algebraic completion of 
$T_{P(z-\zeta)}(V_{cl}) \boxtimes_{P(\zeta)} V_{op}$. 
Moreover, by Proposition 1.2 in \cite{H9} 
(or the nondegeneracy of intertwining operator algebra \cite{nioa}),
such $V$-module map is unique.

Let $\gamma_1$ be a path in $\R_+$  
from $1$ to $z-\zeta$ and $\gamma_2$ a path in 
$\R_+$ from $1$ to $\zeta$. Then we define 
a map $\mu_{cl-op}: T(V_{cl})\boxtimes V_{op} \rightarrow V_{op}$ 
as follows:
\begin{equation} \label{m-co-def}
\mu_{cl-op} := \mu_{cl-op}^{P(z-\zeta)P(\zeta)} \circ \mathcal{T}_{\gamma_2}
\circ (\mathcal{T}_{\gamma_1} \boxtimes \id_{V_{op}}). 
\end{equation}
Since $\mathcal{T}_{\gamma}$ depends on path only homotopically, 
it is clear that above definition of $\mu_{cl-op}$ is independent of 
$z,\zeta$ in $\R_+$ and paths $\gamma_1,\gamma_2$ in $\R_+$.  
In particular, one can choose $\gamma_1$ and $\gamma_2$ to
be the straight line between $1$ and $z-\zeta, \zeta$
respectively.

By Theorem \ref{prop-op-cl-V}, to obtain a categorical formulation of 
open-closed field algebra over $V$ is enough to
study the categorical formulations of the unit
property (\ref{Y-co-unit}), 
Associativity I (\ref{asso-co-op}), 
Associativity II (\ref{asso-co-co}) 
and Commutativity I in Proposition \ref{comm-prop-1}.

We first consider the  property (\ref{Y-co-unit}). 

\begin{prop} \label{lemma-Y-co-unit}
The condition (\ref{Y-co-unit}) 
is equivalent to the following condition: 
\beq \label{id-cat-def}
\mu_{cl-op} \circ ( (T(\iota_{V_{cl}})\circ \varphi_0) 
\boxtimes \id_{V_{op}}) \circ l_{V_{op}}^{-1} =\id_{V_{op}},
\eeq
which can also be expressed by the following graphic equation:
\beq  \label{id-fig}
\epsfxsize  0.27\textwidth
\epsfysize  0.16\textwidth
\epsfbox{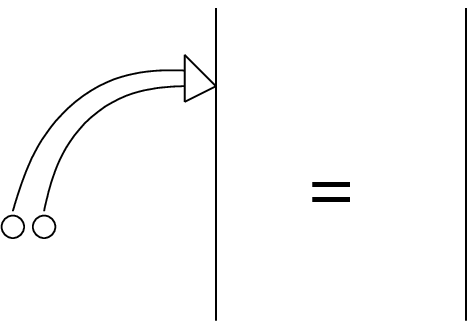}\quad \quad .
\eeq
\end{prop}
\pf
First, (\ref{Y-co-unit}) is equivalent to the following condition:
\beq  \label{Y-co-unit-1}
\mathbb{Y}_{cl-op}(\one_{cl}; z,\zeta)=\id_{V_{op}}, \quad \mbox{for
$z > \zeta >0$.}
\eeq
Recall that
$\one_{cl} = \iota_{V_{cl}}(\one\otimes \one)$. Replacing 
$u^L \otimes u^R$ by $\one_{cl}$ in (\ref{Y-hat-tilde}), 
one can see that both sides of 
equation (\ref{Y-hat-tilde}) are independent of $z$ and $\zeta$.
Hence (\ref{Y-co-unit-1}) holds for all $z, \zeta \in \R_+$. 
Using (\ref{Y-hat-tilde}) and (\ref{m-co-Y-co}), we obtain
\beq  \label{m-co-equ-1}
\overline{\mu_{cl-op}^{P(z-\zeta)P(\zeta)}}(T_{P(z-\zeta)}(\one_{cl}) 
\boxtimes_{P(\zeta)} v) = v
\eeq
for $v\in V_{op}$ and $z,\zeta\in \R_+$. 
Therefore,  we have, for $v\in V_{op}$,
\bea  \label{m-co-equ-2}
&& 
\mu_{cl-op}^{P(z-\zeta)P(\zeta)} \circ T_{P(z-\zeta)}(\iota_{V_{cl}})
\boxtimes_{P(\zeta)} \id_{V_{op}}  ((\one \boxtimes_{P(z-\zeta)} \one)
\boxtimes_{P(\zeta)} v)  \nn
&&\hspace{1cm} = \mu_{cl-op}^{P(z-\zeta)P(\zeta)} 
(T_{P(z-\zeta)}(\one_{cl}) \boxtimes_{P(\zeta)} v) \nn
&&\hspace{1cm} = v,
\eea
which can be further expressed equivalently as 
\bea \label{id-equ-1}
\id_{V_{op}} &=& 
\mu_{cl-op}^{P(z-\zeta)P(\zeta)} \circ 
\big( (T_{P(z-\zeta)}(\iota_{V_{cl}}) \circ \mathcal{T}_{\gamma_1}) 
\boxtimes_{P(\zeta)} \id_{V_{op}}\big)\circ \mathcal{T}_{\gamma_2}\circ
(l_{V}^{-1} \boxtimes \id_{V_{op}})  \circ l_{V_{op}}^{-1}  \nn
&=& \mu_{cl-op}^{P(z-\zeta)P(\zeta)} \circ 
(\mathcal{T}_{\gamma_1} \boxtimes_{P(\zeta)} \id_{V_{op}}) \circ 
\mathcal{T}_{\gamma_2} \circ \big( (T(\iota_{V_{cl}})\circ \varphi_0) 
\boxtimes \id_{V_{op}}\big) \circ l_{V_{op}}^{-1} \nn
&=& \mu_{cl-op}^{P(z-\zeta)P(\zeta)} \circ 
\mathcal{T}_{\gamma_2} \circ 
(\mathcal{T}_{\gamma_1} \boxtimes \id_{V_{op}}) \circ 
\big( (T(\iota_{V_{cl}})\circ \varphi_0) 
\boxtimes \id_{V_{op}}\big) \circ l_{V_{op}}^{-1} \nn
&=& \mu_{cl-op} \circ \big( (T(\iota_{V_{cl}})\circ \varphi_0) 
\boxtimes \id_{V_{op}}\big) \circ l_{V_{op}}^{-1}
\eea
where $\gamma_1$ and $\gamma_2$ are paths in $\R_+$ from 
$1$ to $z-\zeta$ and $\zeta$ respectively. 

Conversely, from (\ref{m-co-equ-1}), (\ref{m-co-equ-2}) and
(\ref{id-equ-1}), it is clear that 
(\ref{id-cat-def}) or (\ref{id-fig}) also
implies (\ref{Y-co-unit}). 
\epf

Now we consider the associativity II 
(recall Proposition \ref{prop-asso-co-co}). 
\begin{prop}  \label{cat-lemma-asso-co-co}
The associativity II is equivalent to the following condition: 
\begin{equation} \label{asso-co-co-cat}
\mu_{cl-op} \circ (\id_{T(V_{cl})} \boxtimes \mu_{cl-op})
= \mu_{cl-op} \circ \big( 
(T(\mu_{cl})\circ \varphi_2) \boxtimes \id_{V_{op}}\big) 
\circ \mathcal{A},
\end{equation}
which can also be expressed by the following graphic equation
\beq    \label{asso-co-co-fig}
\epsfxsize  0.5\textwidth
\epsfysize  0.2\textwidth
\epsfbox{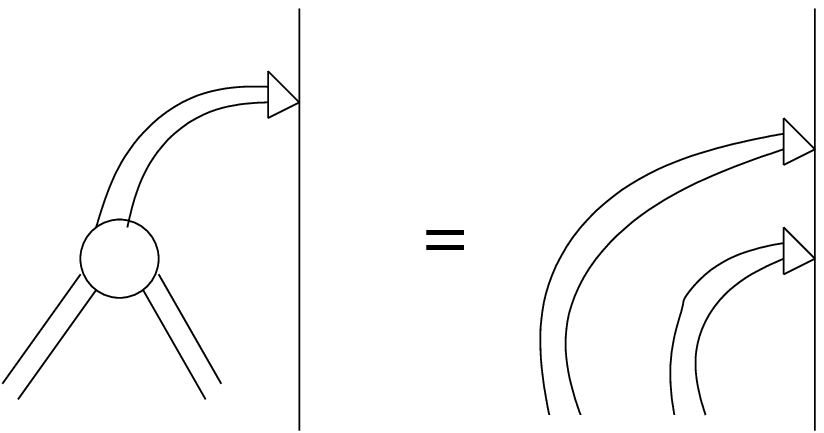}\,\, .
\eeq
\end{prop}
\pf
Using the convergence property of products and iterates 
of intertwining operators of $V$, it is not hard to show 
that, in the domain
\beq  \label{domain-LHS}
D_3: = \{ (z_1,\zeta_1, z_2, \zeta_2)| 
z_1>\zeta_1>z_2>\zeta_2>0, 
2\zeta_2>z_2, 2\zeta_1>z_1 + z_2 \},
\eeq 
we have 
\bea  \label{LHS-asso-co-co}
&&\hspace{-1cm}Y_{cl-op}(u_1^L\otimes u_1^R; z_1, \zeta_1) 
Y_{cl-op}(u_1^L\otimes u_1^R; z_1, \zeta_1)v   \nn
&&\hspace{0cm}=\Y^{(5)}(\Y^{(6)}(u_1^L, z_1-\zeta_1)u_1^R, \zeta_1)
\Y^{(5)}(\Y^{(6)}(u_2^L, z_2-\zeta_2)u_2^R, \zeta_2)v
\eea
for $u_1^L\otimes u_1^R, u_2^L\otimes u_2^R\in V_{cl}$
and $v\in V_{op}$. 

By (\ref{m-co-Y-co-0}(\ref{m-co-Y-co}) 
and the fact that module maps are continuous, we obtain
\bea  \label{asso-co-co-LHS-2}
&&\hspace{-1cm} Y_{cl-op}(u_1^L\otimes u_1^R; z_1, \zeta_1) 
Y_{cl-op}(u_1^L\otimes u_1^R; z_1, \zeta_1)v \nn
&&\hspace{0cm}=
\overline{\mu_{cl-op}^{P(z_1-\zeta_1)P(\zeta_1)} \circ
(\id_{T_{P(z_1-\zeta_1)}(V_{cl})} \boxtimes_{P(\zeta_1)} 
\mu_{cl-op}^{P(z_2-\zeta_2)P(\zeta_2)} } ) 
\big( \nn 
&&\hspace{1cm}(u_1^L\boxtimes_{P(z_1-\zeta_1)} u_1^R)
\boxtimes_{P(\zeta_1)} 
((u_2^L\boxtimes_{P(z_2-\zeta_2)} u_2^R) \boxtimes_{P(\zeta_2)} v) \big)
\eea
for $(z_1, \zeta_1, z_2, \zeta_2)$ in the domain $D_3$. 
Moreover, by the universal property of tensor product and
the nondegeneracy of intertwining operator algebra \cite{nioa},
\beq \label{cat-LHS-asso}
\mu_{cl-op}^{P(z_1-\zeta_1)P(\zeta_1)}\circ
(\id_{T_{P(z_1-\zeta_1)}(V_{cl})} \boxtimes_{P(\zeta_1)} 
\mu_{cl-op}^{P(z_2-\zeta_2)P(\zeta_2)} )
\eeq
is the unique morphism 
in $\hom (T_{P(z_1-\zeta_1)}(V_{cl}) \boxtimes_{P(\zeta_1)} 
(T_{P(z_2-\zeta_2)}(V_{cl}) \boxtimes_{P(\zeta_2)} V_{op}),  V_{op})$
such that the identity (\ref{asso-co-co-LHS-2}) holds in $D_3$.

On the other hand of (\ref{asso-co-co}), 
it is proved in \cite{HKo2} that $\mathbb{Y}$ can be expanded 
as follows: 
\beq 
\mathbb{Y}_{an}(u_1^L\otimes u_1^R; z, \zeta)u_2^L\otimes u_2^R
=\sum_{i=1}^N \Y_i^L(u_1^L, z)u_2^L \otimes \Y_i^R(u_1^R, \zeta)u_2^R 
\eeq
for some $N\in \Z_+$. There is a unique morphism  
$$
\mu_{cl}^{P(z_1-z_2)P(\zeta_1-\zeta_2)} \in 
\hom_{\mathcal{C}_{V\otimes V}} 
(V_{cl}\boxtimes_{P(z_1-z_2)P(\zeta_1-\zeta_2)}V_{cl}, V_{cl})
$$
such that, for $u, v\in V_{cl}$, 
\beq  \label{mu-cl-Y-z-zeta}
\overline{\mu_{cl}^{P(z_1-z_2)P(\zeta_1-\zeta_2)}}
(u\boxtimes_{P(z_1-z_2)P(\zeta_1-\zeta_2)}v)
= \mathbb{Y}_{an}(u; z_1-z_2, \zeta_1-\zeta_2)v.
\eeq
By the convergence property of intertwining operator algebra, 
it not hard to see that, in the domain 
\beq \label{domain-RHS}
D_4 :=\{ (z_1,\zeta_1,z_2,\zeta_2)|
z_1>z_2>\zeta_1>\zeta_2>0,  2\zeta_2>z_1, 
2z_2>z_1+\zeta_1 \},
\eeq 
we have
\bea  \label{RHS-asso-co-co-1}
&&Y_{cl-op}(\mathbb{Y}_{an}(u_1^L\otimes u_1^R; z_1-z_2, 
\zeta_1-\zeta_2)u_2^L\otimes u_2^R; z_2, \zeta_2)v  \nn
&&\hspace{0.5cm}=
\sum_{i=1}^N Y_{cl-op}(\Y_i^L(u_1^L, z_1-z_2)u_2^L \otimes 
\Y_i^R(u_1^R, \zeta_1-\zeta_2)u_2^R; z_2, \zeta_2)v \nn
&&\hspace{0.5cm}=
\sum_{i=1}^N \Y^{(5)}(\Y^{(6)}(\Y_i^L(u_1^L, z_1-z_2)u_2^L, z_2-\zeta_2)
\Y_i^R(u_1^R, \zeta_1-\zeta_2)u_2^R, \zeta_2)v.  \nn
\eea
Combining (\ref{RHS-asso-co-co-1}) with (\ref{mu-cl-Y-z-zeta}), 
(\ref{Y-hat-tilde}) and (\ref{m-co-Y-co}) and the fact that 
module maps are continuous, we obtain the following identity:
\bea  \label{asso-co-co-RHS-2}
&&\hspace{-1cm}\mathbb{Y}_{cl-op}(\mathbb{Y}_{an}(u_1^L\otimes u_1^R; z_1-z_2, 
\zeta_1-\zeta_2)u_2^L\otimes u_2^R; z_2, \zeta_2)v  \nn
&&\hspace{-0.5cm}= 
\overline{\mu_{cl-op}^{P(z_2-\zeta_2)P(\zeta_2)} \circ
(T_{P(z_2-\zeta_2)}(\mu_{cl}^{P(z_1-z_2)P(\zeta_1-\zeta_2)}) 
\boxtimes_{P(\zeta_2)} \id_{V_{op}} })
\big( \nn 
&&\hspace{0.5cm}
((u_1^L\boxtimes_{P(z_1-z_2)} u_2^L)\boxtimes_{P(z_2-\zeta_2)} 
(u_1^R\boxtimes_{P(\zeta_1-\zeta_2)} u_2^R))\boxtimes_{P(\zeta_2)} v\big)  
\eea
for $(z_1,\zeta_1, z_2, \zeta_2)$ in the domain $D_4$.
Moreover, by the universal property of tensor product and 
the nondegeneracy of intertwining operator algebra \cite{nioa},
\begin{equation} \label{cat-RHS-asso}
\mu_{cl-op}^{P(z_2-\zeta_2)P(\zeta_2)} \circ
(T_{P(z_2-\zeta_2)}(\mu_{cl}^{P(z_1-z_2)P(\zeta_1-\zeta_2)}) 
\boxtimes_{P(\zeta_2)} \id_{V_{op}} )
\end{equation}
is the unique morphism in $\text{Hom}(T_{P(z_2-\zeta_2)}(V_{cl}
\boxtimes_{P(z_1-z_2)P(\zeta_1-\zeta_2)} V_{cl})  
\boxtimes_{P(\zeta_2)} V_{op}, V_{op})$ such that 
the identity (\ref{asso-co-co-RHS-2}) holds in $D_4$.

Notice that the domains $D_3$ and $D_4$ are disjoint. 
Now we fix a point $(z_1,\zeta_1,z_2,\zeta_2)$ in the domain 
\beq  \label{domain-LHS-1}
\{ (z_1,\zeta_1, z_2, \zeta_2)| 
z_1>\zeta_1>z_2>\zeta_2>0, 
2\zeta_2>z_1, 2\zeta_1>z_1 + z_2, 2z_2>\zeta_1+\zeta_2 \},
\eeq
which is a subdomain of $D_3$. 
Let $\tilde{z}_2=\zeta_1, \tilde{\zeta}_1=z_2$. 
Then the quadruple 
$(z_1,\tilde{\zeta}_1,\tilde{z}_2, \zeta_2)$ is in the domain $D_4$.

By the analytic properties guaranteed by Theorem \ref{prop-op-cl-V}, 
both sides of associativity (\ref{asso-co-co}) can be uniquely extended to 
the boundary $(z_1,\zeta_1,z_2,\zeta_2)\in 
\R_+^{4}\cap M_{\C^{\times}}^4$. Moreover,
$$
\mathbb{Y}_{cl-op}(u_1^L\otimes u_1^R; z_1, \zeta_1)
\mathbb{Y}_{cl-op}(u_2^L\otimes u_2^R; z_2, \zeta_2)v,
$$
and 
$$
\mathbb{Y}_{cl-op}(\mathbb{Y}_{an}(u_1^L\otimes u_1^R; z_1-\tilde{z}_2,
\tilde{\zeta}_1-\zeta_2)
u_2^L\otimes u_2^R; \tilde{z}_2, \zeta_2)v,
$$
can be obtained from each other by 
analytic continuation along the following path
\beq   \label{path-4}
\begin{picture}(14,2)
\put(3, 0){\resizebox{7cm}{2cm}{\includegraphics{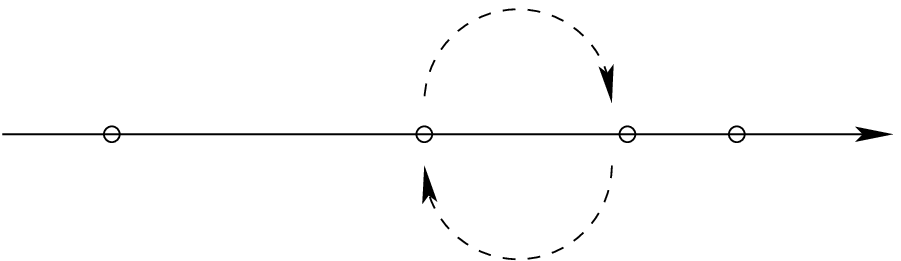}}}
\put(5.8, 1.2){$z_2$}\put(8, 0.5){$\zeta_1$}
\put(5.8, 0.5){$\tilde{\zeta}_1$}\put(8, 1.2){$\tilde{z}_2$}
\put(3.6,0.5){$\zeta_2$}\put(8.8,1.2){$z_1$}
\put(11,0){.}
\end{picture}
\eeq

Meanwhile, if we start from the element
$$
(u_1^L\boxtimes_{P(z_1-\zeta_1)} u_1^R)
\boxtimes_{P(\zeta_1)} 
((u_2^L\boxtimes_{P(z_2-\zeta_2)} u_2^R) \boxtimes_{P(\zeta_2)} v)
$$
in $\overline{T_{P(z_1-\zeta_1)}(V_{cl}) \boxtimes_{P(\zeta_1)} 
(T_{P(z_2-\zeta_2)}(V_{cl}) \boxtimes_{P(\zeta_2)} V_{op})}$ and 
apply associativity isomorphisms repeatedly
and braiding isomorphism on it, we obtain
\bea  \label{AAAR}
&&\xrightarrow{ \overline{
(\mathcal{A}_{P(z_1)P(\zeta_1)}^{P(z_1-\zeta_1)P(\zeta_1)})^{-1} } }
u_1^L\boxtimes_{P(z_1)} (u_1^R \boxtimes_{P(\zeta_1)} 
((u_2^L \boxtimes_{P(z_2-\zeta_2)} u_2^R) \boxtimes_{P(\zeta_2)} v )), \nn
&& \xrightarrow{ \overline{
\mathcal{A}_{P(\zeta_1)P(\zeta_2)}^{P(\zeta_1-\zeta_2)P(\zeta_2)} } }
u_1^L\boxtimes_{P(z_1)}\big( (u_1^R \boxtimes_{P(\zeta_1-\zeta_2)} 
(u_2^L \boxtimes_{P(z_2-\zeta_2)} u_2^R)) \boxtimes_{P(\zeta_2)} v \big), \nn
&&\xrightarrow{ \overline{
\mathcal{A}_{P(\zeta_1-\zeta_2)P(z_2-\zeta_2)}^{P(\zeta_1-z_2)P(z_2-\zeta_2)} } } 
u_1^L\boxtimes_{P(z_1)}\big( ((u_1^R \boxtimes_{P(\zeta_1-z_2)} 
u_2^L) \boxtimes_{P(z_2-\zeta_2)} u_2^R) \boxtimes_{P(\zeta_2)} v \big), \nn
&&\xrightarrow{\overline{\mathcal{R}_-^{P(\zeta_1-z_2)} }} 
u_1^L\boxtimes_{P(z_1)}\big( ( e^{(\zeta_1-z_2)L(-1)} \cdot  \nn
&&\hspace{3cm}
\overline{\mathcal{T}_{\gamma_-}}(u_2^L \boxtimes_{P(-\zeta_1+z_2)} 
u_1^R) \boxtimes_{P(z_2-\zeta_2)} u_2^R) \boxtimes_{P(\zeta_2)} v \big),
\eea
where we have ignored the obvious identity maps and $\gamma_-$ is a path
from $-\zeta_1+z_2$ to $\zeta_1-z_2$ given by 
$\gamma_-(t)= (-\zeta_1+z_2) e^{-\pi it}, t\in [0,1]$.
Now we analytic continuate the last line of (\ref{AAAR})
along the path (\ref{path-4}), we obtain
\beq  \label{after-an-c}
u_1^L\boxtimes_{P(z_1)}\big( ((u_2^L \boxtimes_{P(\zeta_1-z_2)} 
u_1^R) \boxtimes_{P(z_2-\zeta_2)} u_2^R) \boxtimes_{P(\zeta_2)} v \big), 
\eeq
which is nothing but the third line of (\ref{AAAR}) except
$u_1^R$ and $u_2^L$ being exchanged. 
Note that the last line of (\ref{AAAR}) and (\ref{after-an-c}) are 
elements in the algebraic completion of the same $V$-module. 
We denote this $V$-module as $W$. Now we further apply 
associativity morphisms on (\ref{after-an-c}) and 
use $z_2=\tilde{\zeta}_1$ and $\zeta_1=\tilde{z}_2$. 
We then obtain
\bea  \label{after-an-c-AAA}
&&\hspace{-1cm}\xrightarrow{\overline{(
\mathcal{A}_{P(\tilde{z}_2-\zeta_2)}
^{P(\tilde{z}_2-\tilde{\zeta}_2)P(\tilde{\zeta}_1-\zeta_2)} )^{-1}}}
u_1^L\boxtimes_{P(z_1)}\big( (u_2^L \boxtimes_{P(\tilde{z}_2-\zeta_2)} 
(u_1^R \boxtimes_{P(\tilde{\zeta}_1-\zeta_2)} u_2^R) )
\boxtimes_{P(\zeta_2)} v \big) , \nn
&&\hspace{-1cm}\xrightarrow{
\overline{\mathcal{A}_{P(z_1)P(\zeta_2)}^{P(z_1-\zeta_2)P(\zeta_2)}} }
\big( u_1^L\boxtimes_{P(z_1-\zeta_2)} (u_2^L \boxtimes_{P(\tilde{z}_2-\zeta_2)} 
(u_1^R \boxtimes_{P(\tilde{\zeta}_1-\zeta_2)} u_2^R) )\big) 
\boxtimes_{P(\zeta_2)} v  ,   \nn
&&\hspace{-1cm}\xrightarrow{ \overline{ 
\mathcal{A}_{P(z_1-\zeta_2)P(\tilde{z}_2-\zeta_2)}
^{P(z_1-\tilde{z}_2)P(\tilde{z}_2-\zeta_2)} } }
\big( (u_1^L\boxtimes_{P(z_1-\tilde{z}_2)} u_2^L) 
\boxtimes_{P(\tilde{z}_2-\zeta_2)} 
(u_1^R \boxtimes_{P(\tilde{\zeta}_1-\zeta_2)} u_2^R) \big)  
\boxtimes_{P(\zeta_2)} v.
\eea

Let $m$ be the morphism $W\rightarrow V_{op}$ such that 
\bea  \label{equ-m-L}
&&\hspace{-1cm}\mbox{ (\ref{cat-LHS-asso}) } = 
m\circ \mathcal{R}_-^{P(\zeta_1-z_2)} \circ 
\mathcal{A}_{P(\zeta_1-\zeta_2)P(z_2-\zeta_2)}^{P(\zeta_1-z_2)P(z_2-\zeta_2)} 
\circ   \nn
&&\hspace{4cm} 
\mathcal{A}_{P(\zeta_1)P(\zeta_2)}^{P(\zeta_1-\zeta_2)P(\zeta_2)}
\circ
(\mathcal{A}_{P(z_1)P(\zeta_1)}^{P(z_1-\zeta_1)P(\zeta_1)})^{-1}. 
\eea
If we apply $\overline{m}$ on both the 
last line of (\ref{AAAR}) and (\ref{after-an-c}), the two 
images of $\overline{m}$ are clearly the analytic
continuation of each other along the path (\ref{path-4}). 
On the other hand, combining this fact with 
(\ref{asso-co-co-RHS-2}) and (\ref{after-an-c-AAA}), 
we obtain that 
\bea \label{equ-m-R}
m = \mbox{(\ref{cat-RHS-asso})} \circ 
\mathcal{A}_{P(z_1-\zeta_2)P(\tilde{z}_2-\zeta_2)}
^{P(z_1-\tilde{z}_2)P(\tilde{z}_2-\zeta_2)}  \circ
\mathcal{A}_{P(z_1)P(\zeta_2)}^{P(z_1-\zeta_2)P(\zeta_2)} \circ
(\mathcal{A}_{P(\tilde{z}_2-\zeta_2)}
^{P(\tilde{z}_2-\tilde{\zeta}_2)P(\tilde{\zeta}_1-\zeta_2)} )^{-1}
\eea
because the extensions
of both sides of (\ref{equ-m-R}), applied on (\ref{after-an-c}),
give the same element in $\overline{V_{op}}$. 
Therefore, we further obtain 
from (\ref{equ-m-L}) and (\ref{equ-m-R}) the following identity: 
\bea  \label{equ-pre-diag}
&&\hspace{-1cm}\mbox{ (\ref{cat-LHS-asso}) } =
\mbox{(\ref{cat-RHS-asso})} \circ 
\mathcal{A}_{P(z_1-\zeta_2)P(\tilde{z}_2-\zeta_2)}
^{P(z_1-\tilde{z}_2)P(\tilde{z}_2-\zeta_2)}  \circ
\mathcal{A}_{P(z_1)P(\zeta_2)}^{P(z_1-\zeta_2)P(\zeta_2)} \circ
(\mathcal{A}_{P(\tilde{z}_2-\zeta_2)}
^{P(\tilde{z}_2-\tilde{\zeta}_2)P(\tilde{\zeta}_1-\zeta_2)} )^{-1} \nn
&&\hspace{0.7cm} 
\circ \mathcal{R}_-^{P(\zeta_1-z_2)} \circ 
\mathcal{A}_{P(\zeta_1-\zeta_2)P(z_2-\zeta_2)}^{P(\zeta_1-z_2)P(z_2-\zeta_2)} 
\circ   
\mathcal{A}_{P(\zeta_1)P(\zeta_2)}^{P(\zeta_1-\zeta_2)P(\zeta_2)}
\circ
(\mathcal{A}_{P(z_1)P(\zeta_1)}^{P(z_1-\zeta_1)P(\zeta_1)})^{-1}. \nn
\eea

Using the commutative diagram (\ref{A-T-A-T})(\ref{R-T-R-T})
and the definition of $\varphi_2$ (recall (\ref{varphi-2-defn})), 
it is easy to see that (\ref{equ-pre-diag}) implies the 
commutativity of the following diagram:
\beq  \label{cd-1}
\xymatrix{
T_{P(z_1-\zeta_1)}(V_{cl}) \boxtimes_{P(\zeta_1)} 
(T_{P(z_2-\zeta_2)}(V_{cl}) \boxtimes_{P(\zeta_2)} V_{op} )  
\ar[d]_{f_1} \ar[rdd]^{(\ref{cat-LHS-asso})}  &  \\
T(V_{cl})\boxtimes (T(V_{cl}) \boxtimes V_{op}) \ar[d]_{\mathcal{A}}  
\ar@{.>}[rd] &  \\
(T(V_{cl})\boxtimes T(V_{cl}) )\boxtimes V_{op} 
\ar[d]_{\varphi_2\boxtimes \id_{V_{op}} }     &  V_{op} \\
T(V_{cl}\boxtimes V_{cl}) \boxtimes V_{op}
\ar[d]_{g_1} \ar@{.>}[ru]  &   \\
T_{P(\tilde{z}_2-\zeta_2)}(V_{cl}
\boxtimes_{P(z_1-\tilde{z}_2)P(\tilde{\zeta}_1-\zeta_2)} V_{cl})  
\boxtimes_{P(\zeta_2)} V_{op}  
\ar[ruu]_{\hspace{0.4cm}(\ref{cat-RHS-asso})|_{z_2=\tilde{z}_2, 
\zeta_1=\tilde{\zeta}_1}} 
}
\eeq
where 
\bea
&&f_1= \id_{T(V_{cl})} \boxtimes (\mathcal{T}_{\gamma_4} \boxtimes 
\id_{H_{op}} ) 
\circ \mathcal{T}_{\gamma_2} \boxtimes \mathcal{T}_{\gamma_3} 
\circ \mathcal{T}_{\gamma_1} \\
&&g_1= T_{P(z_2-\zeta_2)} (\mathcal{T}_{\gamma_7} \otimes 
\mathcal{T}_{\gamma_8}) 
\boxtimes_{P(\zeta_2)} \id_{H_{op}} \circ 
(\mathcal{T}_{\gamma_6} \boxtimes_{P(\zeta_2)} \id_{V_{op}}) 
\circ \mathcal{T}_{\gamma_5}
\eea
in which $\gamma_i, i=1,\ldots, 4$ are paths in $\R_+$ from 
$\zeta_1, z_1-\zeta_1, \zeta_2, z_2-\zeta_2$ to $1$ respectively and 
$\gamma_i, i=5,\ldots, 8$ are paths in $\R_+$ 
from $1$ to $\zeta_2,z_2-\zeta_2, z_1-z_2, \zeta_1-\zeta_2$
respectively. 

Using (\ref{m-co-def}), it is easy to see that 
\bea \label{asso-parallel}
\mu_{cl-op} \circ (\id_{T(V_{cl})} \boxtimes 
\mu_{cl-op}) &=& (\ref{cat-LHS-asso}) \circ f_1^{-1}, \nn
\mu_{cl-op} \circ (T(\mu_{cl}) \boxtimes \id_{V_{op}})
&=&  (\ref{cat-RHS-asso}) \circ g_1.
\eea
(\ref{asso-parallel}) together with the commutative diagram 
(\ref{cd-1}) implies (\ref{asso-co-co-cat}), which is 
nothing but the commutativity of the subdiagram in 
the middle of (\ref{cd-1}).

Conversely, (\ref{asso-co-co-cat}) implies 
the commutativity of the diagram (\ref{cd-1}). 
It is easy to see that above arguments can be reversed. 
Therefore, we can also obtain the associativity II
(\ref{asso-co-co}) from (\ref{asso-co-co-cat}). 
\epf

We now study the categorical formulations of 
the rest conditions needed in Theorem \ref{prop-op-cl-V}. 
The proof of them are essentially the same 
as that of Proposition \ref{cat-lemma-asso-co-co}. So 
we will only sketch the proofs below. 

\begin{prop}
The associativity I (recall Proposition \ref{prop-asso-co-op})
is equivalent to the following condition: 
\beq \label{asso-co-op-cat}
\mu_{cl-op} (\id_{T(V_{cl})} \boxtimes \mu_{op}) = \mu_{op}(\mu_{cl-op} \boxtimes 
\id_{V_{op}}) \circ \mathcal{A}
\eeq
which can also be expressed in the following graph:
\beq  \label{asso-co-op-fig}
\epsfxsize  0.5\textwidth
\epsfysize  0.2\textwidth
\epsfbox{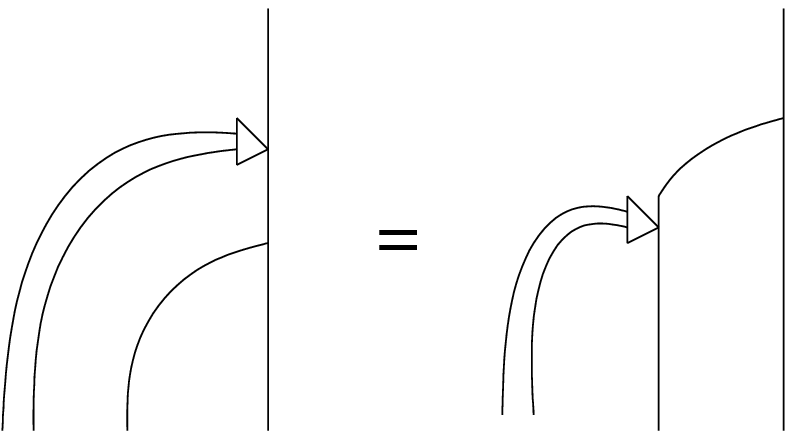}\, \, .
\eeq
\end{prop}
\pf
The left hand side of (\ref{asso-co-op}) gives arise to a morphism 
\beq \label{asso-co-op-LHS}
\mu_{cl-op}^{P(z-\zeta)P(\zeta)} \circ (\id_{T_{P(z-\zeta)}(V_{cl})} 
\boxtimes_{P(\zeta)} \mu_{op}^{P(r)} )
\eeq
in $\hom(T_{P(z-\zeta)}(V_{cl}) \boxtimes_{P(\zeta)} (V_{op} \boxtimes_{P(r)}
V_{op}), V_{op})$. The right hand side of 
(\ref{asso-co-op}) gives arise to a morphism 
\beq \label{asso-co-op-RHS}
\mu_{op}^{P(r)} \circ (\mu_{cl-op}^{P(z-\zeta)P(\zeta-r)}
\boxtimes_{P(\zeta-r)} \id_{V_{op}} )
\eeq
in $\hom( (T_{P(z-\zeta)}(V_{cl}) \boxtimes_{P(\zeta)} V_{op})
\boxtimes_{P(r)} V_{op},V_{op})$. 

For $z, \zeta, r$ in a proper
subdomain of $z>\zeta>r>0$,
using similar arguments as in Proposition 
\ref{cat-lemma-asso-co-co}, we obtain that  
the associativity (\ref{asso-co-op}), for some 
$z, \zeta, r, \tilde{r}\in \R_+$, implies the following 
commutative diagram:
\beq \label{cd-2}
\xymatrix{ 
T_{P(z-\zeta)}(V_{cl}) \boxtimes_{P(\zeta)} (V_{op} \boxtimes_{P(r)}
V_{op}) \ar[d]_{f_2}  \ar[rd]^{(\ref{asso-co-op-LHS})} & \\
T(V_{cl}) \boxtimes (V_{op} \boxtimes V_{op}) 
\ar[d]_{\mathcal{A}} \ar@{.>}[r]  & V_{op}    \\
(T(V_{cl}) \boxtimes V_{op}) \boxtimes V_{op}
\ar[d]_{g_2} \ar@{.>}[ru]   &  \\
(T_{P(z-\zeta)}(V_{cl}) \boxtimes_{P(\zeta-r)} V_{op})
\boxtimes_{P(r)} V_{op})\, , \ar[ruu]_{(\ref{asso-co-op-RHS})} 
}
\eeq
where 
\bea
&& f_2= (\mathcal{T}_{\gamma_2} \boxtimes \mathcal{T}_{\gamma_3} )
\circ \mathcal{T}_{\gamma_1},  \\
&& g_2= \mathcal{T}_{\gamma_6} \boxtimes_{P(\zeta)} \id_{V_{op}})
\boxtimes_{P(r)} \id_{V_{op}}  \circ 
(\mathcal{T}_{\gamma_5} \boxtimes_{P(r)} \id_{V_{op}}) 
\circ \mathcal{T}_{\gamma_4},
\eea
where $\gamma_i, i=1,2,3$ are paths in $\R_+$ from 
$\zeta, z-\zeta$ and $r$ to $1$ respectively and 
$\gamma_i, i=4,5,6$ are paths in $\R_+$ from 
$1$ to $r, \zeta-r$ and $z-\zeta$ respectively. 
The commutativity of outside loop in (\ref{cd-2})
implies immediately the commutativity of the 
subdiagram in the middle of (\ref{cd-2}), which is 
nothing but the identity
(\ref{asso-co-op-cat}) or (\ref{asso-co-op-fig}).

Conversely, using (\ref{cd-2}) and reversing above arguments, 
it is clear that (\ref{asso-co-op-cat}) or (\ref{asso-co-op-fig})
also implies the associativity (\ref{asso-co-op}). 
\epf

Let $V_{cl}=\oplus_{i=1}^N W_i^L\otimes W_i^R$. We define a map 
$T(V_{cl})\boxtimes V_{op} \xrightarrow{\sigma_1} V_{op} \boxtimes T(V_{cl})$ by
\beq
\sigma_1 := \oplus_{i=1}^N (\mathcal{R}_+ \boxtimes \id_{W_i^R})
\circ \mathcal{A} \circ (\id_{W_i^L} \boxtimes \mathcal{R}_-)
\circ  \mathcal{A}^{-1}. 
\eeq

\begin{prop}
The commutativity I (recall Proposition \ref{comm-prop-1}) 
is equivalent to the following identity: 
\beq \label{comm-co-op-cat}
\mu_{cl-op}(\id_{T(V_{cl})} \boxtimes \mu_{op}) = \mu_{op}(\id_{V_{op}} 
\boxtimes \mu_{cl-op}) \circ 
\mathcal{A}^{-1}\circ \sigma_1 \circ \mathcal{A}, 
\eeq
or the following graphic identities:
\beq \label{comm-co-op-fig}
\epsfxsize  0.8\textwidth
\epsfysize  0.2\textwidth
\epsfbox{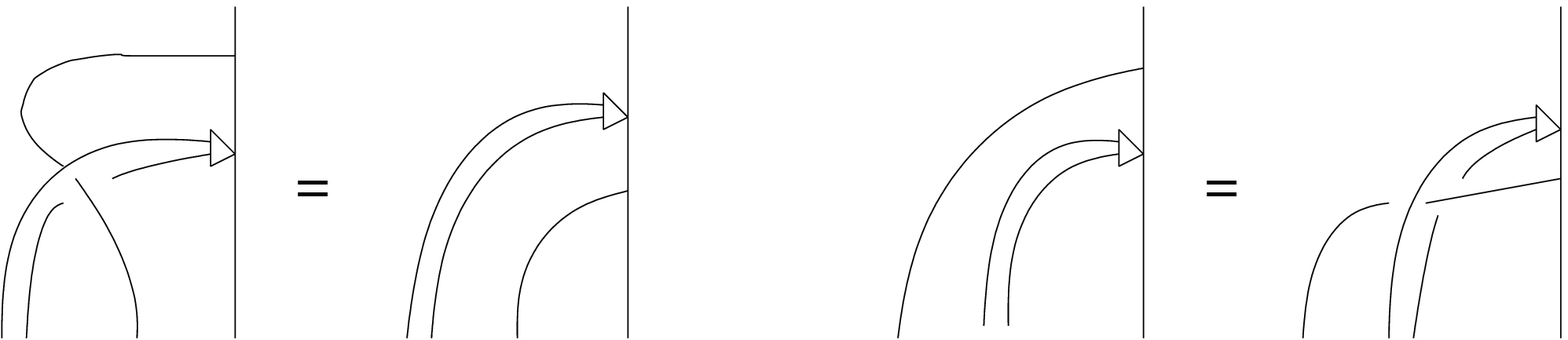}\,\, .
\eeq
\end{prop}
\pf
There is a morphism
\beq \label{comm-co-op-L}
\mu_{cl-op}^{P(z-\zeta)P(\zeta)} \circ 
(\id_{T_{P(z-\zeta)}(V_{cl})} \boxtimes_{P(\zeta)} \mu_{op}^{P(r)})
\eeq
in $\hom(T_{P(z-\zeta)}(V_{cl})\boxtimes_{P(\zeta)} 
(V_{op}\boxtimes_{P(r)} V_{op}))$
associated to (\ref{comm-ocfa-coop-1}). 
There is another morphism
\beq \label{comm-co-op-R}
\mu_{op}^{P(\tilde{r})} \circ (\id_{V_{op}} \boxtimes_{P(\tilde{r})} \mu_{cl-op}^{P(z-\zeta)P(\zeta)})
\eeq
in $\hom(V_{op}\boxtimes (T(V_{cl})\boxtimes V_{op}))$
associated to (\ref{comm-ocfa-coop-2}).

For $z, \zeta, r, r_1$ in a proper
subdomain of $r_1>z>\zeta>r>0$,
using similar arguments as in Proposition 
\ref{cat-lemma-asso-co-co}, 
we obtain that the commutativity I, 
implies the following commutative diagram:
\beq  \label{cd-4}
\xymatrix{
T_{P(z-\zeta)}(V_{cl})\boxtimes_{P(\zeta)} 
(V_{op}\boxtimes_{P(r)} V_{op}) \ar[d]_{f_4}
\ar[rdd]^{(\ref{comm-co-op-L})} &  \\
T(V_{cl}) \boxtimes (V_{op}\boxtimes V_{op})
\ar[d]_{\mathcal{A}} \ar@{.>}[rd]   &  \\
(T(V_{cl}) \boxtimes V_{op}) \boxtimes V_{op}
\ar[d]_{\sigma_1} &  V_{op} \ar[d]^{\id_{V_{op}}} \\
(V_{op} \boxtimes T(V_{cl})) \boxtimes V_{op}
\ar[d]_{\mathcal{A}^{-1}}   &   V_{op} \\
V_{op} \boxtimes ( T(V_{cl}) \boxtimes V_{op}) 
\ar[d]_{g_4} \ar@{.>}[ru]  &  \\
V_{op}\boxtimes_{P(r_1)} (T_{P(z-\zeta)}(V_{cl})\boxtimes_{P(\zeta)} V_{op})
\ar[ruu]_{(\ref{comm-co-op-R})} 
}
\eeq
where 
\bea
f_4 &=& (\mathcal{T}_{\gamma_2} \boxtimes \mathcal{T}_{\gamma_3}) \circ
\mathcal{T}_{\gamma_1},  \\
g_4 &=& (\id_{V_{op}} \boxtimes (\mathcal{T}_{\gamma_5} \boxtimes 
\mathcal{T}_{\gamma_6}) \circ \mathcal{T}_{\gamma_4},
\eea
in which $\gamma_i, i=1,2,3$ are paths in $\R_+$ from 
$\zeta$, $z-\zeta$ and $r$ respectively to $1$ and 
$\gamma_i, i=4,5,6$ are path in $\R_+$ from 
$1$ to $r_1$, $z-\zeta$ and $\zeta$ respectively.
Above commutative diagram immediately implies that 
the subdiagram in the middle of (\ref{cd-4}) is commutative. 
This is nothing but the commutativity
(\ref{comm-co-op-cat}) or the first formula in 
(\ref{comm-co-op-fig}). 
Moreover, it also easy to see that 
the two formula in (\ref{comm-co-op-fig}) are actually equivalent. 

Conversely, 
using commutative diagram (\ref{cd-4}) and 
reversing above arguments, it is clear that 
(\ref{comm-co-op-fig}) implies the commutativity of 
rational $\tilde{\mathfrak{S}}^c$.
\epf

Commutativity II (recall Proposition \ref{comm-prop-2})
is not needed in Theorem \ref{prop-op-cl-V} because it
automatically follows from associativity II and skew symmetry
of $V_{cl}$. It also has a very nice categorical formulation as
given in the following proposition, which follows from
(\ref{T-cl-2-fig}) and (\ref{asso-co-co-fig}) immediately.

\begin{prop}
For an open-closed field algebra over $V$, we have
\beq \label{comm-co-co-cat}
\mu_{cl-op} \circ \id_{T(V_{cl})} \boxtimes \mu_{cl-op}  \nn
= \mu_{cl-op}\circ \id_{T(V_{cl})} \boxtimes \mu_{cl-op} \circ 
\mathcal{A}^{-1} \circ \sigma \circ \mathcal{A},
\eeq
or equivalently
\beq \label{comm-co-co-fig}
\epsfxsize  0.7\textwidth
\epsfysize  0.2\textwidth
\epsfbox{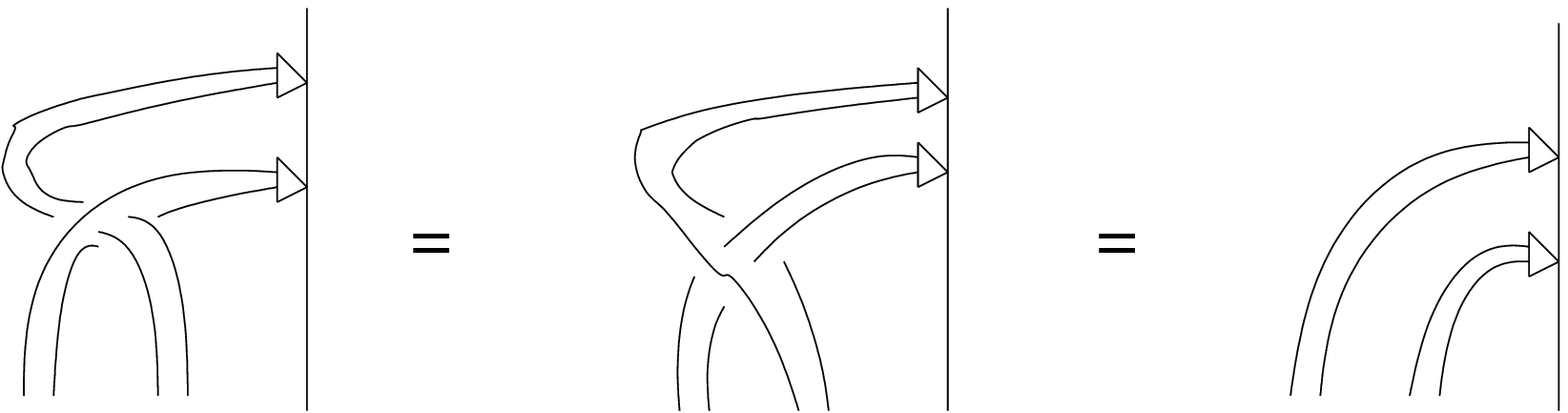} \,\,.
\eeq
\end{prop}

In summary, we have already completely reformulated the 
all the data and conditions in Theorem \ref{prop-op-cl-V}
in the language of tensor category 
as (\ref{m-co-def}), (\ref{id-cat-def}), (\ref{asso-co-co-cat}),
(\ref{asso-co-op-cat}) 
and (\ref{comm-co-op-cat}) 
or equivalently as graphic identities
(\ref{id-fig}), (\ref{asso-co-co-fig}), (\ref{asso-co-op-fig})
and (\ref{comm-co-op-fig}).

We define a morphism $\iota_{cl-op}: T(V_{cl}) \rightarrow V_{op}$ as
the composition of the following maps: 
\beq
T(V_{cl}) \xrightarrow{r_{T(V_{cl})}^{-1}} T(V_{cl}) \boxtimes 
\one_{\mathcal{C}_V} \xrightarrow{\id \boxtimes \, \iota_{V_{op}}}
 T(V_{cl}) \boxtimes V_{op} \xrightarrow{\mu_{cl-op}} V_{op}. 
\eeq
or equivalently as the following graphic formula:
\beq   
\begin{picture}(14,1.5)
\put(3,0.6){$\iota_{cl-op}$}\put(4.3,0.6){$=$}
\put(5, 0){\resizebox{4cm}{1.5cm}
{\includegraphics{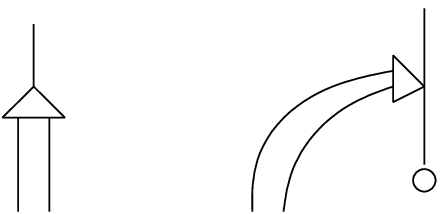}}}
\put(6.2, 0.6){$:=$}\put(9.5,0){.}
\label{iota-*-cat-def}
\end{picture}
\eeq

\begin{lemma} 
$\iota_{cl-op}$ is an algebra morphism from $T(V_{cl})$ to $V_{op}$. 
\end{lemma}
\pf
That $\iota_{cl-op}$ maps identity to identity is
proved as follows:
\beq
\label{iota-unit-fig}
\epsfxsize  0.3\textwidth
\epsfysize  0.12\textwidth
\epsfbox{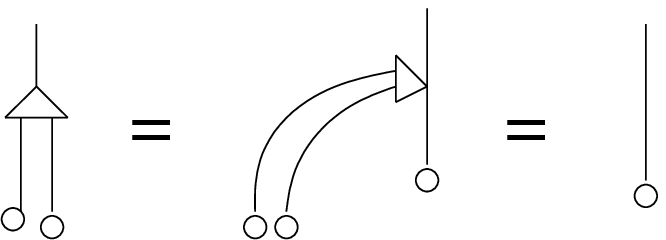}\, \, .
\eeq
The homomorphism property 
$\mu_{op}\circ (\iota_{cl-op} \boxtimes \iota_{cl-op}) =
\iota_{cl-op} \circ \tilde{m}_{cl}$ is proven as follows:
\beq
\label{iota-m-fig}
\epsfxsize  0.8\textwidth
\epsfysize  0.15\textwidth
\epsfbox{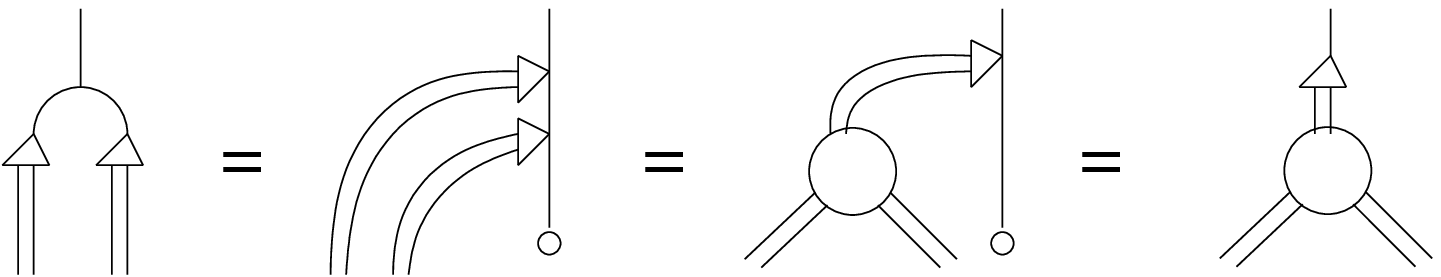} \, \, .
\eeq
\epf

\begin{defn} \label{w-cl-op-alg}  {\rm
An {\it open-closed $\mathcal{C}_V|\mathcal{C}_{V\otimes V}$-algebra}, 
denoted as 
$$
( \, (A_{op}, \mu_{op}, \iota_{op})|(A_{cl},\mu_{cl}, \iota_{cl}),  
\iota_{cl-op}\, )
$$  
or simply $(A_{op}|A_{cl}, \iota_{cl-op})$,  consists of
an associative algebra $(A_{op}, \mu_{op}, \iota_{op})$ in $\mathcal{C}_V$,
a commutative associative algebra with a trivial twist 
$(A_{cl}, \mu_{cl}, \iota_{cl})$ 
in $\mathcal{C}_{V\otimes V}$ and an associative algebra 
homomorphism $\iota_{cl-op}: T(A_{cl}) \rightarrow A_{op}$,
satisfying the following commutativity:
\beq
\mu_{op} \circ (\iota_{cl-op} \boxtimes \id_{V_{op}})     
= \mu_{op} \circ 
(  \id_{V_{op}} \boxtimes \iota_{cl-op}) \circ \sigma_1, 
\eeq
or equivalently,
\beq \label{iota-m-comm-fig}
\epsfxsize  0.4\textwidth
\epsfysize  0.16\textwidth
\epsfbox{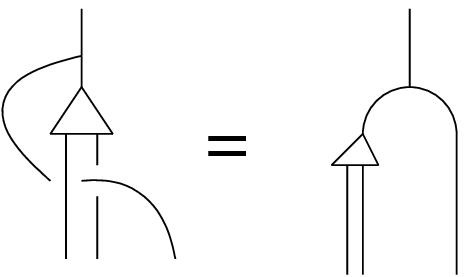}\,\, .
\eeq
}
\end{defn}

\begin{thm} \label{s-alg-cat}
The notion of open-closed field algebra over $V$ is equivalent to that of 
open-closed $\mathcal{C}_V|\mathcal{C}_{V\otimes V}$-algebra,
in the sense that the categories of above two notions are isomorphic. 
\end{thm}
\pf
Given an open-closed field algebra
over $V$, we have shown that it gives 
a triple $(V_{cl}, V_{op},\mu_{cl-op})$, in which
$V_{cl}$ is a commutative associative algebra in 
$\mathcal{C}_{V\otimes V}$ with a trivial twist,   
and $V_{op}$ is an algebra in $\mathcal{C}_V$, 
and $\mu_{cl-op}$ satisfies
(\ref{id-fig}),(\ref{asso-co-co-fig}),
(\ref{asso-co-op-fig}) and (\ref{comm-co-op-fig}). 
Moreover, we have shown that $\iota_{cl-op}$ defined 
by (\ref{iota-*-cat-def}) gives an 
morphism of associative algebra. 
Now we prove (\ref{iota-m-comm-fig}) as follows: 
$$
\label{iota-m-comm-1}
\epsfxsize  0.9\textwidth
\epsfysize  0.16\textwidth
\epsfbox{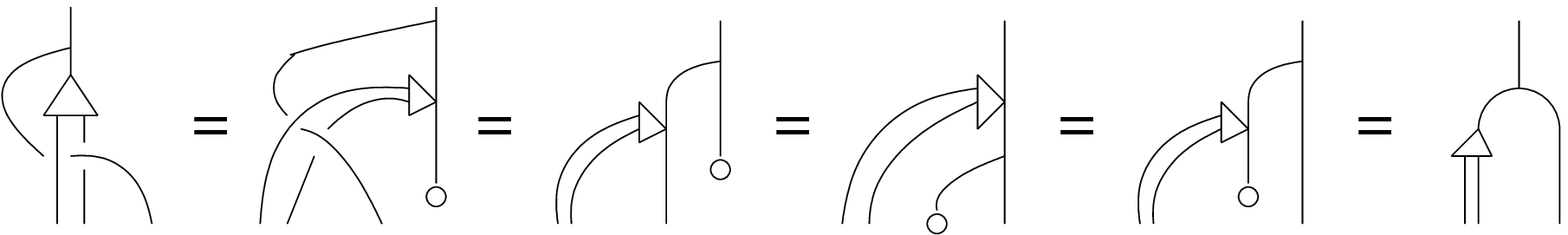}\,\, .
$$
Hence $(V_{op}|V_{cl}, \iota_{cl-op})$ is an open-closed
$\mathcal{C}_V|\mathcal{C}_{V\otimes V}$-algebra. 

It is easy to check that it gives a functor from the category 
of open-closed field algebras over $V$ to that of open-closed 
$\mathcal{C}_V|\mathcal{C}_{V\otimes V}$-algebras.

Conversely, given an open-closed 
$\mathcal{C}_V|\mathcal{C}_{V\otimes V}$-algebra, 
$(V_{op}|V_{cl}, \iota_{cl-op})$, we define a morphism 
$\mu_{cl-op} \in \hom(T(V_{cl})\boxtimes V_{op}, V_{op})$ as 
\beq   
\begin{picture}(14,2)
\put(3,0.8){$\mu_{cl-op}$}\put(4.3,0.8){$=$}
\put(5, 0){\resizebox{5cm}{2cm}
{\includegraphics{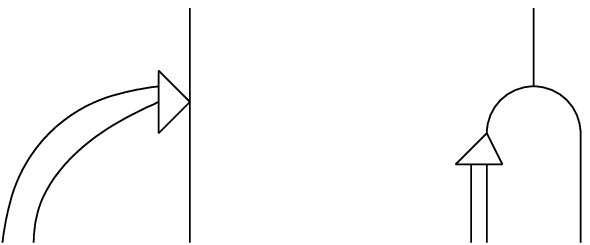}}}
\put(7.8, 0.8){$:=$}\put(10.5,0){.}
\label{m-co-by-iota}
\end{picture}
\eeq
Since $\iota_{cl-op}$ is an algebra homomorphism, 
it maps unit to unit (recall
(\ref{iota-unit-fig})). Thus the identity property of 
$\mu_{cl-op}$ holds. Then the identity property of $\mathbb{Y}_{cl-op}$ 
follows.

Furthermore, we have
$$
\epsfxsize  0.8\textwidth
\epsfysize  0.16\textwidth
\epsfbox{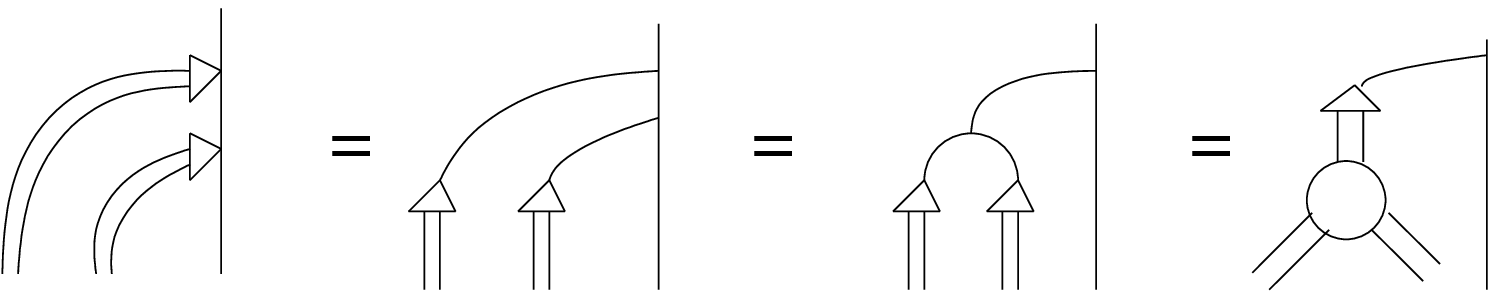}\,\, ,
$$
which gives the associativity II (\ref{asso-co-co-fig}). 
The associativity I (\ref{asso-co-op-fig}) follows from 
$$
\epsfxsize  0.8\textwidth
\epsfysize  0.16\textwidth
\epsfbox{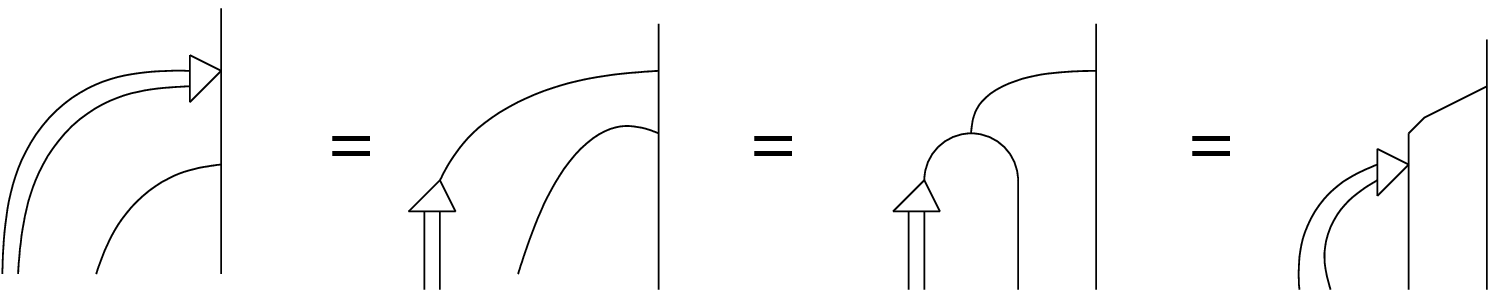}\, \, .
$$
By using (\ref{iota-m-comm-1}), 
we can prove the commutativity (\ref{comm-co-op-fig}) as follows:
$$
\epsfxsize  0.8\textwidth
\epsfysize  0.16\textwidth
\epsfbox{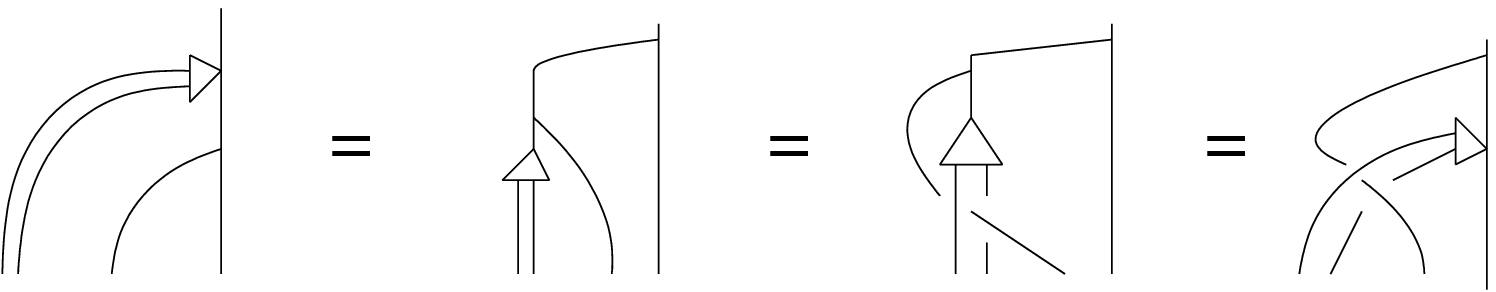}\,\, .
$$
Notice that the other half of (\ref{comm-co-op-fig}) is equivalent to 
the first half. Also notice that 
the commutativity (\ref{comm-co-co-fig}) simply follows from 
(\ref{comm-T-cl}) and the fact that $\iota_{cl-op}$ 
is an algebra homomorphism. 

It is easy to check that this correspondence 
gives a functor from the category of 
open-closed $\mathcal{C}_V|\mathcal{C}_{V\otimes V}$-algebras
to that of open-closed field algebras over $V$. 

That the two relevent categories are isomorphic follows from 
(\ref{iota-*-cat-def}) and (\ref{m-co-by-iota}) easily. 
\epf


It is very easy to construct  
open-closed $\mathcal{C}_V|\mathcal{C}_{V\otimes V}$-algebras.
For example, let $A$ be an associative algebra in $\mathcal{C}_{V\otimes V}$
and $C_l(A)$ the left center of $A$ \cite{O}. 
Let $A_0$ be any subalgebra of $C_l(A)$ and 
$\iota: A_0\hookrightarrow A$ the natural embedding.
Then it is clear that $(T(A)|A_0, T(\iota))$
gives an open-closed $\mathcal{C}_V|\mathcal{C}_{V\otimes V}$-algebra, 
which further gives an open-closed field algebra over $V$ and 
a smooth $\tilde{\mathfrak{S}}^c$-algebra over $V$. 
We will not pursue the construction further in this work. 
Instead, we leave it to \cite{Ko2}, 
in which we will construct explicit examples 
that are most relevent in physics. Moreover, a general theory of 
open-closed $\mathcal{C}_V|\mathcal{C}_{V\otimes V}$-algebras 
and its variants (\cite{Ko2}) will be
developed and its connection to the works
\cite{FFFS}\cite{FS3}\cite{FRS1}-\cite{FRS4}\cite{FjFRS} 
will be explained in future publications.


\noindent {\small \sc Max Planck Institute for Mathematics
in the Sciences, Inselstrasse 22, D-04103, Leipzig, Germany} \\
\noindent {and}   \\
\noindent {\small \sc 
Institut Des Hautes \'{E}tudes Scientifiques, 
Le Bois-Marie, 35, Route De Chartres,
F-91440 Bures-sur-Yvette, France} \\
\noindent {and}   \\
\noindent {\small \sc 
Max-Planck-Institut f\"{u}r Mathematik
Vivatsgasse 7, D-53111 Bonn, Germany} (current address)\\
\noindent {\em E-mail address}: kong@mpim-bonn.mpg.de

\end{document}